\documentclass[A4j,11pt]{article}
\usepackage{graphics}
\usepackage{amsfonts,amsmath,amscd}
\usepackage{amssymb}
\usepackage[notref,notcite]{showkeys}
\usepackage[dvipdfmx]{graphicx}
\usepackage{tikz}
\usepackage{amssymb,amsmath,amsthm,amscd,amsfonts}
\begin{document}

\title{{\bf Homogeneity of infinite dimensional\\
anti-Kaehler isoparametric submanifolds II}}
\author{{\bf Naoyuki Koike}}

\date{}

\maketitle

{\small\textit{Department of Mathematics, Faculty of Science, 
Tokyo University of Science,}}

{\small\textit{1-3 Kagurazaka Shinjuku-ku, Tokyo 162-8601 Japan}}

{\small\textit{E-mail address}: koike@rs.kagu.tus.ac.jp}

%
%
%

\begin{abstract}
In this paper, we prove that, if a full irreducible infinite dimensional anti-Kaehler 
isoparametric submanifold of codimension greater than one has $J$-diagonalizable 
shape operators, then it is an orbit of the action of a Banach Lie group generated by 
one-parameter transformation groups induced by holomorphic Killing vector fields defined entirely 
on the ambient Hilbert space.  
\end{abstract}

\vspace{0.1truecm}


\hspace{1.95truecm}{\small distribution, homogeneous structure, holomorphic Killing vector field}



\vspace{0.3truecm}

\section{Introduction} 
An {\it infinite dimensional isoparametric submanifold} is a proper Fredholm submanifold of finite codimension 
in an infinite dimensional separable Hilbert space over the real number field ${\Bbb R}$ such that its normal 
holonomy group is trivial and that the shape operator for each parallel normal vector field has constant 
eigenvalues, where ``proper Fredholm'' means that the differential of the normal exponential map 
$\exp^{\perp}$ of the submanifold is a Fredholm operator and that the restriction of $\exp^{\perp}$ to unit 
ball normal bundle is proper.  
Throughout this paper, all Hilbert spaces mean infinite dimensional separable Hilbert spaces.  
In 1999, E. Heintze and X. Liu ([HL2]) proved that all full irreducible infinite dimensional isoparametric 
submanifolds of codimension greater than one in a Hilbert space are 
extrinsically homogeneous.  In 2002, by using this result of Heintze-Liu, 
U. Christ ([Ch]) claimed that all irreducible equifocal submanifolds with flat 
section of codimension greater than one in a simply connected symmetric space of compact 
type are extrinsically homogeneous.  
Let $I(V)$ be the group of all isometries of the Hilbert space $V$ and 
$M$ a full irreducible isoparametric submanifolds of codimension greater than one in $V$.  
Set $H:=\{F\in I(V)\,\vert\,F(M)=M\}$.  
The extrinsic homogeneity of $M$ in the result of [HL2] means that $Hx=M$ ($x\in M$).  
Let $I_b(V)$ be the subgroup of $I(V)$ generated by one-parameter transformation groups induced by 
the Killing vector fields defined entirely on $V$.  Note that $I_b(V)$ is a Banach Lie group.  
Set $H_b:=H\cap I_b(V)$, which is a Banach Lie subgroup of $I(V)$.  
Recently, C. Gorodski and E. Heintze ([GH]) proved that $H_bx=M$ holds for any $x\in M$.  
This improved extrinsic homogeneity theorem closed a gap in the proof of the above extrinsic homogeneity 
theorem by U. Christ.  

In [K1], we introduced the notion of a complex equifocal submanifold in a symmetric 
space of non-compact type.  
In [K2], we showed that the study of complex equifocal $C^{\omega}$-submanifolds in symmetric spaces of 
non-compact type is converted to that of anti-Kaehler 
isoparametric submanifolds in the infinite dimensional anti-Kaehler space, where 
$C^{\omega}$ means the real analyticity.  In this paper, we shall investigate an anti-Kaehler isoparametric 
submanifold with $J$-diagonalizable shape opeartors, which was called a {\it proper} 
anti-Kaehler isoparametric submanifold in [K2].  
L. Geatti and C. Gorodski ([GG]) introduced the notion of an isoparametric submanifold 
with diagonalizable Weingarten operators in a finite dimensional pseudo-Euclidean space.  
Note that anti-Kaehler isoparametric submanifolds with 
$J$-diagonalizable shape operators give a subclass of the infinite 
dimensional version of isoparametric submanifolds with diagonalizable Weingarten operators.  
Let $K$ be a maximal compact subgroup of a finite dimensional non-compact semi-simple Lie group $G$ and $H$ 
a symmetric subgroup of $G$.  
Define a Hilbert Lie group $P(G^{\Bbb C},H^{\Bbb C}\times K^{\Bbb C})$ by 
$$P(G^{\Bbb C},H^{\Bbb C}\times K^{\Bbb C}):=\{g\in H^1([0,1],G^{\Bbb C})\,\vert\,\,(g(0),g(1))
\in H^{\Bbb C}\times K^{\Bbb C}\}.$$
Then any principal orbit of the 
$P(G^{\Bbb C},H^{\Bbb C}\times K^{\Bbb C})$-action on $H^0([0,1],\mathfrak g^{\Bbb C})$ 
is an infinite dimensional anti-Kaehler isoparametric submanifold with $J$-diagonalizable shape operators.  
This fact is stated in Remark 1.1 of [K3] and shown by Theorem 1.1 (ii) in [K2] 
and Theorem B in [K3] because the $H$-action on $G/K$ is an action of Hermann type.  
In Example 2 of Section 4, we will state this fact in detail.  
In addition, for an involutive automorphism $\sigma$ of $G$, define a Hilbert Lie group 
$P(G^{\Bbb C},G(\sigma)^{\Bbb C})$ by 
$$P(G^{\Bbb C},G(\sigma)^{\Bbb C}):=\{g\in H^1([0,1],G^{\Bbb C})\,\vert\,\,(g(0),g(1))\in G(\sigma)^{\Bbb C}\},$$
where $G(\sigma):=\{(g,\sigma(g))\,\vert\,g\in G\}$.  
Then any principal orbit of $P(G^{\Bbb C},G^{\Bbb C}(\sigma))$-action on $H^0([0,1],\mathfrak g^{\Bbb C})$ 
also is an infinite dimensional anti-Kaehler isoparametric submanifold with $J$-diagonalizable shape 
operators.  This fact also is shown by Theorem 1.1 in [K2] and Theorem B in [K3] because the 
$G(\sigma)$-action on $G=(G\times G)/\triangle G$ is an action of Hermann type.  
In contrast let $G=KAN$ be the Iwasawa's decomposition of $G$, where $A$ is the abelian part and $N$ is the 
nilpotent part.  
The inverse images of orbits of the natural action $N^{\Bbb C}\curvearrowright G^{\Bbb C}/K^{\Bbb C}$ by 
$\pi\circ\phi$ are infinite dimensional anti-Kaehler isoparametric submanifolds which do not have 
$J$-diagonalizable shape operators, where $\pi$ is the natural projection of $G^{\Bbb C}$ onto 
$G^{\Bbb C}/K^{\Bbb C}$ and $\phi$ is the parallel transport map for $G^{\Bbb C}$.  See [K2] 
(or Example 2 of Section 4) about the definition of $\phi$.  Assume that a $C^{\omega}$-submanifold $M$ in 
$G/K$ has regular complex focal structure satisfying the following two conditions:

\vspace{0.3truecm}

$(\ast_1)\,\,$ 
{\sl The complex focal structure of $M$ is invariant under the parallel translation}

{\sl with respect to the normal connection of $M$}

\vspace{0.15truecm}

\noindent
and

\vspace{0.15truecm}

$(\ast_2)\,\,$
{\sl The complex focal set of $M$ at any point $x(\in M)$ consists of infinitely many}

{\sl complex hyperplanes in the complexified normal space $(T^{\perp}_x M)^{\bf c}$ and the group}

{\sl generated by the complex reflections of order two with respect to the complex}

{\sl hyperplanes is discrete.  Also, for any unit normal vector $v$ of $M$, the  nullity spa-}

{\sl ces of complex focal radii along the normal geodesic $\gamma_v$ with $\gamma'_v(0)=v$ span}

{\sl $\displaystyle{\left(({\rm Ker}\,A_v\cap{\rm Ker}\,R(v))^{\Bbb C}\right)^{\perp}}$.}

\vspace{0.3truecm}

\noindent
Then each connected component of $(\pi\circ\phi)^{-1}(M^{\Bbb C})$ is an anti-Kaehler isoparametric submanifold 
with $J$-diagonalizable shape operators.  

Recently we have proved the following extrinsic homogeneity theorem ([K7]):

\vspace{0.2truecm}

{\sl Let $M$ be a full irreducible anti-Kaehler isoparametric 
$C^{\omega}$-submanifold with $J$-diagona-\newline
lizable shape operators of codimension greater than 
one in an infinite dimensional anti-Kaehler space.  Then $M$ is extrinsically homogeneous.}

\vspace{0.2truecm}

Let $I_h(V)$ be the group of all holomorphic isometries of an infinite dimensional 
anti-Kaehler space $V$ and set $H:=\{F\in I_h(V)\,\vert\,F(M)=M\}$.  
The extrinsic homogeneity of $M$ in the above result means $Hx=M$ ($x\in M$).  
Let $I_h^b(V)$ be the subgroup of $I_h(V)$ generated by one-parameter transformation groups induced by 
holomorphic Killing vector fields defined entirely on $V$.  Note that $I_h^b(V)$ is a Banach Lie group.  
Set $H_b:=H\cap I_h^b(V)$, which is a Banach Lie subgroup of $I_h^b(V)$.  
In this paper, we prove the following extrinsic homogeneity theorem similar to the result of [GH].  

\vspace{0.5truecm}

\noindent
{\bf Theorem A.} {\sl Let $M$ be a full irreducible anti-Kaehler isoparametric 
$C^{\omega}$-submanifold with $J$-diagonalizable shape operators of codimension greater 
than one in the infinite dimensional anti-Kaehler space $V$.  
Then $M=H_bx$ holds for any $x\in M$.}

\vspace{0.5truecm}

The assumption of the $J$-diagonalizabilty of shape operators 
is essential in our method to prove Theorem A.  
It is still an open problem whether any submanifold in the statement of Theorem A is given as a principal 
orbit of the above $P(G^{\Bbb C},H^{\Bbb C}\times K^{\Bbb C})$-action or 
$P(G^{\Bbb C},G(\sigma)^{\Bbb C})$-action for some $G,\,H,\,K$ or some $G,\,\sigma$.  

\section{Basic notions and facts}
In this section, we shall recall some basic notions and facts.  

\vspace{0.25truecm}

\noindent
{\bf 2.1. Some notions associated with anti-Kaehler isoparametric submanifolds} 

Let $(V,\langle\,\,,\,\,\rangle,J)$ be an infinite dimensional anti-Kaehler space and 
$M$ an anti-Kaehler isoparametric submanifold in $V$.  See [K2] \text{red}{and [K7]} about the definitions of 
an infinite dimensional anti-Kaehler space and an anti-Kaehler isoparametric submanifold.  
Denote by $(\langle\,\,,\,\,\rangle,J)$ the anti-Kaehler structure of $M$ and 
$A$ the shape tensor of $M$.  
Fix a unit normal vector $v$ of $M$.  If there exists $X(\not=0)\in TM$ with 
$A_vX=aX+bJX$, then we call the complex number $a+b\sqrt{-1}$ a 
$J$-{\it eigenvalue of} $A_v$ 
(or a $J$-{\it principal curvature of direction} $v$) 
and call $X$ a $J$-{\it eigenvector for} $a+b\sqrt{-1}$.  
Also, we call the space of all $J$-eigenvectors for 
$a+b\sqrt{-1}$ a $J$-{\it eigenspace for} $a+b\sqrt{-1}$.  
The $J$-eigenspaces are orthogonal to one another and they are $J$-invariant, respectively.  
We call the set of all $J$-eigenvalues of $A_v$ the $J$-{\it spectrum of} 
$A_v$ and denote it by ${{\rm Spec}}_JA_v$.  
Let $\{e_i\}_{i=1}^{\infty}$ be an orthonormal system of $T_xM$.  If 
$\{e_i\}_{i=1}^{\infty}\cup\{Je_i\}_{i=1}^{\infty}$ is an orthonormal base 
of $T_xM$, then we call $\{e_i\}_{i=1}^{\infty}$ (rather than 
$\{e_i\}_{i=1}^{\infty}\cup\{Je_i\}_{i=1}^{\infty}$) a 
$J$-{\it orthonormal base}.  
If there exists a $J$-orthonormal base consisting of $J$-eigenvectors of 
$A_v$, then we say that $A_v$ {\it is diagonalized with respect to a $J$-orthonormal base} 
(or $A_v$ {\it is} $J$-{\it diagonalizable}).  
If, for each $v\in T^{\perp}M$, the shape operator $A_v$ is $J$-diagonalizable, then 
we say that $M$ {\it has} $J$-{\it diagonalizable shape operators}.  
Let $M$ be an anti-Kaehler isoparametric submanifold with $J$-diagonalizable shape operators.  
The shape operators $A_v$'s ($v\in T^{\perp}_xM$) are simultaneously diagonalized with respect to 
a $J$-orthonormal base.  
Let $\{E_0\}\cup\{E_i\,\vert\,i\in I\}$ be the family of distributions on $M$ such that, 
for each $x\in M$, 
$\{(E_0)_x\}\cup\{(E_i)_x\,\vert\,i\in I\}$ is the set of all common $J$-eigenspaces of 
$A_v$'s ($v\in T^{\perp}_xM$), where 
$\displaystyle{(E_0)_x=\mathop{\cap}_{v\in T^{\perp}_xM}{\rm Ker}\,A_v}$.  
For each $x\in M$, $T_xM$ is equal to the closure 
$\overline{\displaystyle{(E_0)_x\oplus
\left(\mathop{\oplus}_{i\in I}(E_i)_x\right)}}$ of 
$\displaystyle{(E_0)_x\oplus\left(\mathop{\oplus}_{i\in I}(E_i)_x\right)}$.  
We regard $T^{\perp}_xM$ ($x\in M$) as a complex vector space by 
$J_x\vert_{T^{\perp}_xM}$ and denote the dual space of the complex vector 
space $T^{\perp}_xM$ by $(T^{\perp}_xM)^{\ast_{\bf c}}$.  
Also, denote by $(T^{\perp}M)^{\ast_{\bf c}}$ the complex vector bundle over 
$M$ having $(T^{\perp}_xM)^{\ast_{\bf c}}$ as the fibre over $x$.  
Let $\lambda_i$ ($i\in I$) be the section of $(T^{\perp}M)^{\ast_{\bf c}}$ 
such that $A_v={{\rm Re}}(\lambda_i)_x(v){{\rm id}}+{{\rm Im}}(\lambda_i)_x(v)
J_x$ on $(E_i)_x$ for any $x\in M$ and any $v\in T^{\perp}_xM$.  
We call $\lambda_i$ ($i\in I$) $J$-{\it principal curvatures} of $M$ and 
$E_i$ ($i\in I$) $J$-{\it curvature distributions} of $M$.  
The distribution $E_i$ is integrable and each leaf of $E_i$ is a complex sphere.  
Each leaf of $E_i$ is called a {\it complex curvature sphere}.  
It is shown that there uniquely exists a normal vector field $n_i$ of $M$ with 
$\lambda_i(\cdot)=\langle n_i,\cdot\rangle-\sqrt{-1}\langle Jn_i,\cdot\rangle$.  
We call $n_i$ ($i\in I$) the $J$-{\it curvature normals} of $M$.  
Set ${\it l}^x_i:=(\lambda_i)_x^{-1}(1)$.  
Then the tangential focal set of $M$ at $x$ is equal to 
$\displaystyle{\mathop{\cup}_{i\in I}{\it l}_i^x}$ ([K2, Theorem 2 (i)]).  
We call each ${\it l}_i^x$ a {\it complex focal hyperplane of} $M$ {\it at} 
$x$.  Let $\widetilde v$ be a parallel normal vector field of $M$.  
If $\widetilde v_x$ belongs to at least one ${\it l}_i$, then it is called 
a {\it focal normal vector field} of $M$.  
For a focal normal vetor field $\widetilde v$, the focal map 
$f_{\widetilde v}$ is defined by $f_{\widetilde v}(x):=x+\widetilde v_x\,\,\,(x\in M)$.  The image 
$f_{\widetilde v}(M)$ is called a {\it focal submanifold} of $M$, which we denote by $F_{\widetilde v}$.  
For each $x\in F_{\widetilde v}$, the inverse image $f_{\widetilde v}^{-1}(x)$ 
is called a {focal leaf} of $M$.  
Denote by $T_i^x$ the complex reflection of order $2$ with respect to 
${\it l}_i^x$ (i.e., the rotation of angle $\pi$ having ${\it l}_i^x$ as the 
axis), which is an affine transformation of $T^{\perp}_xM$.  
Let ${\cal W}_x$ be the group generated by $T_i^x$'s ($i\in I$), which is an affine Weyl group.  
This group ${\cal W}_x$ is independent of the choice 
of $x\in M$ (up to group isomorphicness).  
Hence we simply denote it by ${\cal W}$.  We call this group the 
{\it complex Coxeter group associated with} $M$.  
According to Lemma 3.8 of [K4], 
${\cal W}$ is decomposable (i.e., it is decomposed into a non-trivial product of 
two discrete complex reflection groups) if and only if there exist two 
$J$-invariant linear subspaces $P_1$ ($\not=\{{\bf 0}\}$) and $P_2$ ($\not=\{{\bf 0}\}$) 
of $T^{\perp}_xM$ such that $T^{\perp}_xM=P_1\oplus P_2$ (orthogonal 
direct sum), $P_1\cup P_2$ contains all $J$-curvature normals of 
$M$ at $x$ and that $P_i$ ($i=1,2$) contains at least one $J$-
curvature normal of $M$ at $x$, where ${\bf 0}$ is the zero vector of $T^{\perp}_xM$.  
Also, $M$ is irreducible if and only if ${\cal W}$ is not decomposable ([K4, Theorem 1]).  

We note that the notions described in this subsection are defined also for a finite dimensional 
anti-Kaehler space similarly.  

\vspace{0.25truecm}

\noindent
{\bf 2.2. Aks-representation}

Let $L/H$ be an irreducible anti-Kaehler symmetric space and $(\mathfrak l,\tau)$ the 
anti-Kaehler symmetric Lie algebra associated with $L/H$.  See [K5] and [K7] about 
the definitions of these notions.  Also, set $\mathfrak p:={\rm Ker}(\tau+{\rm id})$.  
The space ${\rm Ker}(\tau-{\rm id})$ is equal to the Lie algebra $\mathfrak h$ of $H$ 
and $\mathfrak p$ is identified with $T_{eK}(L/H)$.  
Denote by ${\rm Ad}_L$ be the adjoint representation of $L$.  Define 
$\rho:H\to{\rm GL}(\mathfrak p)$ by 
$\rho(h):={\rm Ad}_L(h)\vert_{\mathfrak p}\,\,(h\in H)$.  
We call this representation $\rho$ an {\it aks}-{\it representation} 
({\it associated with} $L/H$).  Denote by ${\rm ad}_{\mathfrak h}$ the adjoint 
representation of $\mathfrak h$.  Let $\mathfrak a_s$ be a maximal split abelian 
subspace of $\mathfrak p$ (see [R] or [OS] about the definition of a maximal split 
abelian subspace) and 
$\mathfrak p=\mathfrak p_0+\sum\limits_{\alpha\in \triangle_+}\mathfrak p_{\alpha}$ the 
root space decomposition with respect to $\mathfrak a_s$, where 
the space $\mathfrak p_{\alpha}$ is defined by 
$\mathfrak p_{\alpha}:=\{X\in\mathfrak p\,\vert\,{\rm ad}_{\mathfrak l}(a)^2(X)
=\alpha(a)^2X\,\,{\rm for}\,\,{\rm all}\,\,a\in \mathfrak a_s\}$ 
($\alpha\in\mathfrak a_s^{\ast}$) and 
$\triangle_+$ is the positive root system of the root system 
$\triangle:=\{\alpha\in\mathfrak a_s^{\ast}\,\vert\,
\mathfrak p_{\alpha}\not=\{0\}\}$ under some lexicographic ordering of 
$\mathfrak a_s^{\ast}$.  
Set $\mathfrak a:=
\mathfrak p_0\,(\supset\mathfrak a_s),\,j:=J_{eK}$ and 
$\langle\,\,,\,\,\rangle_0:=\langle\,\,,\,\,\rangle_{eH}$.  
Note that $(\mathfrak p,j,\langle\,\,,\,\,\rangle_0)$ is a (finite dimmensional) anti-Kaehler space.  
It is shown that $\langle\,\,,\,\,\rangle_0\vert_{\mathfrak a_s\times\mathfrak a_s}$ is 
positive (or negative) definite, $\mathfrak a=\mathfrak a_s\oplus 
j\mathfrak a_s$ and $\langle\,\,,\,\,\rangle_0\vert_{\mathfrak a_s\times 
j\mathfrak a_s}=0$.  Note that $\mathfrak p_{\alpha}=\{X\in
\mathfrak p\,\vert\,{\rm ad}_{\mathfrak l}(a)^2(X)=\alpha^{\bf c}(a)^2X\,\,{\rm for}
\,\,{\rm all}\,\,a\in\mathfrak a\}$ holds for each $\alpha\in\triangle_+$, where 
$\alpha^{\bf c}$ is the complexification of $\alpha:\mathfrak a_s\to{\bf R}$ 
(which is a complex linear function over $\mathfrak a_s^{\bf c}=\mathfrak a$) and 
$\alpha^{\bf c}(a)^2X$ means ${\rm Re}(\alpha^{\bf c}(a)^2)X
+{\rm Im}(\alpha^{\bf c}(a)^2)jX$.  Let 
${\it l}_{\alpha}:=(\alpha^{\bf c})^{-1}(0)$ ($\alpha\in\triangle$) and 
$D:=\mathfrak a\setminus\displaystyle{
\mathop{\cup}_{\alpha\in\triangle_+}{\it l}_{\alpha}}$.  Elements of $D$ are said to be 
{\it regular}.  Take $x\in D$ and let $M$ be the orbit of the aks-representation $\rho$ 
through $x$.  From $x\in D$, $M$ is a principal orbit of this representation.  
Denote by $A$ the shape tensor of $M$.  
Take $v\in T^{\perp}_xM(=\mathfrak a)$.  Then we have 
$T_xM=\sum\limits_{\alpha\in\triangle_+}\mathfrak p_{\alpha}$ and 
$A_v\vert_{\mathfrak p_{\alpha}}=-\frac{\alpha^{\bf c}(v)}{\alpha^{\bf c}(x)}{\rm id}$ 
($\alpha\in\triangle_+$).  
Let $\widetilde v$ be the parallel normal vector field of $M$ with $\widetilde v_x=v$.  
Then we can show that $A_{\widetilde v_{\rho(h)(x)}}
\vert_{\rho(h)_{\ast x}(\mathfrak p_{\alpha})}=-\frac{\alpha^{\bf c}(v)}{\alpha^{\bf c}(x)}
{\rm id}$ for any $h\in H$.  Hence $M$ is an anti-Kaehler isoparametric submanifold 
with $J$-diagonalizable shape operators.  

\section{Homogeneity theorem} 
In this section, we shall recall the extrinsic homogeneity theorem for 
an anti-Kaehler isoparametric submanifold with $J$-diagonalizable shape operators, 
which was obtained in [K7], and the outline of its proof.  
Let $M$ be an irreducible anti-Kaehler isoparametric submanifold of codimension greater 
than one in an infinite dimensional anti-Kaehler space $(V,\langle\,\,,\,\,\rangle,J)$.  
Denote by the same symbol $(\langle\,\,,\,\,\rangle,J)$ the anti-Kaehler structure of $M$.  
Assume that $M$ has $J$-diagonalizable shape operators.  
We use the notations in Subsection 2.1.  
Denote by ${\it l}_i^x$ the complex focal hyperplane $(\lambda_i)_x^{-1}(1)$ 
of $M$ at $x$.  Also set $({\it l}_i^x)':=(\lambda_i)_x^{-1}(0)$.  
Fix $x_0\in M$.  Set ${\it l}_i:={\it l}_i^{x_0}$ and 
${\it l}'_i:=({\it l}_i^{x_0})'$.  Let $Q(x_0)$ be the set of all points 
of $M$ connected with $x_0$ by a piecewise smooth curve in $M$ each of whose smooth 
segments is contained in some complex curvature sphere (which may depend on the smooth 
segment).  By using the generalized Chow's theorem (see Theorem D of [HL2]), we showed the following fact.  

\vspace{0.5truecm}

\noindent
{\bf Lemma 3.1([K7]).} {\sl The set $Q(x_0)$ is dense in $M$.}  

\vspace{0.5truecm}

Here we note that the generalized Chow's theorem is valid because the base manifold $M$ is a Hibert 
manifold even if the metric of $M$ is a pseudo-Riemannian metric.  
For each complex affine subspace $P$ of $T^{\perp}_{x_0}M$, define $I_P$ by 
$$I_P:=\left\{
\begin{array}{ll}
\displaystyle{\{i\in I\,\vert\,(n_i)_{x_0}\in P\}} & 
\displaystyle{({\bf 0}\notin P)}\\
\displaystyle{\{i\in I\,\vert\,(n_i)_{x_0}\in P\}\cup\{0\}} & 
\displaystyle{({\bf 0}\in P).}
\end{array}\right.$$
Define a distribution $D_P$ on $M$ by 
$D_P:=\displaystyle{\mathop{\oplus}_{i\in I_P}E_i}$, which is integrable.  
Denote by $L^P_x$ the leaf through $x$ of the foliation given by $D_P$, and 
$L^i_x$ the leaf through $x$ of the foliation given by $E_i$.  
According to Lemma 4.3 of [K7], if ${\bf 0}\notin P$, then $I_P$ is finite and 
$\displaystyle{(\mathop{\cap}_{i\in I_P}{\it l}_i)\setminus
(\mathop{\cup}_{i\in I\setminus I_P}{\it l}_i)\not=\emptyset}$, and, 
if ${\bf 0}\in P$, then $I_P$ is infinite or $I_P=\{0\}$ and 
$\displaystyle{(\mathop{\cap}_{i\in I_P\setminus\{0\}}{\it l}'_i)\setminus}$\newline
$\displaystyle{(\mathop{\cup}_{i\in I\setminus I_P}{\it l}'_i)\not=\emptyset}$, where 
$\displaystyle{\mathop{\cap}_{i\in I_P\setminus\{0\}}{\it l}'_i}$ means 
$T^{\perp}_{x_0}M$ when $I_P=\{0\}$.  
Set $(W_P)_x:=x+(D_P)_x\oplus{\rm Span}_{\bf C}\{(n_i)_x\,\vert\,i\in I_P
\setminus\{0\}\}$ ($x\in M$).  Let $\gamma:[0,1]\to M$ be a piecewise smooth curve.  
Throughout this section, we assume that the domains of all piecewise smooth curves are equal to 
$[0,1]$.  If $\dot{\gamma}(t)\perp(D_P)_{\gamma(t)}$ for each $t\in[0,1]$, then $\gamma$ is 
said to be {\it perpendicular to} $D_P$ (or $D_P$-{\it perpendicular}).  
Fix $i_0\in I\cup\{0\}$ and $x_0\in M$.  
For each geodesic $\gamma:[0,1]\to L^{i_0}_{x_0}$ in $L^{i_0}_{x_0}$, we ([K7]) constructed an one-parameter 
family $\{F_{\gamma\vert_{[0,t]}}\}_{t\in[0,1]}$ of holomorphic isometries of $V$ satisfying 
$F_{\gamma\vert_{[0,t]}}(\gamma(0))=\gamma(t)$ and 
$(F_{\gamma\vert_{[0,t]}})_{\ast\gamma(0)}\vert_{T^{\perp}_{\gamma(0)}M}=\tau^{\perp}_{\gamma\vert_{[0,t]}}$ 
($t\in[0,1]$), where $\tau^{\perp}_{\gamma\vert_{[0,t]}}$ is the parallel translation along 
$\gamma\vert_{[0,t]}$ with respect to the normal connection of $M$.  
From Proposition 4.6 of [K7], the following fact holds.  

\vspace{0.5truecm}

\noindent
{\bf Lemma 3.3.} {\sl The holomorphic isometry $F_{\gamma\vert_{[0,t]}}$ preserves $M$ invariantly 
(i.e., $F_{\gamma\vert_{[0,t]}}(M)=M$).  Furthermore, it preserves $E_i$ ($i\in I$) invariantly 
(i.e., $(F_{\gamma\vert_{[0,t]}})_{\ast}(E_i)=E_i$).}

\vspace{0.5truecm}

By using Lemmas 3.2 and 3.3, we can prove the following fact (see the proof of Theorem A in [K7]).  

\vspace{0.5truecm}

\noindent
{\bf Theorem 3.4.} {\sl The submanifold $M$ is extrinsically homogeneous, that is, 
$Hx=M$ ($x\in M$) holds, where $H:=\{F\in I_h(V)\,\vert\,F(M)=M\}$.}

\section{The affine root system associated with an irreducible anti-Kaehler isoparametric submanifold}
In this section, we shall first recall the notions of the Weyl group, the affine Weyl group and 
the root system associated with a certain kind of family of the affine hyperplanes in a finite dimensional 
Euclidean affine space ${\Bbb E}$.  
Denote by $({\Bbb V},\langle\,\,,\,\,\rangle)$ the Euclidean vector space associated with ${\Bbb E}$.  
Let ${\cal H}$ be  a family of affine hyperplanes in ${\Bbb E}$ and 
${\cal W}_{\cal H}$ the group generated by the (orthogonal) reflections with respect to members 
of ${\cal H}$.  Assume that unit normal vectors of the members of ${\cal H}$ span ${\Bbb V}$ and that 
${\cal H}$ is invariant under ${\cal W}_{\cal H}$.  
Then ${\cal H}$ is a finite family of affine hyperplanes having a common point or 
a finite family of equidistant infinite parallel families of affine hyperplanes.  
In the first case, ${\cal W}_{\cal H}$ is a Weyl group and hence ${\cal H}$ is described as 
$${\cal H}=\{\alpha^{-1}(0)\,\vert\,\alpha\in\triangle\}\leqno{(4.1)}$$
for some root system $\triangle(\subset{\Bbb V}^{\ast})$ by translatiing ${\cal H}$ suitably.  
In the second case, ${\cal W}$ is an affine Weyl group and hence ${\cal H}$ is described as 
$${\cal H}=\{\alpha^{-1}(ka_{\alpha})\,\vert\,\alpha\in\triangle\,\,\&\,\,k\in{\Bbb Z}\}
\leqno{(4.2)}$$
for some root system $\triangle(\subset{\Bbb V}^{\ast})$ and some positive constants $a_{\alpha}$ 
by translating and homothetically transforming ${\cal H}$ suitably.  
Set ${\it l}_{\alpha,k}:=\alpha^{-1}(ka_{\alpha})$ 
($(\alpha,k)\in\triangle\times{\Bbb Z}$).  
Define a system ${\cal R}$ by 
$$\begin{array}{l}
\displaystyle{{\cal R}:=\{(v_{\alpha},{\it l}_{\alpha,k})\in{\Bbb V}\times{\cal H}\,\vert\,
(\alpha,k)\in\triangle\times{\Bbb Z}\}}\\
\hspace{1truecm}\displaystyle{\cup\left.\left\{\left(\frac12 v_{\alpha},{\it l}_{\alpha,k}\right)\in
{\Bbb V}\times{\cal H}\,\right\vert\,(\alpha,k)\in\triangle'\times{\Bbb Z}\right\},}
\end{array}
\leqno{(4.3)}$$
where $v_{\alpha}$ is the vector of ${\Bbb V}$ defined by 
$\alpha(\bullet)=\langle v_{\alpha},\bullet\rangle$ and $\triangle'$ is a subset of $\triangle$.  
If ${\cal R}$ is ${\cal W}$-invariant, then 
${\cal R}$ is a root system in the sense of I.G. Macdonald [M] (see Definition 7.3 of [GH] also).  
This root system ${\cal R}$ is called a {\it root system associated with} ${\cal H}$.  
In particular, if ${\cal W}$ is infinite, then it is called an 
{\it affine root system associated with} ${\cal H}$.  
If $\triangle'=\emptyset$ (resp. $\triangle'\not=\emptyset$), then ${\cal R}$ is said to be 
{\it reduced} (resp. {\it non-reduced}).  
Also, if ${\cal W}$ is irreducible (resp. reducible), then ${\cal R}$ is said to be 
{\it irreducible} (resp. {\it reducible}).  
Assume that ${\cal R}$ is a reduced irreducible affine root system of rank greater than one.  
Then the Dynkin diagram of ${\cal R}$ is defined as follows.  
Let $\Pi$ be the simple root system of $\triangle$ with respect to some lexicographic ordering of 
$V^{\ast}$ and $\delta$ be the highest root of $\triangle$ with respect to the lexicographic ordering.  
If ${\cal W}$ is finite (resp. infinite), then 
the family $\{{\it l}_{\alpha,0}\,\vert\,\alpha\in\Pi\}$ 
(resp. $\{{\it l}_{\alpha,0}\,\vert\,\alpha\in\Pi\}\cup\{{\it l}_{\delta,1}\}$) 
is the whole of walls of an alcove $C$ of ${\cal W}$-action.  
For any element $(v_{\alpha},{\it l}_{\alpha,k})$ and $(v_{\alpha'},{\it l}_{\alpha',k'})$ of 
${\cal R}$, 
$\displaystyle{\frac{\vert\vert v_{\alpha}\vert\vert}{\vert\vert v_{\alpha'}\vert\vert}
=1,\,2,\,\frac 12,\,3\,\,{\rm or}\,\,\frac 13}$ holds.  
We assign a white circle to each $\alpha\in\Pi\,\,{\rm or}\,\,\Pi\cup\{\delta\}$ 
and link the white circles corresponding to $\alpha$ and $\alpha'$ 
($\alpha,\alpha'\in\Pi\,\,{\rm or}\,\,\Pi\cup\{\delta\}$) by $1,2$ or $3$ edges 
in correspondence to
$\displaystyle{\frac{\vert\vert v_{\alpha}\vert\vert}{\vert\vert v_{\alpha'}\vert\vert}
=1,\,2^{\pm 1}\,\,{\rm or}\,\,3^{\pm 1}}$.  
Also, in the case where 
$\displaystyle{\frac{\vert\vert v_{\alpha}\vert\vert}{\vert\vert v_{\alpha'}\vert\vert}
=2^{\pm 1}\,\,{\rm or}\,\,3^{\pm 1}}$, we add the arrow pointing to the white circle corresponding to 
the shorter length one of $\alpha$ and $\alpha'$ to the $2$ or $3$ edges.  
The diagram obtained thus is called the {\it Dynkin diagram} of ${\cal R}$.  
All of reduced irreducible affine root systems of rank greater than one are 
$(\widetilde A_r)\,(r\geq 2),\,\,(\widetilde B_r)\,(r\geq 3),\,\,(\widetilde B_r^v)\,(r\geq 3),\,\,
(\widetilde C_r)\,(r\geq 2),\,\,(\widetilde C_r^v)\,(r\geq 2),\,\,(\widetilde D_r)\,(r\geq 4),\,\,
(\widetilde E_6),\,\,(\widetilde E_7),\,\,(\widetilde E_8),\,\,(\widetilde F_4),\,\,
(\widetilde F_4^v),\,\,(\widetilde G_2)$ and $(\widetilde G_2^v)$.  See Table 1 of [GH] in detail.  
Assume that ${\cal R}$ (given by $(4.3)$) is an non-reduced irreducible affine root system of 
rank greater than one.  
Define subsystems ${\cal R}_{\rm red}$ and ${\cal R}_{{\rm red}'}$ by 
$$\begin{array}{l}
\displaystyle{{\cal R}_{\rm red}:=\{(v_{\alpha},{\it l}_{\alpha,k})\in{\Bbb V}\times{\cal H}\,\vert\,
(\alpha,k)\in(\triangle\setminus\triangle')\times{\Bbb Z}\}}\\
\hspace{1truecm}\displaystyle{\cup\left.\left\{\left(\frac12 v_{\alpha},{\it l}_{\alpha,k}\right)\in
{\Bbb V}\times{\cal H}\,\right\vert\,(\alpha,k)\in\triangle'\times{\Bbb Z}\right\}}
\end{array}
\leqno{(4.4)}$$
and
$${\cal R}_{{\rm red}'}:=\{(v_{\alpha},{\it l}_{\alpha,k})\in{\Bbb V}\times{\cal H}\,\vert\,
(\alpha,k)\in\triangle\times{\Bbb Z}\}.\leqno{(4.5)}$$
Then the Dynkin diagram of ${\cal R}$ is defined as follows.  
We add the second smaller concentric white circles to the white circles corresponding to 
$\alpha$'s ($\alpha\in\Pi\cap\triangle'\,\,{\rm or}\,\,(\Pi\cup\{\delta\})\cap\triangle'$) 
in the Dynkin diagram of ${\cal R}_{\rm red}$.  
The diagram obtained thus is called the {\it Dynkin diagram} of ${\cal R}$.  
All of non-reduced irreducible affine root systems of rank greater than one are 
$(\widetilde B_r,\widetilde B_r^v)\,(r\geq 3),\,\,
(\widetilde C_r^v,\widetilde C'_r)\,(r\geq 2),\,\,
(\widetilde C'_r,\widetilde C_r)\,(r\geq 2),\,\,
(\widetilde C_r^v,\widetilde C_r)\,(r\geq 2)$ and 
$(\widetilde C_2,\widetilde C_2^v)$, 
where these notations denote the pairs of types of ${\cal R}_{\rm red}$ and ${\cal R}_{{\rm red}'}$.    
See Table 2 of [GH] in detail.  

Next we shall introduce the notion of the root system associated with 
an anti-Kaehler isoparametric submanifold with $J$-diagonalizable shape operators.  
Let $M$ be an anti-Kaehler isoparametric submanifold with $J$-diagonalizable shape operators in 
an anti-Kaehler space $V$, where $V$ may be of finite dimension.  
We use the notations in the previous section.  
Let $V=V_-\oplus V_+$ be the orthogonal decomposition of $V$ 
such that $\langle\,\,,\,\,\rangle\vert_{V_-\times V_-}$ (resp. 
$\langle\,\,,\,\,\rangle\vert_{V_+\times V_+}$) is negative (resp. positive) definite 
and that $JV_-=V_+$.  Note that such a decomposition is unique.  
Denote by $\nabla$ and $\widetilde{\nabla}$ the Riemannian connections of 
$M$ and $V$, respectively.  
Since the complex Coxeter group associated with $M$ permutes $\{{\it l}_i^x\,\vert\,i\in I\}$ and it 
is discrete, there exist a finite family $\{\mu^x_{\beta}\,\vert\,\beta\in B\}$ of complex linear 
functions over the normal space $T^{\perp}_xM$ (regarded as a complex linear space by 
$J_x$) and a finite family $\{b_{\beta}\,\vert\,\beta\in B\}$ of complex numbers such 
that $\{(\mu_{\beta}^x)^{-1}(1+b_{\beta}j)\,\vert\,\beta\in B,\,\,j\in{\Bbb Z}\}$ 
is equal to $\{{\it l}_i^x\,\vert\,i\in I\}$.  
Set $\lambda_{(\beta,j)}^x:=\frac{1}{1+b_{\beta}j}\mu_{\beta}^x$.  
Note that $(\lambda_{(\beta,j)}^x)^{-1}(1)=(\mu_{\beta}^x)^{-1}(1+b_{\beta}j)$.  
Define sections $\lambda_{(\beta,j)}$ of $(T^{\perp}M)^{\ast_{\Bbb C}}$ by assigning 
$\lambda_{(\beta,j)}^x$ to each $x\in M$.  
Set $B_0:=\{\beta\in B\,\vert\,b_{\beta}=0\}$.  
Then the set of all $J$-principal curvatures of $M$ is equal to 
$$\{\lambda_{(\beta,j)}\,\vert\,(\beta,j)\in(B\setminus B_0)\times{\Bbb Z}\}
\cup\{\lambda_{(\beta,0)}\,\vert\,\beta\in B_0\}.$$
Hence, we have $I=(B_0\times\{0\})\cup((B\setminus B_0)\times{\Bbb Z})$.  
Note that $B=B_0$ when $V$ is of finite dimension.  
Let $TM_+$ be the half-dimensional subdistribution of the tangent bundle 
$TM$ such that $\langle\,\,,\,\,\rangle\vert_{TM_+\times TM_+}$ is 
positive definite and that $\langle TM_+,JTM_+\rangle=0$, 
and set $TM_-:=JTM_+$.  Note that such subdistributions are determined uniquely.  
Similarly, we define the half-dimensional subdistributions $T^{\perp}M_{\pm}$ 
(resp. $(E_i)_{\pm}$) of the normal bundle $T^{\perp}M$ (resp. 
$J$-curvature distributions $E_i$'s ($i\in I\cup\{0\}$)).  
Clearly we have 
$$TM_-=\overline{(E_0)_-\oplus\left(
\mathop{\oplus}_{i\in I}(E_i)_-\right)}$$
and 
$$TM_+=\overline{(E_0)_+\oplus\left(\mathop{\oplus}_{i\in I}(E_i)_+\right)}.$$
Fix $x_0\in M$.  
Set $\mathfrak b:=T^{\perp}_{x_0}M$ and $\mathfrak b_{\pm}:=(T^{\perp}M_{\pm})_{x_0}$.  
Clearly we have $\mathfrak b_-=J_{x_0}\mathfrak b_+$ and 
$\mathfrak b=\mathfrak b_++\mathfrak b_-(\approx\mathfrak b_+^{\bf c})$.  

\vspace{0.5truecm}

\noindent
{\bf Lemma 4.1.} {\sl Let $i_1$ and $i_2$ be elements of $I$ such that 
$(n_{i_1})_{x_0}$ and $(n_{i_2})_{x_0}$ are linearly independent over ${\Bbb C}$.  
Set $\mathfrak b':={\rm Span}_{\Bbb R}\{(n_{i_1})_{x_0},(n_{i_2})_{x_0}\}$.  
Then we have $J_{x_0}\mathfrak b'\cap\mathfrak b'=\{0\}$.}

\vspace{0.5truecm}

\noindent
{\it Proof.} Since $(n_{i_1})_{x_0}$ and $(n_{i_2})_{x_0}$ are linearly independent over ${\Bbb C}$, 
there exists a complex affine line $P$ of $T^{\perp}_{x_0}M$ which passes through 
$(n_{i_1})_{x_0}$ and $(n_{i_2})_{x_0}$ but does not pass through $0$.  
Then $L^P_{x_0}(\subset(W_P)_{x_0})$ is a (finite dimensional) anti-Kaehler isoparametric 
submanifold with $J$-digonalizable shape operators of complex codimension greater two.  
Since the complex codimension of $L^P_{x_0}$ is equal to two, it is irreducible or 
the product of two irreducible anti-Kaehler isoparametric submanifolds 
$L^{P_i}_{x_0}(\subset(W_{P_i})_{x_0})$ ($i=1,2$) with $J$-diagonalizable shape operators of 
complex codimension one, where we note that $(W_P)_{x_0}=(W_{P_1})_{x_0}\oplus(W_{P_2})_{x_0}$.  
Also, note that $L^{P_i}(\subset(W_{P_i})_{x_0})$ ($i=1,2$) are complex spheres because 
they are of complex codimension one.  

First we consider the case where $L^P_{x_0}$ is irreducible.  
Then, according to Theorem 4.4 of [K7], $L^P_{x_0}$ is a principal orbit of the aks-representation 
associated with an irreducible anti-Kaehler symmetric space of complex rank greater than one.  
Denote by $L/H$ this irreducible anti-Kaehler symmetric space.  
We use the notations in Subsection 2.2.  
Let $L^P_{x_0}=\rho(H)\cdot w$, where 
$\rho$ is the aks-representation associated with $L/H$ and $w$ is the element of 
$\mathfrak p$ identified with $x_0$.  Let $\mathfrak a_s$ be the maximal split abelian 
subspace of $\mathfrak p$ containing $w$ and $\mathfrak a$ the Cartan subspace of $\mathfrak p$ 
containing $\mathfrak a_s$.  The space $\mathfrak a$ is identified with the normal space of 
$T^{\perp}_{x_0}L^P_{x_0}$ of $L^P_{x_0}(\subset(W_P)_{x_0})$ at $x_0$.  
Let $\triangle_+$ be the positive root system of the root system $\triangle$ 
(with respect to $\mathfrak a_s$) under some lexicographic ordering of $\mathfrak a_s^{\ast}$.  
For each $\alpha\in\triangle_+$, define the section $\lambda_{\alpha}$ of 
the ${\bf C}$-dual bundle $(T^{\perp}L^P_{x_0})^{\ast_{\Bbb c}}$ of $T^{\perp}L^P_{x_0}$ by 
$$(\lambda_{\alpha})_{\rho(h)(w)}:=-\frac{\alpha^{\Bbb C}\circ\rho(h)_{\ast w}^{-1}}
{\alpha^{\Bbb C}(w)}\quad(h\in H).$$
The set of all $J$-principal curvatures of $L^P_{x_0}$ is equal to 
$\{\lambda_{\alpha}\,\vert\,\alpha\in\triangle_+\}$.  Let $n_{\alpha}$ 
be the $J$-curvature normal corresponding to $\lambda_{\alpha}$.  
Since $(\lambda_{\alpha})_{w}=-\frac{\alpha^{\Bbb C}}{\alpha^{\Bbb C}(w)}$, 
we have $(n_{\alpha})_{x_0}\in\mathfrak a_s$ for any $\alpha\in\triangle_+$.  
This fact implies that $(n_{i_1})_{x_0}$ and $(n_{i_2})_{x_0}$ belong to 
$\mathfrak a_s$.  Hence we obtain $J_{x_0}\mathfrak b'\cap\mathfrak b'=\{0\}$.  

Next we consider the case of $L^P_{x_0}=L^{P_1}_{x_0}\times L^{P_2}_{x_0}
(\subset(W_{P_1})_{x_0}\oplus(W_{P_2})_{x_0})$.  
Then one of $(n_{i_1})_{x_0}$ and $(n_{i_2})_{x_0}$ belongs to 
$T^{\perp}_{x_0}L^{P_1}_{x_0}$ and another belongs to $T^{\perp}_{x_0}L^{P_2}_{x_0}$.  
From this fact, it follows that $J_{x_0}\mathfrak b'\cap\mathfrak b'=\{0\}$.  
This completes the proof.  
\begin{flushright}q.e.d.\end{flushright}

\vspace{0.5truecm}

Define a linear subspace $\mathfrak b_{\Bbb R}$ of $\mathfrak b$ by 
$$\mathfrak b_{\Bbb R}:={\rm Span}_{\Bbb R}\{(n_i)_{x_0}\,\vert\,i\in I\}.$$
From Lemma 4.1, it follows that $J_{x_0}\mathfrak b_{\Bbb R}\cap\mathfrak b_{\Bbb R}=\{0\}$.  
Furthermore, since $M$ is full,  $\mathfrak b_{\Bbb R}$ is a real form of $\mathfrak b$.  
For simplicity denote ${\it l}_i^{x_0}$ by ${\it l}_i$.  
It is easy to show that ${\it l}_i\cap\mathfrak b_{\Bbb R}=
((\lambda_i)_{x_0}\vert_{\mathfrak b_{\Bbb R}})^{-1}(1)$.  Denote by ${\it l}_i^{\Bbb R}$ 
this affine hyperplane ${\it l}_i\cap\mathfrak b_{\Bbb R}$ of $\mathfrak b_{\Bbb R}$.  
Let ${\cal W}_{\Bbb R}$ be the group generated by the reflections with respect to 
${\it l}_i^{\Bbb R}$'s ($i\in I$).  
It is clear that ${\cal W}_{\Bbb R}$ is isomorphic to ${\cal W}$.  
Hence, ${\cal W}_{\Bbb R}$ is an affine Weyl group.  Let $B'$ be the set of all elements $\beta$'s of $B$ 
satisfying the 

\vspace{0.5truecm}

\centerline{
\unitlength 0.1in
\begin{picture}( 53.4700, 16.1200)( -6.7300,-19.9600)
%
\special{pn 8}%
\special{pa 1408 764}%
\special{pa 2496 764}%
\special{pa 2496 1852}%
\special{pa 1408 1852}%
\special{pa 1408 764}%
\special{fp}%
%
\special{pn 8}%
\special{pa 1408 1308}%
\special{pa 2496 1308}%
\special{fp}%
%
\special{pn 8}%
\special{pa 1952 1852}%
\special{pa 1952 764}%
\special{fp}%
%
\special{pn 8}%
\special{pa 1408 1552}%
\special{pa 2496 1062}%
\special{fp}%
%
\special{pn 8}%
\special{pa 2496 1244}%
\special{pa 2306 764}%
\special{fp}%
%
\special{pn 20}%
\special{sh 1}%
\special{ar 2432 1090 10 10 0  6.28318530717959E+0000}%
\special{sh 1}%
\special{ar 2432 1090 10 10 0  6.28318530717959E+0000}%
%
\special{pn 8}%
\special{pa 3286 1270}%
\special{pa 3132 1670}%
\special{pa 4566 1252}%
\special{pa 4674 872}%
\special{pa 4674 872}%
\special{pa 4674 872}%
\special{pa 3286 1270}%
\special{fp}%
%
\special{pn 13}%
\special{pa 4158 1026}%
\special{pa 4374 1308}%
\special{fp}%
%
\special{pn 8}%
\special{pa 4158 1026}%
\special{pa 4076 500}%
\special{fp}%
%
\special{pn 8}%
\special{pa 4374 1308}%
\special{pa 4294 782}%
\special{fp}%
%
\special{pn 8}%
\special{pa 4084 518}%
\special{pa 4290 800}%
\special{fp}%
%
\special{pn 8}%
\special{pa 4376 1316}%
\special{pa 4458 1842}%
\special{fp}%
%
\special{pn 8}%
\special{pa 4256 1542}%
\special{pa 4462 1824}%
\special{fp}%
%
\special{pn 8}%
\special{pa 4256 1542}%
\special{pa 4220 1380}%
\special{fp}%
%
\special{pn 8}%
\special{pa 4158 1026}%
\special{pa 4212 1326}%
\special{fp}%
\put(19.0600,-19.9600){\makebox(0,0)[lt]{$\mathfrak b$}}%
%
\special{pn 8}%
\special{pa 2514 636}%
\special{pa 2342 844}%
\special{dt 0.045}%
\special{sh 1}%
\special{pa 2342 844}%
\special{pa 2400 806}%
\special{pa 2376 804}%
\special{pa 2368 780}%
\special{pa 2342 844}%
\special{fp}%
%
\special{pn 8}%
\special{pa 1924 590}%
\special{pa 2424 1080}%
\special{dt 0.045}%
\special{sh 1}%
\special{pa 2424 1080}%
\special{pa 2390 1020}%
\special{pa 2386 1044}%
\special{pa 2362 1048}%
\special{pa 2424 1080}%
\special{fp}%
%
\special{pn 8}%
\special{pa 1598 646}%
\special{pa 1724 1308}%
\special{dt 0.045}%
\special{sh 1}%
\special{pa 1724 1308}%
\special{pa 1732 1238}%
\special{pa 1714 1256}%
\special{pa 1692 1246}%
\special{pa 1724 1308}%
\special{fp}%
%
\special{pn 8}%
\special{pa 1298 862}%
\special{pa 1598 1462}%
\special{dt 0.045}%
\special{sh 1}%
\special{pa 1598 1462}%
\special{pa 1586 1392}%
\special{pa 1574 1414}%
\special{pa 1550 1410}%
\special{pa 1598 1462}%
\special{fp}%
%
\special{pn 8}%
\special{pa 1298 1788}%
\special{pa 1952 1660}%
\special{dt 0.045}%
\special{sh 1}%
\special{pa 1952 1660}%
\special{pa 1882 1654}%
\special{pa 1900 1670}%
\special{pa 1890 1692}%
\special{pa 1952 1660}%
\special{fp}%
\put(25.3200,-5.9000){\makebox(0,0)[lb]{${\it l}_i$}}%
\put(19.5100,-5.6300){\makebox(0,0)[rb]{${\it l}_i^{\Bbb R}$}}%
\put(16.7000,-6.0900){\makebox(0,0)[rb]{$\mathfrak b_+$}}%
\put(13.0700,-8.2600){\makebox(0,0)[rb]{$\mathfrak b_{\Bbb R}$}}%
\put(12.4300,-17.6000){\makebox(0,0)[rt]{$\mathfrak b_-$}}%
\put(30.2000,-6.1000){\makebox(0,0)[lb]{in fact}}%
\put(44.0200,-5.5400){\makebox(0,0)[lb]{${\it l}_i$}}%
%
\special{pn 8}%
\special{pa 4384 554}%
\special{pa 4166 808}%
\special{dt 0.045}%
\special{sh 1}%
\special{pa 4166 808}%
\special{pa 4226 770}%
\special{pa 4202 768}%
\special{pa 4194 744}%
\special{pa 4166 808}%
\special{fp}%
%
\special{pn 8}%
\special{pa 3840 1706}%
\special{pa 4276 1190}%
\special{dt 0.045}%
\special{sh 1}%
\special{pa 4276 1190}%
\special{pa 4218 1228}%
\special{pa 4242 1230}%
\special{pa 4248 1254}%
\special{pa 4276 1190}%
\special{fp}%
\put(38.4800,-17.5100){\makebox(0,0)[rt]{${\it l}_i^{\Bbb R}$}}%
%
\special{pn 8}%
\special{pa 4602 736}%
\special{pa 4492 1054}%
\special{dt 0.045}%
\special{sh 1}%
\special{pa 4492 1054}%
\special{pa 4534 996}%
\special{pa 4510 1004}%
\special{pa 4496 984}%
\special{pa 4492 1054}%
\special{fp}%
\put(45.8300,-7.0800){\makebox(0,0)[lb]{$\mathfrak b_{\Bbb R}$}}%
%
\special{pn 8}%
\special{ar 3300 1160 860 490  3.4181571 3.4359349}%
\special{ar 3300 1160 860 490  3.4892682 3.5070460}%
\special{ar 3300 1160 860 490  3.5603794 3.5781571}%
\special{ar 3300 1160 860 490  3.6314905 3.6492682}%
\special{ar 3300 1160 860 490  3.7026016 3.7203794}%
\special{ar 3300 1160 860 490  3.7737127 3.7914905}%
\special{ar 3300 1160 860 490  3.8448238 3.8626016}%
\special{ar 3300 1160 860 490  3.9159349 3.9337127}%
\special{ar 3300 1160 860 490  3.9870460 4.0048238}%
\special{ar 3300 1160 860 490  4.0581571 4.0759349}%
\special{ar 3300 1160 860 490  4.1292682 4.1470460}%
\special{ar 3300 1160 860 490  4.2003794 4.2181571}%
\special{ar 3300 1160 860 490  4.2714905 4.2892682}%
\special{ar 3300 1160 860 490  4.3426016 4.3603794}%
\special{ar 3300 1160 860 490  4.4137127 4.4314905}%
\special{ar 3300 1160 860 490  4.4848238 4.5026016}%
\special{ar 3300 1160 860 490  4.5559349 4.5737127}%
\special{ar 3300 1160 860 490  4.6270460 4.6448238}%
\special{ar 3300 1160 860 490  4.6981571 4.7159349}%
\special{ar 3300 1160 860 490  4.7692682 4.7870460}%
\special{ar 3300 1160 860 490  4.8403794 4.8581571}%
\special{ar 3300 1160 860 490  4.9114905 4.9292682}%
\special{ar 3300 1160 860 490  4.9826016 5.0003794}%
\special{ar 3300 1160 860 490  5.0537127 5.0714905}%
\special{ar 3300 1160 860 490  5.1248238 5.1426016}%
\special{ar 3300 1160 860 490  5.1959349 5.2137127}%
\special{ar 3300 1160 860 490  5.2670460 5.2848238}%
\special{ar 3300 1160 860 490  5.3381571 5.3559349}%
\special{ar 3300 1160 860 490  5.4092682 5.4270460}%
\special{ar 3300 1160 860 490  5.4803794 5.4981571}%
\special{ar 3300 1160 860 490  5.5514905 5.5692682}%
\special{ar 3300 1160 860 490  5.6226016 5.6403794}%
\special{ar 3300 1160 860 490  5.6937127 5.7114905}%
\special{ar 3300 1160 860 490  5.7648238 5.7826016}%
%
\special{pn 8}%
\special{pa 4060 940}%
\special{pa 4090 970}%
\special{fp}%
\special{sh 1}%
\special{pa 4090 970}%
\special{pa 4058 910}%
\special{pa 4052 932}%
\special{pa 4030 938}%
\special{pa 4090 970}%
\special{fp}%
\end{picture}%
\hspace{5.1truecm}
}

\vspace{1truecm}

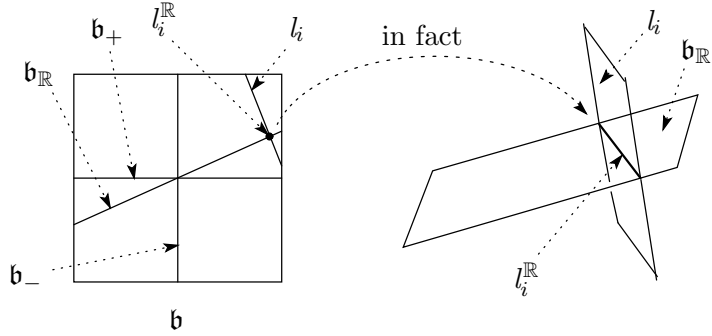
\centerline{{\bf Figure 1.}}

\vspace{0.5truecm}

\noindent
following condition:

\vspace{0.2truecm}

There exists $\hat{\beta}\in B$ such that $(n_{(\beta,0)})_{x_0}$ and $(n_{(\hat{\beta},0)})_{x_0}$ 
are linearly independent 

over ${\Bbb C}$, for the the complex affine line $P$ through 
$(n_{(\beta,0)})_{x_0}$ and $(n_{(\hat{\beta},0)})_{x_0}$, the 

root system associated with $L^P_{x_0}(\subset W_P)$ is of type $({\rm BC}_2)$ 
and the $\frac 12$-multiple 

of the root $\alpha\in\triangle_+$ ($\triangle_+\,:\,$as in the proof of 
Lemma 4.1) corresponding to $\beta$ 

also belongs to $\triangle_+$.

\vspace{0.2truecm}

\noindent
Fix $\displaystyle{Z_0\in\mathop{\cap}_{\beta\in B}{\it l}_{\beta}^{\Bbb R}}$.  
There exists a root system $\triangle_M\,(\subset(\mathfrak b_{\Bbb R})^{\ast})$ such that 
$$\begin{array}{l}
\displaystyle{\left.\left\{-\frac{\alpha}{\alpha(Z_0)}\,\right\vert\,\alpha\in(\triangle_M)_+\right\}
\cup\left.\left\{-\frac{\alpha}{2\alpha(Z_0)}\,\right\vert\,
\alpha\in(\triangle_M)_+\,\,{\rm s.t.}\,\,\frac{\alpha}{2}\in(\triangle_M)_+\right\}}\\
\displaystyle{=\{\lambda_{(\beta,0)}\vert_{\mathfrak b_{\Bbb R}}\,\vert\,\beta\in B\}
\cup\left.\left\{\frac12\lambda_{(\beta,0)}\vert_{\mathfrak b_{\Bbb R}}\,\right\vert\,\beta\in B'
\right\},}
\end{array}$$
where $(\triangle_M)_+$ is the positive root system of $\triangle_M$ under a lexicographic ordering 
of $(\mathfrak b_{\Bbb R})^{\ast}$.  
When $\alpha(\in(\triangle_M)_+)$ corresponds to $\beta\in B$ 
(i.e., $-\frac{\alpha}{\alpha(Z_0)}=\lambda_{(\beta,0)}\vert_{\mathfrak b_{\Bbb R}}$), we denote 
$\lambda_{(\beta,j)},\,n_{(\beta,j)},$\newline
${\it l}_{(\beta,j)}$ and $b_{\beta}$ by 
$\lambda_{(\alpha,j)},\,n_{(\alpha,j)},\,{\it l}_{(\alpha,j)}$ and $b_{\alpha}$, respectively.  
Hence we may denote $(\triangle_M)_+\times{\Bbb Z}$ by $I$.  
In the sequel, $I$ denotes $(\triangle_M)_+\times{\Bbb Z}$.  
Define a system ${\cal R}_M$ by 
$${\cal R}_M:=\{((n_{(\alpha,0)})_{x_0},\,\,{\it l}_{(\alpha,j)}^{\Bbb R})\,\vert\,
\alpha\in(\triangle_M)_+,\,\,j\in{\Bbb Z}\}.$$
This root system ${\cal R}_M$ is a root system associated with ${\cal H}$.  
In particular, if $B_0\not=B$, then it is an affine root system associated with ${\cal H}$.  

\vspace{0.5truecm}

\noindent
{\bf Definition.} We call ${\cal R}_M$ {\it the root system associated with} $M$.  
In particular, if $B\not=B_0$, then we call ${\cal R}_M$ 
{\it the affine root system associated with} $M$.  

\vspace{0.5truecm}

For ${\cal R}_M$, the following fact holds.  

\vspace{0.5truecm}

\noindent
{\bf Proposition 4.2.} {\sl 
If $M$ is irreducible, then ${\cal W}$ is infinite and hence ${\cal R}_M$ is 
the affine root system.}

\vspace{0.5truecm}

\noindent
{\it Proof.} 
To show this statement, we suffice to show that $B\not=B_0$.  
Suppose that $B=B_0$.  Then we have 
$$T_{x_0}M=(E_0)_{x_0}\oplus\left(\mathop{\oplus}_{\beta\in B}
(E_{(\beta,0)})_{x_0}\right).$$
This implies that $M$ is the cylinder over a finite dimensional anti-Kaehler isoparametric 
submanifold of $J$-diagonalizable shape operators.  
This contradicts the fact that $M$ is irreducible.  
Hence we obtain $B\not=B_0$.  
\begin{flushright}q.e.d.\end{flushright}

\vspace{0.5truecm}

\noindent
{\it Example 1.} 
Let $(L,H)$ be an anti-Kaehler symmetric pair and 
$\rho:H\to{\rm GL}(\mathfrak p)$ the aks-representation associated with $(L,H)$, where $\mathfrak p$ 
is as in Subsection 2.2.  We use the notations in Subsection 2.2.  
Let $M$ be the orbit of $\rho(H)$-action through a regular element $x_0(\in\mathfrak a)$ and 
$V$ an infinite dimensional anti-Kaehler space.  
Then the cylinder $M\times V(\subset\mathfrak p\times V)$ over $M$ is 
a (reducible) anti-Kaehler isoparametric submanifold with $J$-diagonalizable shape operators.  
The set ${\cal JPC}_{M\times V}$ of all $J$-principal curvatures of $M\times V$ is given by 
$${\cal JPC}_{M\times V}=\left\{\left.
-\frac{\widetilde{\alpha^{\Bbb C}}}{\alpha(x_0)}\,\right\vert\,\alpha\in\triangle_+\right\},$$
where $\widetilde{\alpha^{\Bbb C}}$ is the parallel section of $(T^{\perp}M)^{\ast_{\Bbb C}}$ with 
$\displaystyle{\left(\widetilde{\alpha^{\Bbb C}}\right)_{x_0}=\alpha^{\Bbb C}}$.  
Hence we have 
$${\cal H}=\{\alpha^{-1}(-\alpha(x_0))\,\vert\,\alpha\in\triangle_+\},$$
and 
$${\cal R}_M=\{((n_{\alpha})_{x_0},\,\,\alpha^{-1}(-\alpha(x_0)))\,\vert\,
\alpha\in\triangle_+\},$$
where $(n_{\alpha})_{x_0}$ is the element of $\mathfrak a_s$ with 
$\alpha(\bullet)=\langle(n_{\alpha})_{x_0},\bullet\rangle$.  
Also, we have $\triangle_M=\triangle$.  
Thus both the types of $\triangle_M$ and ${\cal R}_M$ are equal to that of $\triangle$.  

\vspace{0.5truecm}

\noindent
{\it Example 2.} Let $G/K$ be a symmetric space of non-compact type and 
$H\curvearrowright G/K$ a Hermann type action (i.e., $H$ is a symmetric subgroup of $G$).  
Let $\mathfrak g,\mathfrak k$ and $\mathfrak h$ be the Lie algebras of $G,K$ and $H$, 
and $\theta$ (resp. $\sigma$) the involution of $G$ with 
$({\rm Fix}\,\theta)_0\subset K\subset{\rm Fix}\,\theta$ 
(resp. $({\rm Fix}\,\sigma)_0\subset H\subset{\rm Fix}\,\sigma$).  
Denote by the same symbols the involutions of $\mathfrak g$ induced from 
$\theta$ and $\sigma$.  
Set $\mathfrak p:={\rm Ker}(\theta+{\rm id})$ and $\mathfrak q:={\rm Ker}(\sigma+{\rm id})$.  
Assume that $\theta$ and $\sigma$ commute.  
Then we have $\mathfrak p=\mathfrak p\cap\mathfrak h+\mathfrak p\cap\mathfrak q$.  
Take a maximal abelian $\mathfrak b'$ of $\mathfrak p\cap\mathfrak q$.  
Let $\mathfrak p=\mathfrak z_{\mathfrak p}(\mathfrak b')
+\sum\limits_{\alpha\in\triangle'_+}\mathfrak p_{\alpha}$ be the root space 
decomposition with respect to $\mathfrak b'$, where 
$\mathfrak z_{\mathfrak p}(\mathfrak b')$ is the centralizer of $\mathfrak b'$ 
in $\mathfrak p$, $\triangle'_+$ is the positive root system of the root system 
$\triangle':=\{\alpha\in{\mathfrak b'}^{\ast}\,\vert\,\exists\,X(\not=0)\in
\mathfrak p\,\,{\rm s.t.}\,\,{\rm ad}(b)^2(X)=\alpha(b)^2X\,\,(\forall\,b\in
\mathfrak b')\}$ under some lexicographic ordering of ${\mathfrak b'}^{\ast}$ and 
$\mathfrak p_{\alpha}:=\{X\in\mathfrak p\,\vert\,{\rm ad}(b)^2(X)=\alpha(b)^2X
\,\,(\forall\,b\in\mathfrak b')\}$ ($\alpha\in\triangle'_+$).  Also, let 
${\triangle'}^V_+:=\{\alpha\in\triangle'_+\,\vert\,\mathfrak p_{\alpha}\cap
\mathfrak q\not=\{0\}\}$ and ${\triangle'}^H_+:=\{\alpha\in\triangle'_+\,\vert\,
\mathfrak p_{\alpha}\cap\mathfrak h\not=\{0\}\}$.  
Also, let $\phi:H^0([0,1],\mathfrak g^{\Bbb C})\to G^{\Bbb C}$ be the parallel transport map for 
$G^{\Bbb C}$ and $\pi:G^{\Bbb C}\to G^{\Bbb C}/K^{\Bbb C}$ the natural projection.  
See [K2] about the definition of the parallel transport map for $G^{\Bbb C}$.  
Let $H^{\Bbb C}\curvearrowright G^{\Bbb C}/K^{\Bbb C}$ be the complexified action of the $H$-action, 
$M$ the principal orbit of the $H^{\Bbb C}$-action through ${\rm Exp}\,Z_0$ 
and $\widetilde M$ a connected component of $(\pi\circ\phi)^{-1}(M)$, where 
$Z_0$ is a point of $\mathfrak b:={\mathfrak b'}^{\Bbb C}(=T^{\perp}_{eK^{\Bbb C}}M)$ 
($e\,:\,$the identity element of $G^{\Bbb C}$) and ${\rm Exp}$ is the exponential map of 
$G^{\Bbb C}/K^{\Bbb C}$ at $eK^{\Bbb C}$.  
Note that $\widetilde M$ is a principal orbit of the $P(G^{\Bbb C},H^{\Bbb C}\times K^{\Bbb C})$-action 
stated in Introduction.  
This submanifold $\widetilde M$ is an anti-Kaehler isoparametric submanifold with 
$J$-diagonalizable shape operators in $H^0([0,1],\mathfrak g^{\Bbb C})$.  
In particular, if $G/K$ is irreducible, then $\widetilde M$ is (extrinsically) irreducible.  
Fix $u_0\in(\pi\circ\phi)^{-1}(x_0)\cap\widetilde M$.  
By the similar argument to Section 4 of [K6], it is shown that the set 
${\cal JPC}_{\widetilde M}$ of all $J$-principal curvatures of $\widetilde M$ is given by 
$$\begin{array}{l}
\displaystyle{{\cal JPC}_{\widetilde M}=\left\{\left.
-\frac{\widetilde{\alpha^{\Bbb C}}}{\alpha(Z_0)+k\pi\sqrt{-1}}\,\right\vert\,
\alpha\in{\triangle'}^V_+,\,\,k\in{\Bbb Z}\right\}}\\
\hspace{1.8truecm}\displaystyle{\cup\left\{\left.
-\frac{\widetilde{\alpha^{\Bbb C}}}{\alpha(Z_0)+(k+\frac12)\pi\sqrt{-1}}\,\right\vert\,\alpha\in
{\triangle'}^H_+,\,\,k\in{\Bbb Z}\right\},}
\end{array}\leqno{(4.6)}$$
where $\widetilde{\alpha^{\Bbb C}}$ is the parallel section of 
$(T^{\perp}\widetilde M)^{\ast_{\Bbb C}}$ with 
$\displaystyle{\left(\widetilde{\alpha^{\Bbb C}}\right)_{u_0}=\alpha^{\Bbb C}}$.  
Here the normal space $T^{\perp}_{u_0}\widetilde M$ of $\widetilde M$ at $u_0$ 
is identified with $T^{\perp}_{x_0}M(=\mathfrak b)$ through $(\pi\circ\phi)_{\ast u_0}$.  
Define a complex linear  function $\lambda_{(\alpha,0)}$ over 
$\mathfrak b(={\mathfrak b'}^{\Bbb C})$ 
by $\lambda_{(\alpha,0)}:=-\frac{\widetilde{{\alpha}^{\Bbb C}}}{\alpha(Z_0)}$, which 
is a $J$-principal curvature of $\widetilde M$.  
Let $n_{(\alpha,0)}$ be the $J$-curvature normal of $\widetilde M$ corresponding to 
$\lambda_{(\alpha,0)}$.  From $(4.6)$, we have 
$$\begin{array}{l}
\displaystyle{{\cal H}=\left.\left\{\alpha^{-1}
(-\alpha(Z_0)+k\pi\sqrt{-1})\,\right\vert\,\alpha\in{\triangle'}^V_+,\,k\in{\Bbb Z}\right\}}\\
\hspace{0.7truecm}
\displaystyle{\cup\left.\left\{\alpha^{-1}(-\alpha(Z_0)+(k+\frac 12)\pi\sqrt{-1})\,\right\vert\,
\alpha\in{\triangle'}^H_+,\,k\in{\Bbb Z}\right\}}
\end{array}$$
and 
$$\begin{array}{l}
\displaystyle{{\cal R}_M=\left.\left\{\left((n_{(\alpha,0)})_{u_0},\,\,
\alpha^{-1}(-\alpha(Z_0)+k\pi\sqrt{-1})\right)\,\right\vert\,\alpha\in{\triangle'}^V_+,\,
k\in{\Bbb Z}\right\}}\\
\hspace{0.6truecm}\displaystyle{\cup\left.\left\{\left((n_{(\alpha,0)})_{u_0},\,\,
\alpha^{-1}(-\alpha(Z_0)+(k+\frac 12)\pi\sqrt{-1})\right)\,
\right\vert\,\alpha\in{\triangle'}^H_+,\,k\in{\Bbb Z}\right\},}\\
\hspace{0.6truecm}\displaystyle{\cup\left.\left\{\left(\frac12(n_{(\alpha,0)})_{u_0},\,\,
\alpha^{-1}
(-\alpha(Z_0)+k\pi\sqrt{-1})\right)\,
\right\vert\,\alpha\in({\triangle'}^V_+)',\,k\in{\Bbb Z}\right\},}\\
\hspace{0.6truecm}\displaystyle{\cup\left.\left\{\left(\frac12(n_{(\alpha,0)})_{u_0},\,\,
\alpha^{-1}
(-\alpha(Z_0)+(k+\frac 12)\pi\sqrt{-1})\right)\,
\right\vert\,\alpha\in({\triangle'}^H_+)',\,k\in{\Bbb Z}\right\},}
\end{array}$$
where $({\triangle'}^V_+)':=\{\alpha\in{\triangle'}^V_+\,\vert\,\frac12\alpha\in\triangle'_+\}$ and 
$({\triangle'}^H_+)':=\{\alpha\in{\triangle'}^H_+\,\vert\,\frac12\alpha\in\triangle'_+\}$.  
Also, we have $\triangle_M=\triangle'$.  

\section{Proof of Theorem A}
Let $M(\subset V)$ be as in Theorem A.  
We use the notations in Sections 3 and 4.  
Note that $I=(\triangle_M)_+\times{\Bbb Z}$.  
For simplicity denote ${\cal R}_M$ by ${\cal R}$.  
Let $P$ be a complex affine subspace of $\mathfrak b=T^{\perp}_{x_0}M$ and 
$D_P$ a distribution on $M$ defined 
in Section 3.  Then it is easy to show that $D_P$ is a totally geodesic distribution on $M$.  
We call the integral manifold $L^P_x$ of $D_P$ through $x$ a {\it slice} of $M$.  
Denote by ${\bf 0}$ the origin of $\mathfrak b$.  
If ${\bf 0}\notin P$, then $L^P_x$ is a focal leaf.  
Then, since $L^P_{x_0}$ is a finite dimensional anti-Kaehler isoparametric submanifold 
with $J$-diagonalizable shape operators of codimension greater than one in $(W_P)_{x_0}$, it is 
the product of principal orbits of the aks-representations associated with some irreducible 
anti-Kaehler symmetric spaces by Theorem 4.4 in [K7], where we use also the fact that 
a finite dimensional anti-Kaehler isoparametric (complex) hypesurface is a complex sphere 
(i.e., a principal orbit of the aks-representation associated with an anti-Kaehler symmetric space 
of complex rank one).  
If ${\bf 0}\in P$, then the slice $L^P_{x_0}$ is an infinite dimensional 
anti-Kaehler isoparametric submanifold with $J$-diagonalizable shape operators in 
$(W_P)_{x_0}$.  
Take any $w_0\in(E_i)_{x_0}$ ($i\in I$).  
Let $\gamma:[0,1]\to L^i_{x_0}$ be the geodesic in $L^i_{x_0}$ with 
$\gamma'(0)=w_0$ and 
$\{F_{{\gamma}\vert_{[0,t]}}\}_{t\in{\Bbb R}}$ the one-parameter family of holomorphic isometries of $V$ 
stated in Section 3.  
For simplicity set $F^{w_0}_t:=F_{\gamma_{[0,t]}}$.  Let $X^{w_0}$ be 
the holomorphic Killing field associated with the one-parameter transformation 
group $\{F^{w_0}_t\}_{t\in{\bf R}}$, that is, 
$\displaystyle{X^{w_0}_x:=\left.\frac{d}{dt}\right\vert_{t=0}F^{w_0}_t(x)}$, 
where $x$ moves over the set (which we denote by $U$) of all elements $x$'s 
where $\displaystyle{\left.\frac{d}{dt}\right\vert_{t=0}F^{w_0}_t(x)}$ exists.  
Set $\displaystyle{A^{w_0}:=\left.\frac{d}{dt}\right\vert_{t=0}
(F^{w_0}_t)_{\ast x_0}}$ and $b^{w_0}:=(X^{w_0})_{\bf 0}$, where ${\bf 0}$ in $(X^{w_0})_{\bf 0}$ 
is the zero element of $V$ (i.e., $(X^{w_0})_x=A^{w_0}x+b^{w_0}$).  Clearly we have 
$$\left(\mathop{\oplus}_{i\in I\cup\{0\}}(E_i)_{x_0}\right)
\oplus \mathfrak b\subset U,$$
where we regard the left-hand side as a subspace of $V$ under the identification of 
$T_{x_0}V$ and $V$.  
However, $U$ does not necessarily coincide with the whole of $V$.  
For simplicity we set 
$V'_{x_0}:=\displaystyle{\left(\mathop{\oplus}_{i\in I\cup\{0\}}(E_i)_{x_0}\right)
\oplus \mathfrak b}$ and 
$(V'_{x_0})_T:=\displaystyle{\mathop{\oplus}_{i\in I\cup\{0\}}(E_i)_{x_0}}$.  
Define a map $\overline{\Gamma}_{w_0}:(V'_{x_0})_T\to V$ by 
$\overline{\Gamma}_{w_0}(w):=\displaystyle{\left.\frac{d}{dt}
\right\vert_{t=0}(F^{w_0}_t)_{\ast x_0}(w)(=A^{w_0}w)\,\,\,\,(w\in (V'_{x_0})_T)}$ 
and a map $\Gamma_{w_0}:(V'_{x_0})_T\to T_{x_0}M$ by 
$\Gamma_{w_0}w:=(\overline{\Gamma}_{w_0}w)^T\,\,\,\,(w\in(V'_{x_0})_T)$, 
where $(\cdot)^T$ is the $T_{x_0}M$-component of $(\cdot)$.  
Also, by using $\overline{\Gamma}_w$'s 
($\displaystyle{w\in\mathop{\cup}_{i\in I}(E_i)_{x_0}}$), we define a map 
$\displaystyle{\overline{\Gamma}^{x_0}:
\left(\mathop{\oplus}_{i\in I}(E_i)_{x_0}\right)\times (V'_{x_0})_T\to V}$ 
by setting $\overline{\Gamma}^{x_0}_{w_1}w_2:=\overline{\Gamma}_{w_1}(w_2)$ 
($\displaystyle{w_1\in\mathop{\cup}_{i\in I}(E_i)_{x_0},\,w_2\in (V'_{x_0})_T}$) and 
extending linearly with respect to the first component.  Similarly, by using $\Gamma_w$'s 
($\displaystyle{w\in\mathop{\cup}_{i\in I}(E_i)_{x_0}}$), we define a map 
$\displaystyle{\Gamma^{x_0}:
\left(\mathop{\oplus}_{i\in I}(E_i)_{x_0}\right)\times (V'_{x_0})_T\to T_{x_0}M}$.  
This map $\Gamma^{x_0}$ is called the {\it homogeneous structure} of $M$ at $x_0$.  

In this section, we prove the following fact.  

\vspace{0.5truecm}

\noindent
{\bf Theorem 5.1.} {\sl The holomorphic Killing field $X^{w_0}$ is defined 
on the whole of $V$.}

\vspace{0.5truecm}

For simplicity we denote the extrinsically homogeneous structure $\Gamma^{x_0}$ by $\Gamma$.  
Denote by $h$ the second fundamental form of $M$.  
It is clear that $\overline{\Gamma}_{w_0}w=\Gamma_{w_0}w+h(w_0,w)$ 
($w\in V'_T$) and that $h(w_0,\cdot)$ is defined on the whole of $T_{x_0}M$.  
Hence, in order to show this theorem, 
we suffice to show that $\Gamma_{w_0}(:(V'_{{x_0}})_T\to T_{{x_0}}M)$ is defined (continuously) on 
the whole of $T_{x_0}M$.  
Since $(T_{{x_0}}M,\langle\,\,,\,\,\rangle)$ is an anti-Kaehler space, 
$(T_{{x_0}}M,\,-{\rm pr}_{(T_{{x_0}}M)_-}^{\ast}\langle\,\,,\,\,\rangle
+{\rm pr}_{(T_{{x_0}}M)_+}^{\ast}\langle\,\,,\,\,\rangle)$ is a Hilbert space, where 
${\rm pr}_{(T_{{x_0}}M)_{\pm}}$ is the orthogonal projection of $T_{{x_0}}M$ onto 
$(T_{{x_0}}M)_{\pm}$.  Set 
$\langle\,\,,\,\,\rangle_{\pm}:=-{\rm pr}_{(T_{{x_0}}M)_-}^{\ast}\langle\,\,,\,\,\rangle
+{\rm pr}_{(T_{{x_0}}M)_+}^{\ast}\langle\,\,,\,\,\rangle$.  Denote by 
$\vert\vert\bullet\vert\vert$ the norm of a vector of $T_{{x_0}}M$ 
with respect to $\langle\,\,,\,\,\rangle_{\pm}$ and the operator norm 
of a linear transformation from $(V'_{{x_0}})_T$ to $T_{{x_0}}M$ with respect to 
$\langle\,\,,\,\,\rangle_{\pm}$.  
To show that $\Gamma_{w_0}(:(V'_{{x_0}})_T\to T_{{x_0}}M)$ is defined (continuously) on 
the whole of $T_{x_0}M$, we suffice to show that it is bounded with respect to 
$\vert\vert\bullet\vert\vert$.  In the sequel, we shall prove the boundedness of 
$\Gamma_{w_0}$ with respect to $\vert\vert\bullet\vert\vert$ 
by the similar argument to [GH].  
Even if the proof is similar to that of [GH], we need to discuss it carefully.  
For the domain of $\Gamma$ is an anti-Kaehler space but there exist some parts discussed on 
a special real form of the space.  
Some of facts corresponding to lemmas and propositions in Sections 3-6 and 8 of [GH] are shown 
in the same methods as their proofs in [GH].  We shall state the facts as lemmas without the proof.  

\vspace{0.5truecm}

For $\Gamma$, we can show the following fact.  

\vspace{0.5truecm}

\noindent
{\bf Lemma 5.2.} {\sl Let $i_1\in I$ and $i_2,i_3\in I\cup\{0\}$.  

{\rm (i)} For any $w_k\in(E_{i_k})_{x_0}$ ($k=1,2,3$), we have 
$$\langle\Gamma_{w_1}w_2,w_3\rangle+\langle w_2,\Gamma_{w_1}w_3\rangle=0,$$

{\rm (ii)} For any $w_k\in(E_{i_k})_{x_0}$ ($k=1,2$) and 
any holomorphic isometry $f$ of $V$ preserving $M$ 
invariantly, we have 
$$f_{\ast}\Gamma_{w_1}w_2=\Gamma_{f_{\ast}w_1}f_{\ast}w_2.$$
}

\vspace{0.5truecm}

Also, for $F^{w_0}_t$, we have the following fact.  

\vspace{0.5truecm}

\noindent
{\bf Lemma 5.3.} {\sl Let $L$ be a slice of $M$, 
$i_0$ an element of $I\cup\{0\}$ with $(E_{i_0})_{x_0}\subset T_{x_0}L$ 
and $W$ the complex affine span of $L$.  If $w_0\in(E_{i_0})_{x_0}$, then 
$F^{w_0}_t(L)=L$ holds for all $t\in[0,1]$ and $X^{w_0}$ is tangent to $W$ 
along $W$.  
Furthermore, if $L$ is irreducible and is of rank greater than one, then 
$F^{w_0}_t\vert_W={}^LF^{w_0}_t$ holds for all $t\in[0,1]$, where 
${}^LF^{w_0}_t$ is the one-parameter transformation group of $W$ defined for $L$ in similar to 
$F^{w_0}_t$, and hence the extrinsically homogeneous structure of $L(\subset W)$ 
at ${x_0}$ is the restriction of $\Gamma$.}

\vspace{0.5truecm}

These lemmas are proved in the methods of the proofs of Lemmas 3.4 and 3.5 of [GH], respectively.  
Let $\widetilde v$ be a (non-focal) parallel normal vector field of 
$M$, $\eta_{\widetilde v}:M\to V$ 
the end-point map for $\widetilde v$ 
(i.e., $\eta_{\widetilde v}(u):=\exp^{\perp}(\widetilde v_u)
\,\,\,(u\in M))$ and 
$M_{\widetilde v}$ the parallel submanifold for 
$\widetilde v$ (i.e., the image of $\eta_{\widetilde v}$).  
Denote by ${}^{\widetilde v}\Gamma$ the extrinsically homogeneous structure of 
$M_{\widetilde v}$ at $\eta_{\widetilde v}({x_0})$.  
Then we have the following fact.  

\vspace{0.5truecm}

\noindent
{\bf Lemma 5.4.} {\sl For any $w_1\in(E_{i_1})_{x_0}$ ($i_1\in I$) and any 
$w_2\in(E_{i_2})_{x_0}$ ($i_2\in I\cup\{0\}$), we have 
$${}^{\widetilde v}\Gamma_{(\eta_{\widetilde v})_{\ast}w_1}w_2
=(\eta_{\widetilde v})_{\ast}(\Gamma_{w_1}w_2),$$
where we note that $T_{x_0}M
=T_{\eta_{\widetilde v}({x_0})}M_{\widetilde v}$ 
under the parallel translation in $V$.  
Also, we have $(\eta_{\widetilde v})_{\ast}w_1
=(1-(\lambda_{i_1})_{x_0}(\widetilde v_0))w_1$.}

\vspace{0.5truecm}

\noindent
{\it Proof.} From $(\eta_{\widetilde v})_{\ast{x_0}}
={\rm id}-A_{\widetilde v_0}$, 
the second relation follows directly, where $A$ is the shape tensor of $M$.   
Since $(\eta_{\widetilde v})_{\ast{x_0}}$ maps the $J$-curvature 
distributions of $M$ to those of $M_{\widetilde v}$, 
$\eta_{\widetilde v}$ maps the complex curvature spheres of 
$M$ through ${x_0}$ to those of $M_{\widetilde v}$ through 
$\eta_{\widetilde v}({x_0})$.  
On the other hand, since $F^{w_1}_t$ preserves $M$ inavariantly and its differential 
at a point of $M$ induces the parallel translation with respect to the normal connection 
of $M$, we have $\eta_{\widetilde v}\circ F^{w_1}_t
\vert_M=F^{w_1}_t\circ\eta_{\widetilde v}$.  
By using these facts and the properties of $F^{w_1}_t$, we can show that 
$F^{w_1}_t$ coincides with $F^{(\eta_{\widetilde v})_{\ast}w_1}_t$.  
From this fact, the first relation follows.  
\hspace{2.2truecm}q.e.d.

\vspace{0.5truecm}

We have the following fact for a principal orbit of an aks-representation 
of complex rank greater than one.  

\vspace{0.5truecm}

\noindent
{\bf Lemma 5.5.} {\sl Let $N$ be a principal orbit of an aks-representaion 
of complex rank greater than one, $\{n_i\,\vert\,i\in I\}$ the set of all 
$J$-curvature normals of $N$, $E_i$ the $J$-curvature distribution 
corresponding to $n_i$ and $\Gamma$ the extrinsically homogeneous structure of $N$ at 
$x$.  If the $2$-dimensional complex affine subspace $P$ through 
$n_{i_1},n_{i_2}$ and $n_{i_3}$ which does not pass through ${\bf 0}$, then, 
for any $w_k\in(E_{i_k})_x$ ($k=1,2,3$), we have 
$$\Gamma_{w_1}\Gamma_{w_2}w_3-\Gamma_{w_2}\Gamma_{w_1}w_3
=\Gamma_{(\Gamma_{w_1}w_2-\Gamma_{w_2}w_1)}w_3.$$}

\vspace{0.5truecm}

\noindent
{\it Proof.} 
Let $L/H$ be an irreducible anti-Kaehler symmetric space and $(\mathfrak l,\tau)$ the 
anti-Kaehler symmetric Lie algebra associated with $L/H$.  
We use the notations in Subsection 2.2.  Note that $I=\triangle_+\times\{0\}(=\triangle_+)$.  
Let $N$ be the principal orbit of the aks-representation 
$\rho:={\rm Ad}_L\vert_{\mathfrak p}:H\to GL(\mathfrak p)$ 
through a regular element $x(\in D)$.  
Take any $\alpha\in\triangle_+$ and any $w\in(E_{\alpha})_x(=\mathfrak p_{\alpha})$.  
Then, according to the proof of Lemma 4.6.3 of [K7], 
the holomorphic isometry $F^w_t$ is equal to $\rho(\exp_L(t\overline w))$, 
where $\overline w$ is 
where $\overline w$ is the element of $\mathfrak h_{\alpha}$ such that 
${\rm ad}_{\mathfrak l}(a)(\overline w)=w$ for all $a\in\mathfrak a$, where 
$\displaystyle{\mathfrak h_{\alpha}:=\{X\in \mathfrak h\,\vert\,{\rm ad}_{\mathfrak l}(a)^2(X)
=\alpha^{\Bbb C}(a)^2X\,\,{\rm for}\,\,{\rm all}\,\,a\in\mathfrak a\}}$.  
Hence we have 
$$\Gamma_w={\rm ad}_{\mathfrak l}(\overline w).\leqno{(5.1)}$$  
Therefore we can derive the desired relation in the method of the proof of Proposition 3.8 of [GH].  
\hspace{12.75truecm}q.e.d.

\vspace{0.5truecm}

For each $i\in I$, deonte by $W_i$ the complex affine subspace 
${x_0}+((E_i)_{x_0}\oplus{\rm Span}_{\bf C}\{(n_i)_{x_0}\})$ of $V$.  
Also, let $f_i$ be the focal map having $L^{E_i}_u$'s 
($u\in M$) as fibres, $\Phi_i$ the normal holonomy group 
of the focal submanifold $f_i(M)$ at $f_i({x_0})$ and 
$(\Phi_i)_{x_0}$ the isotropy group of $\Phi_i$ at ${x_0}$.  
This group $(\Phi_i)_{x_0}$ preserves $(E_i)_{x_0}$ invariantly.  
The irreducible decomposition 
of the action $(\Phi_i)_{x_0}\curvearrowright(E_i)_{x_0}$ 
is given by the form $(E_i)_{x_0}=(E_i)'_{x_0}\oplus (E_i)''_{x_0}$, 
where ${\rm dim}_{\bf C}(E_i)''_{x_0}=0,1\,\,{\rm or}\,\,3$, and 
${\rm dim}_{\bf C}(E_i)'_{x_0}$ is even in case of 
${\rm dim}_{\bf C}(E_i)''_{x_0}=1\,\,{\rm or}\,\,3$.  
Set $m_i:={\rm dim}_{\bf C}E_i$.  Note that $\Phi_i$ is orbit equivalent to 
the aks-representation associated with one of the following irreducible 
complex rank one anti-Kaehler symmetric spaces:
$$\begin{array}{c}
SO(m_i+2,{\Bbb C})/SO(m_i+1,{\Bbb C}),\,\,
SL(\frac{m_i+1}{2}+1,{\Bbb C})/SL(\frac{m_i+1}{2},{\Bbb C})\cdot
{\Bbb C}_{\ast},\\
Sp(\frac{m_i+1}{4}+1,{\Bbb C})/Sp(1,{\Bbb C})\times 
Sp(\frac{m_i+1}{4},{\Bbb C})
\end{array}$$
and that 
$${\rm dim}_{\bf C}(E_i)''_{x_0}=\left\{
\begin{array}{ll}
0 & ((\Phi_i)_{x_0}=SO(m_i+1,{\Bbb C}))\\
1 & ((\Phi_i)_{x_0}=SL(\frac{m_i+1}{2},{\Bbb C})\cdot{\Bbb C}_{\ast})\\
3 & ((\Phi_i)_{x_0}=Sp(1,{\Bbb C})\times Sp(\frac{m_i+1}{4},{\Bbb C})).
\end{array}\right.$$
By using Lemma 5.3 and $(5.1)$, we can derive the following fact corresponding to Proposition 3.11 of [GH].  

\vspace{0.5truecm}

\noindent
{\bf Lemma 5.6.} {\sl Let $i\in I$.  Then we have 
$$\begin{array}{c}
\Gamma_{(E_i)''_{x_0}}(E_i)''_{x_0}=0,\,\,\,\,
\Gamma_{(E_i)'_{x_0}}(E_i)''_{x_0}
\subset(E_i)'_{x_0},\\
\Gamma_{(E_i)''_{x_0}}(E_i)'_{x_0}\subset
(E_i)'_{x_0}\,\,\,\,{\rm and}\,\,\,\,\Gamma_{(E_i)'_{x_0}}(E_i)'_{x_0}
\subset(E_i)''_{x_0}.
\end{array}$$}

\vspace{0.5truecm}

Also, we have the following facts corresponding to Propositions 3.12 and 3.13 of [GH].  

\vspace{0.5truecm}

\noindent
{\bf Lemma 5.7.} {\sl For $i_1\in I$ and $i_2\in I\cup\{0\}$ with 
$i_2\not=i_1$, 
we have $\langle\Gamma_{(E_{i_1})_{x_0}}(E_{i_2})_{x_0},
(E_{i_2})_{x_0}\rangle=0$.}

\vspace{0.5truecm}

\noindent
{\bf Lemma 5.8.} {\sl Let $i_1\in I$ and $i_2,i_3\in I\cup\{0\}$.  
For $w_k\in(E_{i_k})_{x_0}$ ($k=1,2,3$), we have 
$(\overline{\nabla}_{w_1}\widetilde h)(w_2,w_3)=\langle\Gamma_{w_1}w_2,w_3
\rangle((n_{i_2})_{x_0}-(n_{i_3})_{x_0})$ and 
$\Gamma_{w_1}w_2=\widetilde{\nabla}_{w_1}\widetilde w_2\,\,\,\,
({\rm mod}\,(E_{i_2})_{x_0})$, where $\overline{\nabla}$ is the connection 
of the tensor bundle $T^{\ast}M\otimes T^{\ast}M\otimes T^{\perp}M$ induced 
from $\nabla$ and the normal connection $\nabla^{\perp}$ of 
$M$, and $\widetilde w_2$ is a local section of $E_{i_2}$ 
with $(\widetilde w_2)_{x_0}=w_2$.}

\vspace{0.5truecm}

Let $i_1,i_2,i_3\in I\cup\{0\}$ with $i_2\not=i_3$.  Then we define 
$\displaystyle{\frac{n_{i_1}-n_{i_3}}{n_{i_2}-n_{i_3}}}$ by 
$$\frac{n_{i_1}-n_{i_3}}{n_{i_2}-n_{i_3}}
:=\left\{
\begin{array}{ll}
b & \displaystyle{\left(\begin{array}{c}
\displaystyle{({\rm when}\,\,(n_{i_1})_{{x_0}}-(n_{i_3})_{{x_0}}=b((n_{i_2})_{{x_0}}-(n_{i_3})_{{x_0}})}\\
\displaystyle{{\rm for}\,\,{\rm some}\,\,b\in{\Bbb C}}
\end{array}\right)}\\
0 & \displaystyle{\left(\begin{array}{c}
\displaystyle{{\rm when}\,\,(n_{i_1})_{{x_0}}-(n_{i_3})_{{x_0}}\,\,{\rm and}\,\,
(n_{i_2})_{{x_0}}-(n_{i_3})_{{x_0}})}\\
\displaystyle{{\rm are}\,\,{\rm linearly}\,\,{\rm independent}\,\,{\rm over}\,\,{\Bbb C}}
\end{array}\right).}
\end{array}\right.$$
Note that this value is independent of the choice of ${x_0}\in M$.  
Denote by $w^k$ the $(E_k)_{x_0}$-component of $w\in T_{x_0}M$.  
We can derive the following fact corresponding to Proposition 3.15 of [GH] from the first relation in Lemma 5.8 
and the Codazzi equation.  

\vspace{0.5truecm}

\noindent
{\bf Lemma 5.9.} {\sl Let $i_1,i_2\in I$ and $i_3\in I\cup\{0\}$ with 
$i_3\not=i_2$.  
For any $w_k\in(E_{i_k})_{x_0}$ ($k=1,2$), we have 
$$(\Gamma_{w_1}w_2)^{i_3}=\frac{n_{i_1}-n_{i_3}}{n_{i_2}-n_{i_3}}
(\Gamma_{w_2}w_1)^{i_3}.$$}

\vspace{0.5truecm}

Also, we have the following fact correspondnig to Lemma 3.16 of [GH].  

\vspace{0.5truecm}

\noindent
{\bf Lemma 5.10.} {\sl {\rm (i)} Let $i_1\in I$ and $i_2,i_3\in I\cup\{0\}$.  
If $(\Gamma_{w_1}w_2)^{i_3}\not=0$ for some $w_1\in(E_{i_1})_{x_0}$ and 
$w_2\in(E_{i_2})_{x_0}$, then $(n_{i_1})_{x_0},(n_{i_2})_{x_0}$ and 
$(n_{i_3})_{x_0}$ are contained in a complex affine line.  

{\rm (ii)} Let $i_1,i_2,i_3\in I$.  The condition 
$(\Gamma_{(E_{i_1})_{x_0}}(E_{i_2})_{x_0})^{i_3}\not=0$ is symmetric in 
$i_1,i_2,i_3$.}

\vspace{0.5truecm}

Also, we have the following fact corresponding to Theorem 4.1 of [GH].  

\vspace{0.5truecm}

\noindent
{\bf Lemma 5.11.} {\sl $\sum\limits_{i_1,i_2\in I\,\,{\rm s.t.}\,\,i_1\not=i_2}
\Gamma_{(E_{i_1})_{x_0}}(E_{i_2})_{x_0}$ is dense in $T_{x_0}M$ and includes 
$\sum\limits_{i\in I}(E_i)_{x_0}$.}

\vspace{0.5truecm}

By using this lemma, we can derive the following fact corresponding to Corollary 4.2 of [GH].  

\vspace{0.5truecm}

\noindent
{\bf Lemma 5.12.} {\sl {\rm (i)} For each $i_1\in I$, we have 
$\sum\limits_{i_2,i_3\in I\,\,{\rm s.t.}\,\,n_{i_2},n_{i_3}\notin
{\rm Span}_{\bf C}\{n_{i_1}\}}(\Gamma_{(E_{i_2})_{x_0}}
(E_{i_3})_{x_0})^{i_1}=(E_{i_1})_{x_0}$.  

{\rm (ii)} 
$\sum\limits_{i_1,i_2\in I\,\,{\rm s.t.}\,\,n_{i_1},n_{i_2}\,:\,{\rm lin.}\,\,
{\rm dep.}}(\Gamma_{(E_{i_1})_{x_0}}(E_{i_2})_{x_0})^0$ is dense in 
$(E_0)_{x_0}$, where "lin. dep." means "linearly dependent".}

\vspace{0.5truecm}

\noindent
{\bf Notation.} In the sequel, for $w\in(E_i)_{x_0}$ ($i\in I\cup\{0\}$), 
$\widetilde w$ denotes a local section of $E_i$ with $\widetilde w_{x_0}=w$.  

\vspace{0.5truecm}

For $w_1\in(E_{i_1})_{x_0}$ and $w_2\in(E_{i_2})_{x_0}$ 
($i_1,i_2\in I\cup\{0\}$), define 
$\nabla'_{\widetilde w_1}\widetilde w_2$ by 
$(\nabla'_{\widetilde w_1}\widetilde w_2)_x
:=(\nabla_{\widetilde w_1}\widetilde w_2)_x-\Gamma^x_{(\widetilde w_1)_x}
(\widetilde w_2)_x$, where $x$ moves over the common domain of $\widetilde w_1$ 
and $\widetilde w_2$.  
Denote by $R$ the curvature tensor of $M$.  
Let $i_1,i_2,i_3\in I,\,\,i_4\in I\cup\{0\}$ and 
$w_k\in(E_{i_k})_{x_0}$ ($k=1,\cdots,4$).  According to the Gauss equation, 
we have 
$$\langle R(w_1,w_2)w_3,w_4\rangle
=(\langle w_1,w_4\rangle\langle w_2,w_3\rangle
-\langle w_1,w_3\rangle\langle w_2,w_4\rangle)\langle n_{i_1},n_{i_2}\rangle.
\leqno{(5.2)}$$
Also, from the definition of $\nabla'$, we have 
$$\begin{array}{l}
\langle R(w_1,w_2)w_3,w_4\rangle
=\langle\Gamma_{w_1}w_3,\Gamma_{w_2}w_4\rangle
-\langle\Gamma_{w_2}w_3,\Gamma_{w_1}w_4\rangle
-\langle(\nabla_{[\widetilde w_1,\widetilde w_2]}\widetilde w_3)_{x_0},
w_4\rangle\\
\hspace{1truecm}+w_1\langle(\nabla_{\widetilde w_2}\widetilde w_3)_{x_0},
w_4\rangle-\langle(\nabla'_{\widetilde w_2}\widetilde w_3)_{x_0},
(\nabla_{\widetilde w_1}\widetilde w_4)_{x_0}\rangle
-\langle\Gamma_{w_2}w_3,(\nabla'_{\widetilde w_1}\widetilde w_4)_{x_0}
\rangle\\
\hspace{1truecm}-w_2\langle(\nabla_{\widetilde w_1}\widetilde w_3)_{x_0},
w_4\rangle+\langle(\nabla'_{\widetilde w_1}\widetilde w_3)_{x_0},
(\nabla_{\widetilde w_2}\widetilde w_4)_{x_0}\rangle
+\langle\Gamma_{w_1}w_3,(\nabla'_{\widetilde w_2}\widetilde w_4)_{x_0}
\rangle.
\end{array}\leqno{(5.3)}$$

\vspace{0.5truecm}

For $\nabla'$ and $\Gamma$, we have the following relations.  

\vspace{0.5truecm}

\noindent
{\bf Lemma 5.13.} {\sl Let $i_1,i_2,i_3\in I$ and $i_4\in I\cup\{0\}$.  

{\rm (i)} For any $w_k\in(E_{i_k})_{x_0}$ ($k=1,2,3$), we have 
$$w_1\langle\widetilde w_2,\widetilde w_3\rangle
=\langle(\nabla'_{\widetilde w_1}\widetilde w_2)_{x_0},\,\widetilde w_3
\rangle+\langle w_2,\,(\nabla'_{\widetilde w_1}\widetilde w_3)_{x_0}
\rangle.$$

{\rm (ii)} If $i_1\not=i_2$, then we have 
$\nabla'_{\widetilde w_1}\widetilde w_2
=(\nabla_{\widetilde w_1}\widetilde w_2)^{i_2}$ for any 
$w_k\in(E_{i_k})_{x_0}$ ($k=1,2$).  

{\rm (iii)} For any $w_k\in(E_{i_k})_{x_0}$ ($k=1,2,3$), we have 
$$\left(\nabla'_{\widetilde w_1}(\widetilde{(\Gamma_{w_2}w_3)^{i_3})})
\right)_{x_0}
=\left(\Gamma_{(\nabla'_{\widetilde w_1}\widetilde w_2)_{x_0}}w_3
\right)^{i_3}
+\left(\Gamma_{w_2}(\nabla'_{\widetilde w_1}\widetilde w_3)_{x_0}
\right)^{i_3}.$$}

\vspace{0.5truecm}

\noindent
{\it Proof.} The relations in (i) and (ii) are trivial.  From (ii) of Lemma 5.2, the relation in (iii) is shown 
in the method of the proof of Lemma 5.2 of [GH].  
\hspace{2.75truecm}q.e.d.

\vspace{0.5truecm}

Let $i_1\in I$ and $i_2\in I\cup\{0\}$.  
For $w\in T_{x_0}M,\,w_1\in(E_{i_1})_{x_0}$ and 
$w_2\in(E_{i_2})_{x_0}$, we define 
$\langle\Gamma_ww_1,w_2\rangle$ by 
$$\begin{array}{l}
\displaystyle{\langle\Gamma_ww_1,w_2\rangle:=-\sum_{i\in I}
\langle\Gamma_{w_1}w_2,\frac{n_i-n_{i_2}}{n_{i_1}-n_{i_2}}w^i\rangle}\\
\displaystyle{\left(=\lim_{m\to\infty}\sum_{i\in I\,\,{\rm s.t.}\,\,\vert w^i\vert>\frac 1m}
\langle\Gamma_{w_1}w_2,\frac{n_i-n_{i_2}}{n_{i_1}-n_{i_2}}w^i\rangle\right).}
\end{array}
\leqno{(5.4)}$$
According to (i) of Lemma 5.2 and Lemma 5.9, this defintion is valid.  
From the relation in (iii) of Lemma 5.13, we can show the following fact in the method of the proof of 
Theorem 5.7 of [GH].  

\vspace{0.5truecm}

\noindent
{\bf Lemma 5.14.} {\sl Let $i_1,i_2,i_3\in I$ and $i_4\in I\cup\{0\}$ with 
$i_4\not=i_3$.  For any $w_k\in(E_{i_k})_{x_0}$ ($k=1,\cdots,4$), we have 
$$\langle\left([\Gamma_{w_1},\Gamma_{w_2}]
-\Gamma_{\Gamma_{w_1}w_2-\Gamma_{w_2}w_1}\right)w_3,w_4\rangle
=-(\langle w_1,w_4\rangle\langle w_2,w_3\rangle
-\langle w_1,w_3\rangle\langle w_2,w_4\rangle)
\langle n_{i_1},n_{i_2}\rangle.$$}

\vspace{0.5truecm}

\noindent
By using Lemmas 5.9 and 5.14, we can show the following fact.  

\vspace{0.5truecm}

\noindent
{\bf Lemma 5.15.} {\sl Let $(i_1,i_2,i_3)$ be an element of 
$I^2\times(I\cup\{0\})$ such that 
there exists no complex affine line containing 
$(n_{i_1})_{x_0},(n_{i_2})_{x_0}$ and $(n_{i_3})_{x_0}$, and $i_4$ 
an element of $I$.  For any $w_k\in(E_{i_k})_{x_0}$ ($k=1,\cdots,4$), 
we have 
$$\begin{array}{l}
\langle\Gamma_{w_1}w_2,\Gamma_{w_4}w_3\rangle
=\langle\Gamma_{w_4}w_2,\Gamma_{w_1}w_3\rangle
+c\langle\Gamma_{w_1}w_4,\Gamma_{w_2}w_3\rangle,
\end{array}$$
where $c$ is a constant.  
Furthermore, if $i_1=i_4$ or the intersection of the complex affine line 
through $(n_{i_1})_{x_0}$ and $(n_{i_4})_{x_0}$ and 
the complex affine line through and $(n_{i_2})_{x_0}$ and $(n_{i_3})_{x_0}$ 
contains no $J$-curvature normal, then we have $c=0$.  On the other hand, 
if their intersection contains a $J$-curvature normal $(n_{i_5})_{x_0}$, then we have 
$$c=\frac{n_{i_3}-n_{i_5}}{n_{i_2}-n_{i_3}}\times
\frac{n_{i_1}-n_{i_4}}{n_{i_1}-n_{i_5}}.$$
}

\vspace{0.5truecm}

\noindent
We can show the following fact in the method of the proof of Corollary 5.11 of [GH].  

\vspace{0.5truecm}

\noindent
{\bf Lemma 5.16.} {\sl Let $i_1,i_2,i_3\in I$ satisfying $i_3\not=i_1,i_2$ and 
$\frac{n_{i_2}}{n_{i_3}}\not=-\frac{n_{i_1}-n_{i_2}}{n_{i_1}-n_{i_3}}$.  
Assume that $\langle(\Gamma_{(E_{i_1})_{x_0}}(E_{i_2})_{x_0})^{i_4},\,
\Gamma_{(E_{i_1})_{x_0}}(E_{i_3})_{x_0}\rangle=0$ for any $i_4\in I$ 
and $(\Gamma_{(E_{i_1})_{x_0}}(E_{i_2})_{x_0})^{i_3}=0$ (these 
conditions hold 
if $\Gamma_{(E_{i_1})_{x_0}}(E_{i_2})_{x_0}\subset(E_0)_{x_0}$).  
Then we have $\langle\Gamma_{(E_{i_1})_{x_0}}(E_{i_2})_{x_0},\,
\Gamma_{(E_{i_1})_{x_0}}(E_{i_3})_{x_0}\rangle$\newline
$=0$.}

\vspace{0.5truecm}

\noindent
Also, we can derive the following fact.  

\vspace{0.5truecm}

\noindent
{\bf Lemma 5.17.} {\sl Let $i_1,i_2\in I$ with $i_1\not=i_2$.  
For any $w_k\in(E_{i_k})_{x_0}$ ($k=1,2$), we have 
$$\sum_{i_3\in(I\cup\{0\})\setminus\{i_1\}}
{\rm Re}\left(\frac{n_{i_2}-n_{i_3}}{n_{i_1}-n_{i_3}}\right)
\vert\vert(\Gamma_{w_1}w_2)^{i_3}\vert\vert^2
=\frac12\langle n_{i_1},n_{i_2}\rangle\,\langle w_1,w_1\rangle\,
\vert\vert w_2\vert\vert^2.$$}

\vspace{0.5truecm}

\noindent
{\it Proof.} 
Let $w_2=(w_2)_-+(w_2)_+\,\,\,\,((w_2)_-\in((E_{i_2})_-)_{x_0},\,\,
(w_2)_+\in((E_{i_2})_+)_{x_0})$.  
In similar to Corollary 5.13 of [GH], we can show 
$$\begin{array}{l}
\hspace{0truecm}\displaystyle{\sum_{i_3\in(I\cup\{0\})\setminus\{i_1\}}
\left\langle(\Gamma_{w_1}(w_2)_{\varepsilon})^{i_3},\frac{n_{i_2}-n_{i_3}}
{n_{i_1}-n_{i_3}}(\Gamma_{w_1}(w_2)_{\varepsilon})^{i_3}\right\rangle}\\
\displaystyle{=\frac12\langle n_{i_1},n_{i_2}\rangle\,\langle w_1,w_1\rangle
\,\langle(w_2)_{\varepsilon},(w_2)_{\varepsilon}\rangle,}
\end{array}\leqno{(5.5)}$$
where $\varepsilon=-\,\,{\rm or}\,\,+$.  
On the other hand, since $F^{w_1}_t$'s preserve $E_i$'s invariantly and they are holomorphic 
isometries, $\Gamma_{w_1}$ preserves $((E_i)_-)_{x_0}$'s and $((E_i)_+)_{x_0}$ invariantly, 
respectively.  Hence we have 
$\Gamma_{w_1}(w_2)_{\varepsilon}=(\Gamma_{w_1}w_2)_{\varepsilon}$.  
Also, it is clear that 
$((\Gamma_{w_1}w_2)_{\varepsilon})^{i_3}=((\Gamma_{w_1}w_2)^{i_3})_{\varepsilon}$.  
From these relations, we have 
$$\left\langle(\Gamma_{w_1}(w_2)_{\varepsilon})^{i_3},
\frac{n_{i_2}-n_{i_3}}{n_{i_1}-n_{i_3}}
(\Gamma_{w_1}(w_2)_{\varepsilon})^{i_3}\right\rangle
={\rm Re}\left(\frac{n_{i_2}-n_{i_3}}{n_{i_1}-n_{i_3}}\right)
\langle((\Gamma_{w_1}w_2)^{i_3})_{\varepsilon},\,
((\Gamma_{w_1}w_2)^{i_3})_{\varepsilon}\rangle.$$
By summing the $(-1)$-multiples of $(5.5)$'s for $\varepsilon=\pm$ and using this relation, we have 
the desired relation.  
\begin{flushright}q.e.d.\end{flushright}

\vspace{0.5truecm}

By using Lemmas 5.3, 5.7, 5.10 and 5.17, we can show the following fact.  

\vspace{0.5truecm}

\noindent
{\bf Lemma 5.18.} {\sl Assume that the complex Coxeter group ${\cal W}$ associated with $M$ 
is of type $\widetilde A,\,\widetilde D$ or $\widetilde E$.  
Let $i_1$ and $i_2$ be elements of $I$ such that $n_{i_1}$ and $n_{i_2}$ are linearly independent.  

{\rm (i)} If $n_{i_1}$ and $n_{i_2}$ are orthogonal, then we have 
$\Gamma_{w_1}w_2=0$ for any $w_k\in(E_{i_k})_{x_0}$ ($k=1,2$).  

{\rm (ii)} If $n_{i_1}$ and $n_{i_2}$ are not orthogonal, then we have 
$\displaystyle{\vert\vert\Gamma_{w_1}w_2\vert\vert\leq\frac12
\vert\vert w_1\vert\vert\,\vert\vert w_2\vert\vert\,
\vert\vert n_{i_1}\vert\vert}$ 
for any $w_k\in(E_{i_k})_{x_0}$ ($k=1,2$).}

\vspace{0.5truecm}

\noindent
{\it Proof.} Let $P$ be the complex affine line in 
$\mathfrak b$ through $(n_{i_1})_{x_0}$ and 
$(n_{i_2})_{x_0}$.  Since $n_{i_1}$ and $n_{i_2}$ are linearly independent, 
we have $0\notin P$.  Hence the slice $L^P_{x_0}$ is a finite dimensional 
anti-Kaehler isoparametric submanifold with $J$-diagonalizable shape operators 
(of codimension two in $(W_P)_{x_0}$).  
Hence, since ${\cal W}$ is isomorphic to an affine Weyl group 
of type $\widetilde A,\,\widetilde D$ or $\widetilde E$, 
the root system (which we denote by $\triangle_P$) associated with $L^P_{x_0}$ is of 
type $A_1\times A_1$ or $A_2$.  
First we shall show the statement (i).  
Assume that $(n_{i_1})_{x_0}$ and $(n_{i_2})_{x_0}$ are orthogonal.  
Then $\triangle_P$ is of type $A_1\times A_1$ and hence $P$ contains no other 
$J$-curvature normal.  
By using this fact and Lemma 5.3, we can show 
$\Gamma_{w_1}w_2={}^{L^P_{x_0}}\Gamma_{w_1}w_2=0$ for any 
$w_k\in(E_{i_k})_{x_0}$ ($k=1,2$), where ${}^{L^P_{x_0}}\Gamma$ is the 
extrinsically homogeneous structure of $L^P_{x_0}$.  
Next we shall show the statement (ii).  
Assume that $(n_{i_1})_{x_0}$ and $(n_{i_2})_{x_0}$ are not orthogonal.  
Then $\triangle_P$ is of type $A_2$ and hence 
there exists $i_3\in I\setminus\{i_1,i_2\}$ with $(n_{i_3})_{x_0}\in P$.  
The set 
${\it l}_{i_1}\cap{\it l}_{i_2}\cap{\it l}_{i_3}\cap
{\rm Span}_{\bf C}\{(n_{i_1})_{x_0},(n_{i_2})_{x_0}\}$ 
consists of the only one point.  Denote by $p_0$ this point.  
Let $e_1,e_2$ and $e_3$ be a unit normal vector of 
${\it l}_{i_1},{\it l}_{i_2}$ and ${\it l}_{i_3}$, respectively.  
Since $\triangle_P$ is of type $(A_2)$, 
we may assume that $e_3=e_1+e_2$ by replacing some of these vectors 
to the $(-1)$-multiples of them if necessary.  
Since $\frac{(n_{i_1})_{x_0}}{\langle(n_{i_1})_{x_0},(n_{i_1})_{x_0}
\rangle}\in{\it l}_{i_1}$, we have $(n_{i_1})_{x_0}
=\frac{e_1}{\langle\overrightarrow{0p_0},e_1\rangle}$, 
where $0$ is the origin of $\mathfrak b$.  
Similarly we have $(n_{i_2})_{x_0}=\frac{e_2}{\langle\overrightarrow{0p_0},
e_2\rangle}$ 
and 
$(n_{i_3})_{x_0}=\frac{e_3}{\langle\overrightarrow{0p_0},e_3\rangle}$.  
By using these facts, Lemmas 5.7, 5.10 and 5.17, we can show 
$$\begin{array}{l}
\displaystyle{\vert\vert\Gamma_{w_1}w_2\vert\vert^2
=\vert\vert(\Gamma_{w_1}w_2)^{i_3}\vert\vert^2
\leq\frac12{\rm Re}\left(\frac{n_{i_1}-n_{i_3}}{n_{i_2}-n_{i_3}}\right)
\vert\langle n_{i_1},n_{i_2}\rangle\vert\,\vert\vert w_1\vert\vert^2\,
\vert\vert w_2\vert\vert^2}\\
\hspace{4.45truecm}\displaystyle{\leq\frac14
\vert\vert w_1\vert\vert^2\,\vert\vert w_2\vert\vert^2\,
\vert\vert n_{i_1}\vert\vert^2.}
\end{array}$$
Thus we obtain the desired relation.  
\begin{flushright}q.e.d.\end{flushright}

\vspace{0.5truecm}

By using Lemmas 5.3, 5.4, 5.7 and 5.10, we can show the following fact.  

\vspace{0.5truecm}

\noindent
{\bf Lemma 5.19.} {\sl We have 
$$\mathop{\sup}_{i\in I}\,\,\mathop{\sup}_{P\in{\cal H}_i}
\,\,\mathop{\sup}_{(w_1,w_2)\in(E_i)_{x_0}\times(D_P)_{x_0}}
\frac{\vert\vert\Gamma_{w_1}w_2\vert\vert}{\vert\vert w_1\vert\vert\,
\vert\vert w_2\vert\vert\,\vert\vert(n_i)_{x_0}\vert\vert}\,\,<\,\,\infty,$$
where ${\cal H}_i$ is the set of all complex affine subspaces $P$ in 
$T_{x_0}M$ with $0\notin P$ and $(n_i)_{x_0}\in P$.}

\vspace{0.5truecm}

\noindent
{\it Proof.} Let ${\cal H}_i^{\rm irr}$ be the set of all elements $P$ of 
${\cal H}_i$ such that $L^P_{x_0}(\subset(W_P)_{x_0})$ is irreducible.  
First we shall show 
$$\mathop{\sup}_{i\in I}\,\,\mathop{\sup}_{P\in {\cal H}_i^{\rm irr}}
\,\,\mathop{\sup}_{(w_1,w_2)\in(E_i)_{x_0}\times(D_P)_{x_0}}
\frac{\vert\vert\Gamma_{w_1}w_2\vert\vert}{\vert\vert w_1\vert\vert\,
\vert\vert w_2\vert\vert\,\vert\vert(n_i)_{x_0}\vert\vert}\,\,<\,\,\infty.
\leqno{(5.6)}$$
Fix $i_0\in I$ and $P_0\in{\cal H}_{i_0}^{\rm irr}$.  
If the complex codimension of $L^{P_0}_{x_0}(\subset(W_{P_0})_{x_0})$ is 
equal to one, then we can take $P'_0\in{\cal H}_{i_0}^{\rm irr}$ such that 
$P_0\subset P'_0$ and that the complex codimension of 
$L^{P'_0}_{x_0}(\subset(W_{P'_0})_{x_0})$ is greater than one.  
Then we have 
$$\begin{array}{l}
\hspace{0.4truecm}
\displaystyle{\mathop{\sup}_{(w_1,w_2)\in(E_{i_0})_{x_0}\times
(D_{P_0})_{x_0}}
\frac{\vert\vert\Gamma_{w_1}w_2\vert\vert}{\vert\vert w_1\vert\vert\,
\vert\vert w_2\vert\vert\,\vert\vert(n_{i_0})_{x_0}\vert\vert}}\\
\displaystyle{\leq\,\mathop{\sup}_{(w_1,w_2)\in(E_{i_0})_{x_0}\times
(D_{P'_0})_{x_0}}
\frac{\vert\vert\Gamma_{w_1}w_2\vert\vert}{\vert\vert w_1\vert\vert\,
\vert\vert w_2\vert\vert\,\vert\vert(n_{i_0})_{x_0}\vert\vert}.}
\end{array}$$
and hence 
$$\begin{array}{l}
\hspace{0.5truecm}
\displaystyle{\mathop{\sup}_{i\in I}\,
\mathop{\sup}_{P\in{\cal H}_i^{\rm irr}}
\,\,\mathop{\sup}_{(w_1,w_2)\in(E_i)_{x_0}\times(D_P)_{x_0}}
\frac{\vert\vert\Gamma_{w_1}w_2\vert\vert}{\vert\vert w_1\vert\vert\,
\vert\vert w_2\vert\vert\,\vert\vert(n_i)_{x_0}\vert\vert}}\\
\displaystyle{=\mathop{\sup}_{i\in I}\,\,\mathop{\sup}_{P\in{\cal H}_i^{{\rm irr},\geq2}}
\,\,\mathop{\sup}_{(w_1,w_2)\in(E_i)_{x_0}\times(D_P)_{x_0}}
\frac{\vert\vert\Gamma_{w_1}w_2\vert\vert}{\vert\vert w_1\vert\vert\,
\vert\vert w_2\vert\vert\,\vert\vert(n_i)_{x_0}\vert\vert},}
\end{array}\leqno{(5.7)}$$
where ${\cal H}_i^{{\rm irr},\geq2}$ is the set of all elements $P$'s 
of ${\cal H}_i^{{\rm irr}}$ such that the complex codimension of 
$L^P_{x_0}(\subset(W_P)_{x_0})$ is greater than one.  
Fix $\alpha_1\in(\triangle_M)_+$ and $P_1\in{\cal H}_{(\alpha_1,0)}^{{\rm irr},\geq2}$.  
Take any $j_1\in{\Bbb Z}$.  For each $P\in{\cal H}_{(\alpha_1,0)}^{\rm irr}$, 
there exists $P'\in{\cal H}_{(\alpha_1,j_1)}^{\rm irr}$ such that 
$\{\alpha\in (\triangle_M)_+\,\vert\,\exists\,j\in{\Bbb Z}\,\,{\rm s.t.}\,\,(n_{(\alpha,j)})_{x_0}\in P\}
=\{\alpha\in (\triangle_M)_+\,\vert\,\exists\,j\in{\Bbb Z}\,\,{\rm s.t.}\,\,(n_{(\alpha,j)})_{x_0}\in P'\}$.  
Then, since ${\rm dim}_{\Bbb C}(W_P)_{x_0}={\rm dim}_{\Bbb C}(W_{P'})_{x_0}$, and since 
the root systems associated with $L^{P}_{x_0}$ and $L^{P'}_{x_0}$ coincide, they 
are ragraded as principal orbits of the aks-representation of the same irreducible anti-Kaehler 
symmetric space.  That is, $L^{P'}_{x_0}$ is regarded as 
a parallel submanifold of $L^P_{x_0}$ under a suitable identification of $(W_P)_{x_0}$ and 
$(W_{P'})_{x_0}$.  Therefore, by using Lemmas 5.3 and 5.4, we can show 
$$\begin{array}{l}
\hspace{0.5truecm}\displaystyle{
\mathop{\sup}_{P\in{\cal H}_{(\alpha_1,0)}^{\rm irr}}
\,\,\mathop{\sup}_{(w_1,w_2)\in(E_{(\alpha_1,0)})_{x_0}\times(D_P)_{x_0}}
\frac{\vert\vert\Gamma_{w_1}w_2\vert\vert}{\vert\vert w_1\vert\vert\,
\vert\vert w_2\vert\vert\,\vert\vert(n_{(\alpha_1,0)})_{x_0}\vert\vert}}\\
\displaystyle{
=\mathop{\sup}_{P\in{\cal H}_{(\alpha_1,j_1)}^{\rm irr}}
\,\,\mathop{\sup}_{(w_1,w_2)\in(E_{(\alpha_1,j_1)})_{x_0}\times(D_P)_{x_0}}
\frac{\vert\vert\Gamma_{w_1}w_2\vert\vert}{\vert\vert w_1\vert\vert\,
\vert\vert w_2\vert\vert\,\vert\vert(n_{(\alpha_1,j_1)})_{x_0}\vert\vert}.}
\end{array}$$
Hence it follows from the arbitrariness of $j_1$ that 
$$\begin{array}{l}
\hspace{0.5truecm}\displaystyle{
\mathop{\sup}_{i\in I}\,\,\mathop{\sup}_{P\in{\cal H}_i^{\rm irr}}
\,\,\mathop{\sup}_{(w_1,w_2)\in(E_i)_{x_0}\times(D_P)_{x_0}}
\frac{\vert\vert\Gamma_{w_1}w_2\vert\vert}{\vert\vert w_1\vert\vert\,
\vert\vert w_2\vert\vert\,\vert\vert(n_i)_{x_0}\vert\vert}}\\
\displaystyle{=\mathop{\sup}_{\alpha\in (\triangle_M)_+}\,\,
\mathop{\sup}_{P\in{\cal H}_{(\alpha,0)}^{\rm irr}}
\,\,\mathop{\sup}_{(w_1,w_2)\in(E_{(\alpha,0)})_{x_0}\times(D_P)_{x_0}}
\frac{\vert\vert\Gamma_{w_1}w_2\vert\vert}{\vert\vert w_1\vert\vert\,
\vert\vert w_2\vert\vert\,\vert\vert(n_{(\alpha,0)})_{x_0}\vert\vert}
\,\,<\,\,\infty.}
\end{array}$$
Thus we obtain $(5.6)$.  
For simplicity set 
$$C:=\mathop{\sup}_{i\in I}\,\,\mathop{\sup}_{P\in {\cal H}_i^{\rm irr}}
\,\,\mathop{\sup}_{(w_1,w_2)\in(E_i)_{x_0}\times(D_P)_{x_0}}
\frac{\vert\vert\Gamma_{w_1}w_2\vert\vert}{\vert\vert w_1\vert\vert\,
\vert\vert w_2\vert\vert\,\vert\vert(n_i)_{x_0}\vert\vert}.$$
Fix $i_0\in I$ and $P_0\in{\cal H}_{i_0}\setminus{\cal H}_{i_0}^{\rm irr}$.  
Let $L^{D_{P_0}}_{x_0}=L_1\times\cdots\times L_k$ be the irreducible 
decomposition of $L^{D_{P_0}}_{x_0}$.  
Take any $i_1,i_2\in I$ with $(n_{i_1})_{x_0},(n_{i_2})_{x_0}\in P_0$.  
If $(n_{i_1})_{x_0}$ and $(n_{i_2})_{x_0}$ are not orthogonal, then 
$(E_{i_1})_{x_0}\oplus(E_{i_2})_{x_0}\subset T_{x_0}L_a$ for some 
$a\in\{1,\cdots,k\}$.  Hence we have 
$$\mathop{\sup}_{(w_1,w_2)\in(E_{i_1})_{x_0}\times(E_{i_2})_{x_0}}
\frac{\vert\vert\Gamma_{w_1}w_2\vert\vert}{\vert\vert w_1\vert\vert\,
\vert\vert w_2\vert\vert\,\vert\vert(n_{i_1})_{x_0}\vert\vert}\leq C.$$
If $(n_{i_1})_{x_0}$ and $(n_{i_2})_{x_0}$ are orthogonal, then the 
complex affine line through $(n_{i_1})_{x_0}$ and $(n_{i_2})_{x_0}$ 
does not contain other $J$-curvature normal.  Hence it follows from Lemma 5.7 and (i) of Lemma 5.10 
that $\Gamma_{(E_{i_1})_{x_0}}(E_{i_2})_{x_0}=0$.  
Therefore, we obtain 
$$\mathop{\sup}_{i\in I}\,\,\mathop{\sup}_{P\in {\cal H}_i}
\,\,\mathop{\sup}_{(w_1,w_2)\in(E_i)_{x_0}\times(D_P)_{x_0}}
\frac{\vert\vert\Gamma_{w_1}w_2\vert\vert}{\vert\vert w_1\vert\vert\,
\vert\vert w_2\vert\vert\,\vert\vert(n_i)_{x_0}\vert\vert}=C.$$
This completes the proof.  
\begin{flushright}q.e.d.\end{flushright}

\vspace{0.5truecm}

By using Lemma 5.19, we can show the following fact.  

\vspace{0.5truecm}

\noindent
{\bf Lemma 5.20.} {\sl Let $i_0=(\alpha_0,j_0)\in I$ and 
$w\in(E_{i_0})_{x_0}$.  
Then $\Gamma_w$ can be extended continuously to $T_{x_0}
M$ 
if and only if the restriction of $\Gamma_w$ to 
$\displaystyle{\mathop{\oplus}_{j\in{\Bbb Z}}(E_{(\alpha_0,j)})_{x_0}}$ 
can be extended continuously to 
$\displaystyle{\overline{\mathop{\oplus}_{j\in{\Bbb Z}}
(E_{(\alpha_0,j)})_{x_0}}}$.}

\vspace{0.5truecm}

\noindent
{\it Proof.} Set $V_0:=(E_0)_{x_0}$, 
$\displaystyle{V_1:=\mathop{\oplus}_
{i\in I\setminus\{(\alpha_0,j)\,\vert\,j\in{\Bbb Z}\}}(E_i)_{x_0}}$ 
and 
$\displaystyle{V_2:=\mathop{\oplus}_
{j\in{\Bbb Z}}(E_{(\alpha_0,j)})_{x_0}}$.  
Clearly we have $T_{x_0}M=V_0\oplus\overline V_1\oplus 
\overline V_2$.  
Since $\Gamma_w$ is a closed operator by the definition and since 
$(E_0)_{x_0}$ is closed in the domain of $\Gamma_w$, 
$\Gamma_w\vert_{(E_0)_{x_0}}$ also is a closed operator.  
Hence, according to the closed graph theorem, $\Gamma_w\vert_{(E_0)_{x_0}}$ is bounded 
(hence continuous).  Easily we can show 
$$V_1=\mathop{\oplus}_{{\it l}}\left(\mathop{\oplus}_
{i\in I\setminus\{(\alpha_0,j)\,\vert\,j\in{\Bbb Z}\}\,\,{\rm s.t.}\,\,
(n_i)_{x_0}\in{\it l}}(E_i)_{x_0}\right),$$
where ${\it l}$ runs over the set of all complex affine lines in 
$\mathfrak b\setminus\{0\}$ through 
$(n_{i_0})_{x_0}$.  
For simplicity set 
$$V_{1,{\it l}}:=
\mathop{\oplus}_{i\in I\setminus\{(\alpha_0,j)\,\vert\,j\in{\Bbb Z}\}\,\,{\rm s.t.}\,\,
(n_i)_{x_0}\in{\it l}}(E_i)_{x_0}.$$
According to Lemma 5.19, for each ${\it l}$, 
we have 
$$\mathop{\sup}_{w'\in V_{1,{\it l}}}
\frac{\vert\vert\Gamma_ww'\vert\vert}{\vert\vert w'\vert\vert}\leq 
C\vert\vert(n_{i_0})_{x_0}\vert\vert\,\vert\vert w\vert\vert,$$
where $C$ is the positive constant as in the proof of Lemma 5.19, 
and hence 
$$\mathop{\sup}_{w'\in V_1}\frac{\vert\vert\Gamma_ww'\vert\vert}
{\vert\vert w'\vert\vert}\leq C\vert\vert(n_{i_0})_{x_0}\vert\vert\,
\vert\vert w\vert\vert.$$
Therefore the restriction of $\Gamma_w$ to $V_1$ is bounded and hence it can be extended continuously 
to $\overline V_1$.  
From these facts, the statement of this lemma follows.  
\begin{flushright}q.e.d.\end{flushright}

\vspace{0.5truecm}

According to Lemma 6.4 of [GH], we have the following fact.  

\vspace{0.5truecm}

\noindent
{\bf Lemma 5.21.} {\sl Let $W$ be a Hilbert space, 
$\displaystyle{W=\overline{\mathop{\oplus}_{i\in{\Bbb Z}}W_i}}$ 
the orthogonal decomposition of 
$W$ and $f$ a linear map from 
$\displaystyle{\mathop{\oplus}_{i\in{\Bbb Z}}W_i}$ to $W$.  
Assume that there exists a positive constant $C$ such that 
$\vert\vert f(w)\vert\vert\leq C\vert\vert w\vert\vert$ 
for all $\displaystyle{w\in\mathop{\cup}_{i\in{\Bbb Z}}W_i}$ and that 
there exist injective maps $\mu_a:{\Bbb Z}\to{\Bbb Z}$ ($a=1,\cdots,r$) 
such that $\langle f(W_i),f(W_j)\rangle=0$ for any 
$j\notin\{\mu_1(i),\cdots,\mu_r(i)\}$.  Then we have 
$\vert\vert f\vert\vert\leq \sqrt rC$ and hence $f$ can be extended 
continuously to $W$.}

\vspace{0.5truecm}

Easily we can show that 
$$\frac{n_{(\alpha,j_1)}-n_{(\alpha,j_3)}}{n_{(\alpha,j_2)}-n_{(\alpha,j_3)}}
=\frac{j_1-j_3}{j_2-j_3}\times\frac{1+j_2b_{\alpha}{\bf i}}
{1+j_1b_{\alpha}{\bf i}}.\leqno{(5.8)}$$
By using $(5.8)$ and Lemma 5.17, we can show the following fact.  

\vspace{0.5truecm}

\noindent
{\bf Lemma 5.22.} {\sl Let $\alpha\in (\triangle_M)_+$ and $j_1,j_2\in{\Bbb Z}$.  
For any $w_1\in(E_{(\alpha,j_1)})_{x_0}$ and any 
$w_2\in(E_{(\alpha,j_2)})_{x_0}$, we have 
$$\begin{array}{l}
\hspace{0.5truecm}\displaystyle{
\sum_{j\in{\bf Z}\setminus\{j_1\}}\frac{j-j_2}{j-j_1}
\vert\vert(\Gamma_{w_1}w_2)^{(\alpha,j)}\vert\vert^2
+\vert\vert(\Gamma_{w_1}w_2)^0\vert\vert^2}\\
\displaystyle{=\frac12\left({\rm Re}\left(\frac{1+j_1b_{\alpha}{\bf i}}
{1+j_2b_{\alpha}{\bf i}}\right)\right)^{-1}\,\,
{\rm Re}\left(\frac{1}{(1+j_1b_{\alpha}{\bf i})(1+j_2b_{\alpha}{\bf i})}\right)
\,\langle(n_{(\alpha,0)})_{x_0},(n_{(\alpha,0)})_{x_0}\rangle\,
\langle w_1,w_1\rangle\,\,\vert\vert w_2\vert\vert^2.}
\end{array}$$}

\vspace{0.5truecm}

Also, we can show the following fact.  

\vspace{0.5truecm}

\noindent
{\bf Lemma 5.23.} {\sl Let $P$ be the complex affine line through ${\bf 0}$ and 
$(n_{(\alpha_0,0)})_{x_0}$ for some $\alpha_0\in(\triangle_M)_+$.  

{\rm (i)} If the affine root system ${\cal R}$ is of type 
$(\widetilde A_m)$ ($m\geq 2$), $(\widetilde D_m)$ ($m\geq 4$), 
$(\widetilde E_m)$ ($m=6,7,8$) or $(\widetilde F_4)$, then 
there exists a (complex) $2$-dimensional complex affine 
subspace $P'$ including $P$ such that the affine root system 
associated with $L^{P'}_{x_0}(\subset(W_{P'})_{x_0})$ is of type 
($\widetilde A_2$).  

{\rm (ii)} If the affine root system ${\cal R}$ is of type 
$(\widetilde B_m),\,(\widetilde B_m^v)$ or $(\widetilde B_m,\widetilde B_m^v)$ 
($m\geq 2$), then there exists a (complex) $2$-dimensional complex affine 
subspace $P'$ including $P$ such that the affine root system 
associated with $L^{P'}_{x_0}(\subset(W_{P'})_{x_0})$ is of type 
"($\widetilde A_2$) or $(\widetilde C_2)$", 
"($\widetilde A_2$) or $(\widetilde C_2^v)$" or 
"($\widetilde A_2$) or $(\widetilde C_2,\widetilde C_2^v)$", respectively.

{\rm (iii)} If the affine root system ${\cal R}$ is of type 
$(\widetilde C_m),\,(\widetilde C_m^v),\,(\widetilde C'_m),\,
(\widetilde C_m^v,\widetilde C'_m),\,(\widetilde C'_m,\widetilde C_m),$\newline
$(\widetilde C_m^v,\widetilde C_m)$ or $(\widetilde C_m,\widetilde C_m^v)$ 
($m\geq 2$), then there exists a (complex) $2$-dimensional complex affine 
subspace $P'$ including $P$ such that the affine root system 
associated with $L^{P'}_{x_0}(\subset(W_{P'})_{x_0})$ is of type 
"($\widetilde A_2$) or $(\widetilde C_2)$", 
"($\widetilde A_2$) or $(\widetilde C_2^v)$", 
"($\widetilde A_2$) or $(\widetilde C'_2)$", 
"($\widetilde A_2$) or $(\widetilde C_2^v,\widetilde C'_2)$", 
"($\widetilde A_2$) or $(\widetilde C'_2,\widetilde C_2)$", 
"($\widetilde A_2$) or $(\widetilde C^v_2,\widetilde C_2)$" or 
"($\widetilde A_2$) or $(\widetilde C_2,\widetilde C_2^v)$", respectively.  
}

\vspace{0.5truecm}

\noindent
{\it Proof.} First we shall show the statement (i).  
Let $\Pi(\subset(\triangle_M)_+)$ be a simple root system of $\triangle_M$.  
Without loss of generality, we may assume that $\alpha_0$ is one of the elements of $\Pi$.  
Since ${\cal R}$ is of $(\widetilde A_m)$ ($m\geq 2$), 
$(\widetilde D_m)$ ($m\geq 4$), $(\widetilde E_m)$ ($m=6,7,8$) or 
$(\widetilde F_4)$, it follows from their Dynkin diagrams that there exists 
$\alpha_1\in\Pi$ such that the angle between 
$(n_{(\alpha_0,0)})_{x_0}$ and $(n_{(\alpha_1,0)})_{x_0}$ is equal to $\frac{2\pi}{3}$.  
Let $P_1$ be the complex affine line through 
$(n_{(\alpha_0,0)})_{x_0}$ and $(n_{(\alpha_1,0)})_{x_0}$, and 
$P'$ the (complex) $2$-dimensional complex affine subspace through 
${\bf 0},\,\,(n_{(\alpha_0,0)})_{x_0}$ and $(n_{(\alpha_1,0)})_{x_0}$.  
It is clear that $P_1\subset P'$.  
Also, it is easy to show that the root system associated with 
$L^{P_1}_{x_0}$ is of type ($A_2$) and hence the affine root system 
associated with $L^{P'}_{x_0}$ is of type ($\widetilde A_2$).  
This completes the proof of the statement (i).  

Next we shall show the statement (ii).  
Since $\triangle_M$ is of type $(B_m)$, the positive root system 
$(\triangle_M)_+$ is described as 
$$(\triangle_M)_+=\{\theta_a\,\vert 1\leq a\leq m\}
\cup\{\theta_a\pm\theta_b\,\vert 1\leq a<b\leq m\}$$
for an orthonormal base $\theta_1,\cdots,\theta_m$ of the dual space 
$\mathfrak b^{\ast}$ of $\mathfrak b$, the simple root system $\Pi$ is equal 
to $\{\theta_i-\theta_{i+1}\,\vert\,1\leq i\leq n-1\}\cup\{\theta_n\}$ and 
the highest root is equal to $\theta_1+\theta_2$, where we need to replace 
the inner product $\langle\,\,,\,\,\rangle\vert_{\mathfrak b_{\Bbb R}\times\mathfrak b_{\Bbb R}}$ 
to its suitable constant-multiple.  
Without loss of generality, we may assume that $\alpha_0$ is one of 
the elements of $\Pi$.  In the case where $\alpha_0$ is other than $\theta_n$, 
there exists $\alpha_1\in\Pi$ such that the angle between $(n_{(\alpha_0,0)})_{x_0}$ and 
$(n_{(\alpha_1,0)})_{x_0}$ is equal to $\frac{2\pi}{3}$.  
Let $P_1$ be the complex affine line through 
$(n_{(\alpha_0,0)})_{x_0}$ and $(n_{(\alpha_1,0)})_{x_0}$, and 
$P'$ the (complex) $2$-dimensional complex affine subspace through 
${\bf 0},\,\,(n_{(\alpha_0,0)})_{x_0}$ and $(n_{(\alpha_1,0)})_{x_0}$.  
Then it is shown that 
the root system associated with $L^{P_1}_{x_0}$ is of type $(A_2)$ 
and hence the affine root system associated with $L^{P'}_{x_0}$ is of type $(\widetilde A_2)$.  
In the case where $\alpha_0$ is equal to $\theta_n$, 
we can take $\alpha_1\in\Pi$ such that the angle between 
$(n_{(\alpha_0,0)})_{x_0}$ and $(n_{(\alpha_1,0)})_{x_0}$ is equal to $\frac{3\pi}{4}$.  
Let $P_1$ be the complex affine line through $(n_{(\alpha_0,0)})_{x_0}$ and 
$(n_{(\alpha_1,0)})_{x_0}$, and $P'$ the (complex) $2$-dimensional complex affine subspace through 
${\bf 0},\,\,(n_{(\alpha_0,0)})_{x_0}$ and $(n_{(\alpha_1,0)})_{x_0}$.  
Then it is shown that, in correspondence to ${\cal W}$ is of type 
$(\widetilde B_m),\,(\widetilde B_m^v)$ or $(\widetilde B_m,\widetilde B_m^v)$ 
($m\geq 2$), the root system associated with 
$L^{P_1}_{x_0}$ is of type $(C_2),\,(C_2^v)$ or $(C_2,C_2^v)$ 
and hence the affine root system associated with $L^{P'}_{x_0}$ is of type 
$(\widetilde C_2),\,(\widetilde C_2^v)$ or $(\widetilde C_2,\widetilde C_2^v)$.  

Next we shall show the statement (iii).  
Since $\triangle_M$ is of type $(C_m)$, the positive root system 
$(\triangle_M)_+$ is described as 
$$(\triangle_M)_+=\{2\theta_a\,\vert 1\leq a\leq m\}
\cup\{\theta_a\pm\theta_b\,\vert 1\leq a<b\leq m\}$$
for an orthonormal base $\theta_1,\cdots,\theta_m$ of the dual space 
$\mathfrak b^{\ast}$, the simple root system $\Pi$ is equal to 
$\{\theta_i-\theta_{i+1}\,\vert\,1\leq i\leq n-1\}\cup\{2\theta_n\}$ and 
the highest root is equal to $2\theta_1$, where we need to replace 
the inner product 
$\langle\,\,,\,\,\rangle\vert_{\mathfrak b_{\Bbb R}\times\mathfrak b_{\Bbb R}}$ to its suitable 
constant-multiple.  
Without loss of generality, we may assume that $\alpha_0$ is one of the 
elements of $\Pi$.  
In the case where $\alpha_0$ is other than $2\theta_n$, 
there exists $\alpha_1\in(\triangle_M)_+$ such that the angle 
between $(n_{(\alpha_0,0)})_{x_0}$ and $(n_{(\alpha_1,0)})_{x_0}$ is equal to $\frac{2\pi}{3}$.  
Let $P_1$ be the complex affine line through $(n_{(\alpha_0,0)})_{x_0}$ and 
$(n_{(\alpha_1,0)})_{x_0}$, and $P'$ the (complex) $2$-dimensional complex affine subspace through 
${\bf 0},\,\,(n_{(\alpha_0,0)})_{x_0}$ and $(n_{(\alpha_1,0)})_{x_0}$.  
Then it is shown that the root system associated with $L^{P_1}_{x_0}$ is of type $(A_2)$ 
and hence the affine root system associated with $L^{P'}_{x_0}$ is of type $(\widetilde A_2)$.  
In the case where $\alpha_0$ is equal to $2\theta_n$, 
we can take $\alpha_1\in(\triangle_M)_+$ such that the angle 
between $(n_{(\alpha_0,0)})_{x_0}$ and $(n_{(\alpha_1,0)})_{x_0}$ is equal to $\frac{3\pi}{4}$.  
Let $P_1(\subset \mathfrak b^{\Bbb C})$ be the complex affine line through 
$(n_{(\alpha_0,0)})_{x_0}$ and $(n_{(\alpha_1,0)})_{x_0}$, and 
$P'$ the (complex) $2$-dimensional complex affine subspace through 
${\bf 0},\,\,(n_{(\alpha_0,0)})_{x_0}$ and $(n_{(\alpha_1,0)})_{x_0}$.  
Then it is shown that, in correspondence to ${\cal W}$ is of type 
$(\widetilde C_m),\,(\widetilde C_m^v),\,(\widetilde C'_m),\,
(\widetilde C_m^v,\widetilde C'_m),\,(\widetilde C'_m,\widetilde C_m),\,
(\widetilde C_m^v,\widetilde C_m)$ or $(\widetilde C_m,\widetilde C_m^v)$ 
($m\geq 2$), 
the root system associated with $L^{P_1}_{x_0}$ is of type 
$(C_2),\,(C_2^v),\,(C'_2),\,(C_2^v,C'_2),\,(C'_2,C_2),\,(C_2^v,C_2)$ or 
$(C_2,C_2^v)$ and hence the affine root system associated 
with $L^{P'}_{x_0}$ is of type 
$(\widetilde C_2),\,(\widetilde C_2^v),\,(\widetilde C'_2),\,
(\widetilde C_2^v,\widetilde C'_2),$\newline
$(\widetilde C'_2,\widetilde C_2),\,(\widetilde C_2^v,\widetilde C_2)$ or 
$(\widetilde C_2,\widetilde C_2^v)$.  
\begin{flushright}q.e.d.\end{flushright}

\vspace{0.5truecm}

Also, we can show the following fact.  

\vspace{0.5truecm}

\noindent
{\bf Lemma 5.24.} {\sl If the affine root system ${\cal R}$ is of type 
$(\widetilde G_2)$ and if $\langle n_{i_1},n_{i_2}\rangle=0$, then 
$\Gamma_{w_{i_1}}w_{i_2}=0$ for any $w_{i_1}\in (E_{i_1})_{x_0}$ and $w_{i_2}\in (E_{i_2})_{x_0}$.}

\vspace{0.5truecm}

\noindent
{\it Proof.} Let $i_k=(\alpha_k,j_k)$ ($k=1,2$).  Let $P$ be the complex affine 
line through $(n_{i_1})_{x_0}$ and $(n_{i_2})_{x_0}$.  
Since $\langle n_{i_1},n_{i_2}\rangle=0$, we have 
$\langle(n_{i_1})_{x_0},(n_{i_2})_{x_0}\rangle=0$.  
If there does not exist further $i_3\in I$ with $(n_{i_3})_{x_0}\in P$, 
then the root system associated with the slice $L^P_{x_0}$ is of type $(A_1\times A_1)$.  
Hence we have $\Gamma_{(E_{i_1})_{x_0}}(E_{i_2})_{x_0}=0$.  
Otherwise, it is shown that $\{i\in I\,\vert\,(n_i)_{x_0}\in P\}$ consists 
of exactly six elements because $\triangle_M$ is of type $(G_2)$, where 
we note that 
$\{i\in I\,\vert\,(n_i)_{x_0}\in P\}
=\{i\in I\,\vert\,(n_i)_{x_0}\in P\cap\mathfrak b_{\Bbb R}\}$ and that 
each $P\cap\mathfrak b_{\Bbb R}$ is a real affine line in $\mathfrak b_{\Bbb R}$.  
The root system $\triangle_P$ associated with the slice 
$L^P_{x_0}(\subset(W_P)_{x_0})$ is of type $(G_2)$.  The slice $L^P_{x_0}$ 
is regarded as a principal orbit of the isotropy action of an 
anti-Kaehler symmetric space $L/H$ whose root system is of type $(G_2)$.  
Let $\mathfrak l=\mathfrak h+\mathfrak p$ be the canonical decomposition of 
the Lie algebra $\mathfrak l$ of $L$ associated with the symmetric pair $(L,H)$.  
The space $\mathfrak p$ is identified with $(W_P)_{x_0}$ and 
the normal space of $L^P_{x_0}(\subset(W_P)_{x_0})$ at $x_0$ is 
identified with a maximal abelian subspace $\mathfrak b'$ of $\mathfrak p$.  
Denote by $\mathfrak p_{\overline{\alpha}}(\subset\mathfrak p)$ and 
$\mathfrak h_{\overline{\alpha}}(\subset\mathfrak h)$ be the root spaces for 
$\overline{\alpha}\in\triangle_P$.  
The restriction $\overline{\alpha}_k:=\alpha_k\vert_{\mathfrak b'}$ of $\alpha_k$ to 
$\mathfrak b'$ ($k=1,2$) are elements of $\triangle_P$, where $\mathfrak b'$ is regaraded as a 
linear subspace of $\mathfrak b$ under the identification of $\mathfrak b'$ and 
the normal space $T^{\perp}_{x_0}L^P_{x_0}$ of $L^P_{x_0}$ in $(W_P)_{x_0}$.  
For any $w_k\in(E_{i_k})_{x_0}$ ($k=1,2$), we have 
$$\Gamma_{w_1}w_2\in[\mathfrak h_{\overline{\alpha}_1},\,\mathfrak p_{\overline{\alpha}_2}]
\subset\mathfrak p_{\overline{\alpha}_1+\overline{\alpha}_2}
+\mathfrak p_{\overline{\alpha}_1-\overline{\alpha}_2}.$$
Since $\overline{\alpha}_1$ and $\overline{\alpha}_2$ are orthogonal and $\triangle_P$ is of type 
$(G_2)$, we have $\overline{\alpha}_1\pm\overline{\alpha}_2\notin\triangle_P$.  
Hence we have $\Gamma_{w_1}w_2=0$.  
This completes the proof.  
\begin{flushright}q.e.d.\end{flushright}

\vspace{0.5truecm}

%
%

By using Lemmas 5.6, 5.7, 5.10, 5.11, 5.14, 5.23, 5.24 and Lemma 8.3 of [GH], we can show the following fact.  

\vspace{0.5truecm}

\noindent
{\bf Theorem 5.25.} {\sl If ${\cal R}$ is of type $(\widetilde A_m)$ 
($m\geq2$), $(\widetilde D_m)$ ($m\geq4$), $(\widetilde E_6),\,
(\widetilde E_7),\,(\widetilde E_8),\,(\widetilde F_4)$ or 
$(\widetilde G_2)$, then $\Gamma_{(E_{(\alpha,j_1)})_{x_0}}(E_{(\alpha,j_2)})_{x_0}
\subset(E_0)_{x_0}$ holds for any $\alpha\in(\triangle_M)_+$ and $j_1,j_2\in{\Bbb Z}$.}

\vspace{0.5truecm}

\noindent
{\it Proof.} According to Lemma 5.23, we may assume that ${\cal R}$ is of type 
$(\widetilde A_2)$ or $(\widetilde G_2)$.  
Furthermore, according to Lemma 5.6, we may assume that $j_1\not=j_2$.  
Set $i_k:=(\alpha,j_k)$ ($k=1,2$).  
Suppose that $(\Gamma_{w_1}w_2)^{i_3}\not=0$ for some $w_k\in(E_{i_k})_{x_0}$ ($k=1,2$) 
and some $i_3\in I$.  
Take $w_k\in (E_{i_k})_{x_0}$ ($k=1,2$) with $(\Gamma_{w_1}w_2)^{i_3}\not=0$.  
Let $P$ be the complex affine line through ${\bf 0}$ and $(n_{i_1})_{x_0}$.  
Since $L^P_{x_0}$ is totally geodesic in $M$, we have $(E_{i_3})_{x_0}\subset T_{x_0}M$ and hence 
$(n_{i_3})_{x_0}\in P$.  
Hence $i_3$ is expressed as $i_3=(\alpha,j_3)$ for some $j_3\in{\Bbb Z}$.  
According to Lemma 5.7, we have $j_3\not=j_1,j_2$.  
According to Lemma 5.11, there exists $i_4,i_5\in I$ such that 
$(n_{i_4})_{x_0}$ and $(n_{i_5})_{x_0}$ are ${\Bbb C}$-linearly 
independent and that $\langle(\Gamma_{w_1}w_2)^{i_3},\Gamma_{w_5}w_4\rangle
\not=0$ for some $w_4\in(E_{i_4})_{x_0}$ and some 
$w_5\in(E_{i_5})_{x_0}$.  
Since $\langle(\Gamma_{w_1}w_2)^{i_3},\Gamma_{w_5}w_4\rangle\not=0$, we have 
$(\Gamma_{w_5}w_4)^{i_3}\not=0$.  Hence it follows from Lemma 5.10 that 
$(n_{i_3})_{x_0},(n_{i_4})_{x_0}$ and $(n_{i_5})_{x_0}$ are contained 
in a complex affine line $P_1$.  Since $P\cap P_1=\{(n_{i_3})_{x_0}\}$, 
it follows from Lemma 5.10 that $\langle\Gamma_{w_1}w_2,\Gamma_{w_5}w_4\rangle
=\langle(\Gamma_{w_1}w_2)^{i_3},(\Gamma_{w_5}w_4)^{i_3}\rangle\not=0$.  Also, it is clear that 
arbitrarily chosen three of $(n_{i_1})_{x_0},(n_{i_2})_{x_0},(n_{i_4})_{x_0}$ and 
$(n_{i_5})_{x_0}$ are not contained in any complex affine line.  
Hence, it follows from Lemma 5.15 that 
$$\begin{array}{l}
\langle\Gamma_{w_1}w_2,\Gamma_{w_5}w_4\rangle
=\langle\Gamma_{w_5}w_2,\Gamma_{w_1}w_4\rangle
+c\langle\Gamma_{w_1}w_5,\Gamma_{w_2}w_4\rangle,
\end{array}$$
where $c$ is as in Lemma 5.15.  
Hence we have 
$$({\rm I})\quad\langle\Gamma_{w_5}w_2,\Gamma_{w_1}w_4\rangle\not=0\qquad{\rm or}\qquad
({\rm II})\quad\langle\Gamma_{w_1}w_5,\Gamma_{w_2}w_4\rangle\not=0.$$
We consider the case of (I).  
According to Lemma 5.10, this fact implies that 
the complex affine line through $(n_{i_2})_{x_0}$ and $(n_{i_5})_{x_0}$ intersects with 
the complex affine line through and $(n_{i_1})_{x_0}$ and $(n_{i_4})_{x_0}$ and 
the only intersection point is equal to $(n_{i_6})_{x_0}$ for some $i_6\in I$.  
Then, since $(n_{i_1})_{x_0},(n_{i_2})_{x_0}$ and $(n_{i_3})_{x_0}$ are $C$-linearly dependent 
pairwisely, the complex focal hyperplanes ${\it l}_{i_1},{\it l}_{i_2}$ and ${\it l}_{i_3}$ are 
mutually parallel.  
Note that they are complex lines because we assume that ${\cal R}$ is 
of type $(\widetilde A_2)$ or $(\widetilde G_2)$.  
Hence the (real) lines ${\it l}_{i_1}^{\Bbb R},{\it l}_{i_2}^{\Bbb R}$ and 
${\it l}_{i_3}^{\Bbb R}$ (in $\mathfrak b_{\Bbb R}$) are mutually parallel.  
Also, since $(n_{i_3})_{x_0},(n_{i_4})_{x_0}$ and $(n_{i_5})_{x_0}$ 
are contained in a complex line which does not pass $0$, 
we have ${\it l}_{i_3},{\it l}_{i_4}$ and ${\it l}_{i_5}$ have a common point.  
Hence the lines ${\it l}_{i_3}^{\Bbb R},{\it l}_{i_4}^{\Bbb R}$ and 
${\it l}_{i_5}^{\Bbb R}$ have a common point.  Denote by $p_{345}$ this common point.  
Similarly, since $(n_{i_2})_{x_0},(n_{i_5})_{x_0}$ and 
$(n_{i_6})_{x_0}$ are contained in a complex line which does not pass $0$, 
we have ${\it l}_{i_2},{\it l}_{i_5}$ and ${\it l}_{i_6}$ have a common 
point.  Hence the lines ${\it l}_{i_2}^{\Bbb R},{\it l}_{i_5}^{\Bbb R}$ and ${\it l}_{i_6}^{\Bbb R}$ 
have a common point.  Denote by $p_{256}$ this common point.  
Also, since $(n_{i_1})_{x_0},(n_{i_4})_{x_0}$ and $(n_{i_6})_{x_0}$ 
are contained in a complex line which does not pass $0$, 
${\it l}_{i_1},{\it l}_{i_4}$ and ${\it l}_{i_6}$ have a common point.  
Hence the lines ${\it l}_{i_1}^{\Bbb R},{\it l}_{i_4}^{\Bbb R}$ and 
${\it l}_{i_6}^{\Bbb R}$ have a common point.  Denote by $p_{146}$ this common point.  
These three intersection points $p_{345},p_{256}$ and $p_{146}$ 
lie in no line in $\mathfrak b_-$ because of $i_4\not=i_5$.  
On the other hand, in the case where ${\cal R}$ is of type $(\widetilde A_2$), it is clear that 
the angle between arbitrarily chosen two of ${\it l}_{i_k}^{\Bbb R}$ ($k=1,\cdots,6$) is equal to 
an integer-multiple of $\frac{\pi}{6}$ other than $\frac{\pi}{2}$.  
Also, in the case where ${\cal R}$ is of type $(\widetilde G_2$), 
it follows from Lemmas 5.10 and 5.24 that 
the angle between arbitrarily chosen two of ${\it l}_{i_k}^{\Bbb R}$ ($k=1,\cdots,6$) is equal to 
an integer-multiple of $\frac{\pi}{6}$ other than $\frac{\pi}{2}$.  
Hence, it follows from (i) of Lemma 5.25 that 
$p_{345},p_{256}$ and $p_{146}$ lie in a line in $\mathfrak b_{\Bbb R}$.  
Thus a contradiction arises.  
Similarly, in case of (II), we can drive a contradiction.  
Therefore we obtain $(\Gamma_{w_1}w_2)^{i_3}=0$.  
It follows from the arbitrariness of $i_3$ that $\Gamma_{w_1}w_2\in(E_0)_{x_0}$.  
This completes the proof.  
\begin{flushright}q.e.d.\end{flushright}

\vspace{0.5truecm}

From Lemmas 5.17 and 5.21 and Theorem 5.25, we have the following fact .  

\vspace{0.5truecm}

\noindent
{\bf Proposition 5.26.} {\sl If ${\cal R}$ is one of the following types:
$$
(\widetilde A_m)\,\,(m\geq2),\,\,\,\,(\widetilde D_m)\,\,(m\geq4),\,\,\,\,
(\widetilde E_6),\,\,\,(\widetilde E_7),\,\,\,(\widetilde E_8),\,\,\,
(\widetilde F_4),\,\,\,(\widetilde F_4^v),\,\,\,(\widetilde G_2),\,\,\,
(\widetilde G_2^v),
$$
then $\Gamma_w$ can be extended continuously to $T_{x_0}M$ for any 
$\displaystyle{w\in\mathop{\cup}_{i\in I}E_i}$.}

\vspace{0.5truecm}

\noindent
{\it Proof.} Let $\alpha\in(\triangle_M)_+$ and 
$j_1,j_2\in{\Bbb Z}$.  
Set $i_k:=(\alpha,j_k)$ ($k=1,2$).  From Lemma 5.17 and Theorem 5.25, we have 
$$\vert\vert\Gamma_{w_1}w_2\vert\vert^2
=\frac12{\rm Re}\left(\frac{n_{i_1}-0}{n_{i_2}-0}\right)
\langle n_{i_1},n_{i_2}\rangle\,\langle w_1,w_1\rangle\,
\vert\vert w_2\vert\vert^2$$
for any $w_k\in(E_{i_k})_{x_0}$ ($k=1,2$).  
Clearly we have 
$$\sup_{j\in{\bf Z}}\left\vert
{\rm Re}\left(\frac{n_{i_1}-0}{n_{(\alpha,j)}-0}\right)
\langle n_{i_1},n_{(\alpha,j)}\rangle\right\vert\,<\,\infty.$$
Denote by $C$ this supremum.  Then we have 
$$\vert\vert\Gamma_{w_1}w_2\vert\vert\leq
\sqrt{\frac C2}\vert\vert\vert w_1\vert\vert\,\vert\vert w_2\vert\vert.$$
Hence, it follows from the arbitrarinesses of $w_2$ and $j_2$ that 
$$\vert\vert\Gamma_{w_1}w\vert\vert\leq
\sqrt{\frac C2}\vert\vert\vert w_1\vert\vert\,\vert\vert w\vert\vert$$
for any $\displaystyle{w\in\mathop{\cup}_{j\in{\bf Z}}
(E_{(\alpha,j)})_{x_0}}$.  
On the other hand, since $\Gamma_{(E_{i_1})_{x_0}}(E_{(\alpha,j)})_{x_0}
\subset(E_0)_{x_0}$ ($j\in{\Bbb Z}$) by Theorem 5.25, it follows from 
Lemma 5.16 that 
$$\langle\Gamma_{(E_{i_1})_{x_0}}(E_{(\alpha,j)})_{x_0},\,
\Gamma_{(E_{i_1})_{x_0}}(E_{(\alpha,j')})_{x_0}\rangle=0$$
for any $j'\in{\Bbb Z}$ satisfying $j'\not=j_1,j,2j_1-j$.  Therefore, by using 
Lemma 5.21, we can show that 
$$\vert\vert\Gamma_{w_1}w\vert\vert\leq 
\sqrt{\frac{3C}{2}}\vert\vert\vert w_1\vert\vert\,\vert\vert w\vert\vert$$
for any $\displaystyle{w\in\mathop{\oplus}_{j\in{\bf Z}}
(E_{(\alpha,j)})_{x_0}}$.  Thus the restriction of $\Gamma_{w_1}$ to 
$\displaystyle{\mathop{\oplus}_{j\in{\bf Z}}(E_{(\alpha,j)})_{x_0}}$ is 
bounded and hence it can be extended continuosuly to 
$\displaystyle{\overline{\mathop{\oplus}_{j\in{\bf Z}}(E_{(\alpha,j)})_{x_0}}}$.  
Therefore, according to Lemma 5.20, $\Gamma_{w_1}$ can be extended continuously to $T_{x_0}M$.  
\begin{flushright}q.e.d.\end{flushright}

\vspace{0.5truecm}

From Lemmas 5.10, 5.11, 5.15 5.21, 5.23, Theorem 5.25 and Lemma 8.3 of [GH], 
we have the following fact.  

\vspace{0.5truecm}

\noindent
{\bf Lemma 5.27.} {\sl For any $\alpha\in(\triangle_M)_+$ and 
any $j_1,j_2\in{\Bbb Z}$, we have 
$$\Gamma_{(E_{(\alpha,j_1)})_{x_0}}(E_{(\alpha,j_2)})_{x_0}\subset
(E_0)_{x_0}\oplus(E_{(\alpha,2j_1-j_2)})_{x_0}
\oplus(E_{(\alpha,2j_2-j_1)})_{x_0}
\oplus(E_{(\alpha,\frac{j_1+j_2}{2})})_{x_0},$$
where the last term is omitted in the case where $j_1+j_2$ is odd.}

\vspace{0.5truecm}

\noindent
{\it Proof.} For simplicity set $i_k:=(\alpha,j_k)$ ($k=1,2$).  
According to Lemma 5.23 and Theorem 5.25, we suffice to show in the case where (${\cal R}$) is of 
type $(\widetilde C_2), (\widetilde C_2^v), (\widetilde C'_2), (\widetilde C_2^v,\widetilde C'_2)$, 
$(\widetilde C'_2,\widetilde C_2),(\widetilde C^v_2,\widetilde C_2)$ or 
$(\widetilde C_2,\widetilde C_2^v)$.  
Let $P$ be the complex affne line through ${\bf 0}$ and $(n_{(\alpha,0)})_{x_0}$.  
Since $L^P_{x_0}$ is totally geodesic in $M$, we have 
$$\Gamma_{(E_{(\alpha,j_1)})_{x_0}}(E_{(\alpha,j_2)})_{x_0}\subset
\overline{(E_0)_{x_0}\oplus\left(\mathop{\oplus}_{j\in{\bf Z}}(E_{(\alpha,j)})_{x_0}
\right)}.$$
Assume that $(\Gamma_{w_1}w_2)^{(\alpha,j_3)}\not=0$ for some 
$w_k\in(E_{i_k})_{x_0}$ ($k=1,2$) and some $j_3\in{\Bbb Z}$.  
Set $i_3:=(\alpha,j_3)$.  Then it follows from Lemma 5.7 that $j_3\not=j_1,j_2$.  
According to Lemma 5.11, there exist $i_k=(\alpha_k,j_k)$ $(k=4,5$) such that 
$\langle(\Gamma_{w_1}w_2)^{i_3},\Gamma_{w_4}w_5\rangle\not=0$ for some 
$w_k\in(E_{i_k})_{x_0}$ ($k=4,5$).  
As in the proof of Theorem 5.25, we can show 
$$({\rm I})\quad\langle\Gamma_{w_5}w_2,\Gamma_{w_1}w_4\rangle\not=0\qquad{\rm or}\qquad
({\rm II})\quad\langle\Gamma_{w_1}w_5,\Gamma_{w_2}w_4\rangle\not=0$$
in terms of Lemmas 5.10 and 5.15.  
We consider the case of (I).  According to Lemma 5.10, this fact implies that 
the complex affine line through $(n_{i_2})_{x_0}$ and $(n_{i_5})_{x_0}$ intersects with 
the complex affine line through $(n_{i_1})_{x_0}$ and $(n_{i_4})_{x_0}$ and 
the only intersection point is equal to $(n_{i_6})_{x_0}$ for some $i_6\in I$.  
Then, as in the proof of Theorem 5.25, we can show that 
${\it l}_{i_1}^{\Bbb R},{\it l}_{i_2}^{\Bbb R}$ and ${\it l}_{i_3}^{\Bbb R}$ are mutually parallel, 
that ${\it l}_{i_3}^{\Bbb R},{\it l}_{i_4}^{\Bbb R}$ and ${\it l}_{i_5}^{\Bbb R}$ have 
the common point (which we denote by $p_{345}$), that 
${\it l}_{i_2}^{\Bbb R},{\it l}_{i_5}^{\Bbb R}$ and ${\it l}_{i_6}^{\Bbb R}$ have 
the common point (which we denote by $p_{256}$) and that 
${\it l}_{i_1}^{\Bbb R},{\it l}_{i_4}^{\Bbb R}$ and ${\it l}_{i_6}^{\Bbb R}$ have 
the common point (which we denote by $p_{146}$).   
These three intersection points $p_{345},p_{256}$ and $p_{146}$ are 
lie in no line in $\mathfrak b_{\Bbb R}$ because of $i_4\not=i_5$.  
Hence, it follows from (ii) of Lemma 5.25 that one of 
${\it l}_{i_1}^{\Bbb R},{\it l}_{i_2}^{\Bbb R},{\it l}_{i_3}^{\Bbb R}$ lies in the half way 
distant between the other two, that is, one of $j_1,j_2,j_3$ is equal to the half of the sum 
of the other two (i.e., $j_3=\frac{j_1+j_2}{2},2j_1-j_2$ or $2j_2-j_1$).  
Thus we obtain the desired relation.  
Similarly, in case of (II), we can derive the desired relation.  
\begin{flushright}q.e.d.\end{flushright}

\vspace{0.5truecm}

By using Lemmas 5.16, 5.21 and 5.27, we can show the following fact in the method of the proof of 
Corollary 8.7 of [GH].  

\vspace{0.5truecm}

\noindent
{\bf Lemma 5.28.} {\sl Let $\alpha\in(\triangle_M)_+$ and 
$j_k\in{\Bbb Z}$ ($k=1,2,3$) with $j_1\not=j_2$.  
Then we have 
$\langle\Gamma_{(E_{(\alpha,j_1)})_{x_0}}(E_{(\alpha,j_2)})_{x_0},\,
\Gamma_{(E_{(\alpha,j_1)})_{x_0}}(E_{(\alpha,j_3)})_{x_0}\rangle=0$ 
if $j_3$ is not one of 
$$4j_2-3j_1,\,\,2j_2-j_1,\,\,j_2,\,\,\frac{j_1+j_2}{2},\,\,
\frac{3j_1+j_2}{4},\,\,\frac{3j_1-j_2}{2},\,\,2j_1-j_2,\,\,3j_1-2j_2.$$
} 

\vspace{0.5truecm}

Let $P$ be a complex affine line in $\mathfrak b$ 
containig exactly four $J$-curvature normals 
$(n_{(\alpha_k,j_k)})_{x_0}$ ($k=1,\cdots,4$) at $x_0$ and 
$\mathfrak b'$ the (complex) $2$-dimensional complex linear subspace of $\mathfrak b$ spanned by 
$(n_{(\alpha_k,j_k)})_{x_0}$ ($k=1,\cdots,4$).  
Set $i_k:=(\alpha_k,j_k)$ ($k=1,\cdots,4$).  
Then the root system (which we denote by $\triangle_P$) of the slice 
$L^P_{x_0}$ is of type $(B_2)$ or $(BC_2)$.  
Hence $\triangle_P$ is given by 
$$\triangle_P=\left\{
\begin{array}{ll}
\{\pm\alpha_k\vert_{\mathfrak b'\cap\mathfrak b_{\Bbb R}}\,\vert\,k=1,\cdots,4\} & 
({\rm when}\,\,\triangle_P:(B_2){\rm -type})\\
\{\pm\alpha_k\vert_{\mathfrak b'\cap\mathfrak b_{\Bbb R}}\,\vert\,k=1,\cdots,4\}
\cup\{\pm 2\alpha_k\vert_{\mathfrak b'\cap\mathfrak b_{\Bbb R}}\,\vert\,k=1,2\} & 
({\rm when}\,\,\triangle_P:(BC_2){\rm -type}),
\end{array}\right.$$
where we need to permute $i_1,\cdots,i_4$ suitably if necessary.  
If $\triangle_P$ is of type $(B_2)$, then 
$E_{i_k}$ ($k=1,\cdots,4$) are irreducible with respect to $(\Phi_{i_k})_{x_0}$, respectively, 
where $\Phi_{i_k}$ is the normal holonomy group of the focal submanifold $f_{i_k}(M)$ corresponding 
to $E_{i_k}$ at $x_0$ and $(\Phi_{i_k})_{x_0}$ is the isotropy group of $\Phi_{i_k}$ at $x_0$.  
Also, if $\triangle_P$ is of type $(BC_2)$, then $E_{i_k}$ ($k=1,2$) are reducible with respect to 
$(\Phi_{i_k})_{x_0}$, respectively, and $E_{i_k}$ ($k=3,4$) are irreducible with respect to 
$(\Phi_{i_k})_{x_0}$, respectively.  
We can show the following lemma in the method of the proof of Lemma 8.8 of [GH].  

\vspace{0.5truecm}

\noindent
{\bf Lemma 5.29.} {\sl Let $P$ be as above and 
$(E_{i_k})_{x_0}=(E'_{i_k})_{x_0}\oplus(E''_{i_k})_{x_0}$ the irreducible decomposition of 
the action $(\Phi_{i_k})_{x_0}\curvearrowright(E_{i_k})_{x_0}$, where 
${\rm dim}_{\Bbb C}(E''_{i_k})_{x_0}=0,1$ or $3$.  

{\rm (i)} If $\triangle_P$ is of type $(B_2)$, 
then we have $\Gamma_{(E_{i_3})_{x_0}}(E_{i_4})_{x_0}=0$.  

{\rm (ii)} If $\triangle_P$ is of type $(BC_2)$, 
then the $(E_{i_k})'_{x_0}$-component of 
$\Gamma_{(E_{i_3})_{x_0}}(E_{i_4})_{x_0}$ vanishes, where $k=1,2$.  

{\rm (iii)} If $\triangle_P$ is of type $(BC_2)$, 
then we have $\Gamma_{(E_{i_1})''_{x_0}}(E_{i_2})_{x_0}
=\Gamma_{(E_{i_1})_{x_0}}(E_{i_2})''_{x_0}=0$.}

\vspace{0.5truecm}

By using Lemmas 5.10, 5.23, 5.27, 5.29, Theorem 5.25 and Lemma 8.3 of [GH], we can show the 
following fact corresponding to Theorem 8.12 and Proposition 8.13 of [GH].  

\vspace{0.5truecm}

\noindent
{\bf Lemma 5.30.} {\sl 
{\rm (i)} If $E_{(\alpha,j_1)}$ is irreducible and if 
$j_1-j_2$ is divisible by $4$ or the affine root system ${\cal R}$ associated 
with $M$ is not of type $(\widetilde C_n)$ ($n\geq2$), then 
we have $\Gamma_{(E_{(\alpha,j_1)})_{x_0}}(E_{(\alpha,j_2)})_{x_0}\subset
(E_0)_{x_0}$.  

{\rm (ii)} If $E_{(\alpha,j_1)}$ is irreducible and if 
$j_1-j_2$ is even, then we have $\Gamma_{(E_{(\alpha,j_1)})_{x_0}}
(E_{(\alpha,j_2)})_{x_0}\subset(E_0)_{x_0}\oplus
(E_{(\alpha,\frac{j_1+j_2}{2})})_{x_0}$.  

{\rm (iii)} If $E_{(\alpha,j_1)}$ is reducible and if 
$j_1-j_2$ is even ($j_1\not=j_2$), then we have 
$$(\Gamma_{(E''_{(\alpha,j_1)})_{x_0}}
(E_{(\alpha,j_2)})_{x_0})^{(\alpha,\frac{j_1+j_2}{2})}=0.$$
Furthermore, if $j_1-j_2$ is divisible by $4$, then 
$E_{(\alpha,\frac{j_1+j_2}{2})}$ is reducible and 
the $(E'_{(\alpha,\frac{j_1+j_2}{2})})_{x_0}$-component of each element of 
$\Gamma_{(E_{(\alpha,j_1)})_{x_0}}(E_{(\alpha,j_2)})_{x_0}$ vanishes.}

\vspace{0.5truecm}

For $\alpha\in(\triangle_M)_+$, we set 
$$C_{\alpha}:=\sup_{j,j'\in{\Bbb Z}}\left\vert{\rm Re}\left(
\frac{1+j'b_{\alpha}{\bf i}}
{1+jb_{\alpha}{\bf i}}\right)^{-1}\times
{\rm Re}\left(\frac{1}{(1+jb_{\alpha}{\bf i})(1+j'b_{\alpha}{\bf i})}\right)
\right\vert^{\frac12}.$$
Clearly we have $C_{\alpha}\,\,<\,\,\infty$.  
By using Lemmas 5.22 and 5.27, we can show the following fact.  

\vspace{0.5truecm}

\noindent
{\bf Lemma 5.31.} {\sl Let $i_k=(\alpha,j_k)$ ($k=1,2$) and 
$w_k\in(E_{i_k})_{x_0}$ ($k=1,2$).  If $j_1-j_2$ is not divisible by $2^m$, then we have 
$$\vert\vert\Gamma_{w_1}w_2\vert\vert\leq 2^{m-1}C_{\alpha}\,
\vert\vert(n_{(\alpha,0)})_{x_0}\vert\vert\,\,\vert\vert w_1\vert\vert\,\,\vert\vert w_2\vert\vert,$$
where $m$ is a positive integer.}

\vspace{0.5truecm}

\noindent
{\it Proof.} From Lemmas 5.22 and 5.27, we have 
$$\begin{array}{l}
\hspace{0.5truecm}\displaystyle{2\vert\vert(\Gamma_{w_1}w_2)^{(\alpha,2j_1-j_2)}
\vert\vert^2+\frac12\vert\vert(\Gamma_{w_1}w_2)^{(\alpha,2j_2-j_1)}
\vert\vert^2}\\
\hspace{0.5truecm}\displaystyle{
-\vert\vert(\Gamma_{w_1}w_2)^{(\alpha,\frac{j_1+j_2}{2})}\vert\vert^2
+\vert\vert(\Gamma_{w_1}w_2)^0\vert\vert^2}\\
\displaystyle{=\frac12{\rm Re}\left(\frac{1+j_2b_{\alpha}{\bf i}}
{1+j_1b_{\alpha}{\bf i}}\right)^{-1}
{\rm Re}\left(\frac{1}{(1+j_1b_{\alpha}{\bf i})(1+j_2b_{\alpha}{\bf i})}\right)}\\
\hspace{0.5truecm}\displaystyle{\times\langle(n_{(\alpha,0)})_{x_0},(n_{(\alpha,0)})_{x_0}\rangle
\langle w_1,w_1\rangle\vert\vert w_2\vert\vert^2.}
\end{array}$$
By multiplying $2$ to both sides and adding 
$3\vert\vert(\Gamma_{w_1}w_2)^{(\alpha,\frac{j_1+j_2}{2})}\vert\vert^2$ to both 
sides, we obtain 
$$\begin{array}{l}
\displaystyle{\vert\vert\Gamma_{w_1}w_2\vert\vert^2
\leq\left\vert{\rm Re}\left(\frac{1+j_2b_{\alpha}{\bf i}}{1+j_1b_{\alpha}{\bf i}}\right)^{-1}
{\rm Re}\left(\frac{1}{(1+j_1b_{\alpha}{\bf i})(1+j_2b_{\alpha}{\bf i})}\right)\right\vert}\\
\hspace{2.3truecm}\displaystyle{\times
\vert\vert(n_{(\alpha,0)})_{x_0}\vert\vert^2\vert\vert w_1\vert\vert^2\vert\vert w_2\vert\vert^2}\\
\hspace{2.3truecm}\displaystyle{
+3\vert\vert(\Gamma_{w_1}w_2)^{(\alpha,\frac{j_1+j_2}{2})}\vert\vert^2}\\
\hspace{1.8truecm}\displaystyle{\leq C^2_{\alpha}
\vert\vert(n_{(\alpha,0)})_{x_0}\vert\vert^2\vert\vert w_1\vert\vert^2\vert\vert w_2\vert\vert^2
+3\vert\vert(\Gamma_{w_1}w_2)^{(\alpha,\frac{j_1+j_2}{2})}\vert\vert^2.}
\end{array}\leqno{(5.9)}$$
We use the induction on $m$.  In case of $m=1$, the statement of this lemma is derived from 
$(5.9)$ directly.  
Now we assume that the statement of this lemma holds for $m(\geq1)$ and 
that $j_1-j_2$ is not divisible by $2^{m+1}$.  
Set $w:=(\Gamma_{w_1}w_2)^{(\alpha,\frac{j_1+j_2}{2})}$.  
Since $F^{w_1}_t$'s are holomorphic isometries, 
$\Gamma_{w_1}$ preserves $(T_{x_0}M)_-$ and 
$(T_{x_0}M)_+$ invariantly, respectively.  
Hence we have $\Gamma_{w_1}(({w_2})_{\varepsilon})
=(\Gamma_{w_1}w_2)_{\varepsilon}$ ($\varepsilon=-$ or $+$).  
Also, it follows from the definitions of 
$(T_{x_0}M)_{\varepsilon}$ ($\varepsilon=-$ or $+$) 
that $((\Gamma_{w_1}w_2)_{\varepsilon})^{(\alpha,\frac{j_1+j_2}{2})}
=((\Gamma_{w_1}w_2)^{(\alpha,\frac{j_1+j_2}{2})})_{\varepsilon}$ 
($\varepsilon=-$ or $+$).  
From (i) of Lemma 5.2 and these relations, we have 
$$\langle(\Gamma_{w_1}w_2)_{\varepsilon},w_{\varepsilon}\rangle
=\langle\Gamma_{w_1}w_2,w_{\varepsilon}\rangle
=-\langle(w_2)_{\varepsilon},(\Gamma_{w_1}w)_{\varepsilon}\rangle.$$
Hence we have 
$$\langle\Gamma_{w_1}w_2,w\rangle_{\pm}
=-\langle w_2,\Gamma_{w_1}w\rangle_{\pm}.\leqno{(5.10)}$$
Since $j_1-\frac{j_1+j_2}{2}$ is not divisible by $2^m$, it follows from $(5.10)$ and the assumption 
in the induction that 
$$\begin{array}{l}
\hspace{0.5truecm}
\displaystyle{\vert\vert(\Gamma_{w_1}w_2)^{(\alpha,\frac{j_1+j_2}{2})}
\vert\vert^2=\langle\Gamma_{w_1}w_2,w\rangle_{\pm}
=-\langle w_2,\Gamma_{w_1}w\rangle_{\pm}}\\
\displaystyle{\leq\vert\vert w_2\vert\vert\,\,\vert\vert\Gamma_{w_1}w\vert\vert
\leq 2^{m-1}C_{\alpha}
\vert\vert(n_{(\alpha,0)})_{x_0}\vert\vert\,\,\vert\vert w_1\vert\vert\,\,
\vert\vert w\vert\vert\,\,\vert\vert w_2\vert\vert,}
\end{array}$$
that is, 
$$\begin{array}{l}
\displaystyle{
\vert\vert(\Gamma_{w_1}w_2)^{(\alpha,\frac{j_1+j_2}{2})}\vert\vert
\leq 2^{m-1}C_{\alpha}\vert\vert(n_{(\alpha,0)})_{x_0}\vert\vert\,\,\vert\vert w_1\vert\vert\,\,
\vert\vert w_2\vert\vert.}
\end{array}$$
From this inequality and $(5.9)$, we obtain 
$$\vert\vert\Gamma_{w_1}w_2\vert\vert\leq 2^mC_{\alpha}
\vert\vert(n_{(\alpha,0)})_{x_0}\vert\vert\,\,\vert\vert w_1\vert\vert\,\,\vert\vert w_2\vert\vert.$$
Thus the statement of this lemma holds for $m+1$.  
Therefore the statement of this lemma is true for all $m\in{\Bbb Z}$.  
\begin{flushright}q.e.d.\end{flushright}

\vspace{0.5truecm}

By using Lemmas 5.7, 5.19, 5.21, 5.22, 5.27, 5.28, 5.30 and 5.31, we shall prove 
Theorem 5.1.  

\vspace{0.5truecm}

\noindent
{\it Proof of Theorem 5.1.} Let $i=(\alpha,j)\in I$ and $w\in(E_i)_{x_0}$.  
We suffice to show that $\Gamma_w$ is bounded in order to show that $X^w$ is 
defined on the whole of $V$.  

\noindent
{\bf (Step I)} First we shall show that, in the case where $j'$ is an integer with $j'\not=j$ 
such that $j'-j$ is devided by $4$, there exists a positive constant 
$\bar C_{\alpha}$ depending on only $\alpha$ such that 
$$\vert\vert(\Gamma_ww')^{(\alpha,\frac{j+j'}{2})}\vert\vert\leq{\bar C}_{\alpha}
\vert\vert(n_{(\alpha,0)})_{x_0}\vert\vert\,\,\vert\vert w\vert\vert\,\,\vert\vert w'\vert\vert
\leqno{(5.11)}$$
holds for any $w'\in(E_{(\alpha,j')})_{x_0}$.  
If $(E_i)_{x_0}$ is irreducible with respect to $(\Phi_i)_{x_0}$ or 
"$(E_i)_{x_0}$ is reducible with respect to $(\Phi_i)_{x_0}$ and $w\in(E''_i)_{x_0}$", 
then the left-hand side of $(5.11)$ vanishes by (i) and (iii) of Lemma 5.30.  
In the sequel, we consider the case where $(E_i)_{x_0}$ is reducible and where 
$w\in(E'_i)_{x_0}$.  Set $i':=(\alpha,j'),\,i'':=(\alpha,\frac{j+j'}{2})$ and 
$w'':=(\Gamma_ww')^{i''}$.  
According to (iii) of Lemma 5.30, we have 
$w''\in(E''_{i''})_{x_0}$.  
In similar to $(5.10)$, we have 
$$\langle\Gamma_ww',w''\rangle_{\pm}
=-\langle w',\Gamma_ww''\rangle_{\pm}.\leqno{(5.12)}$$
From this relation, we have 
$$\vert\vert(\Gamma_ww')^{i''}\vert\vert^2
=\langle\Gamma_ww',w''\rangle_{\pm}
=-\langle w',(\Gamma_ww'')^{i'}\rangle_{\pm}
\leq\vert\vert w'\vert\vert\,\,\vert\vert(\Gamma_ww'')^{i'}\vert\vert.\leqno{(5.13)}$$
On the other hand, it follows from Lemma 5.27 that 
$$\Gamma_ww''=(\Gamma_ww'')^0+(\Gamma_ww'')^{i'}
+(\Gamma_ww'')^{(\alpha,(3j-j')/2)}+(\Gamma_ww'')^{(\alpha,(3j+j')/4)}.$$
Hence, by using Lemma 5.22, we can show 
$$\begin{array}{l}
\hspace{0.5truecm}\displaystyle{\frac12\vert\vert(\Gamma_ww'')^{i'}\vert\vert^2
+2\vert\vert(\Gamma_ww'')^{(\alpha,(3j-j')/2)}\vert\vert^2}\\
\hspace{0.5truecm}\displaystyle{
-\vert\vert(\Gamma_ww'')^{(\alpha,(3j+j')/4)}\vert\vert^2
+\vert\vert(\Gamma_ww'')^0\vert\vert^2}\\
\displaystyle{\leq\frac12\left\vert{\rm Re}\left(\frac{1+jb_{\alpha}{\bf i}}
{1+((j+j')/2)b_{\alpha}{\bf i}}\right)^{-1}
{\rm Re}\left(\frac{1}{(1+((j+j')/2)b_{\alpha}{\bf i})(1+jb_{\alpha}{\bf i})}\right)\right\vert}\\
\hspace{0.5truecm}\displaystyle{\times
\vert\vert(n_{(\alpha,0)})_{x_0}\vert\vert^2\vert\vert w''\vert\vert^2\vert\vert w\vert\vert^2.}
\end{array}$$
Also, it follows from (iii) of Lemma 5.30 that $(\Gamma_{w''}w)^{(\alpha,(3j+j')/4)}=0$.  
Hence we obtain 
$$\begin{array}{l}
\hspace{0.5truecm}\displaystyle{\vert\vert(\Gamma_ww'')^{i'}\vert\vert}\\
\displaystyle{\leq\left\vert
{\rm Re}\left(\frac{1+jb_{\alpha}{\bf i}}{1+((j+j')/2)b_{\alpha}{\bf i}}\right)^{-1}
{\rm Re}\left(\frac{1}{(1+((j+j')/2)b_{\alpha}{\bf i})(1+jb_{\alpha}{\bf i})}\right)
\right\vert^{\frac12}}\\
\hspace{0.5truecm}\displaystyle{\times\vert\vert(n_{(\alpha,0)})_{x_0}\vert\vert\,\,
\vert\vert w''\vert\vert\,\,\vert\vert w\vert\vert.}
\end{array}\leqno{(5.14)}$$
Easily we can show 
$$\mathop{\sup}_{j,j'\in{\bf Z}}
\left\vert{\rm Re}\left(\frac{1+jb_{\alpha}{\bf i}}{1+((j+j')/2)b_{\alpha}{\bf i}}\right)^{-1}\,
{\rm Re}\left(\frac{1}{(1+((j+j')/2)b_{\alpha}{\bf i})(1+jb_{\alpha}{\bf i})}\right)
\right\vert^{\frac12}<\infty.$$
Denote by $\bar C_{\alpha}$ this supremum.  
From $(5.13)$ and $(5.14)$, it follows that the inequality $(5.11)$ holds 
for this constant $\bar C_{\alpha}$.  

\noindent
{\bf (Step II)} From the fact shown in (Step I), Lemmas 5.19, 5.21, 5.28, 5.30 and 5.31, it follows 
that there exists a positive constant $\widehat C_{\alpha}$ depending on only $\alpha$ such that 
$$\vert\vert\Gamma_ww'\vert\vert\leq\widehat C_{\alpha}\vert\vert w\vert\vert\,\,
\vert\vert w'\vert\vert$$
for any $w'\in(E_0)_{x_0}^{\perp}$.  
Assume that $w'\in(E_0)_{x_0}$.  Then, since 
$\Gamma_ww'\in(E_0)_{x_0}^{\perp}$ by Lemma 5.7, 
we can find a sequence $\{w''_k\}$ in 
$\displaystyle{\mathop{\oplus}_{\hat i\in I}(E_{\hat i})_{x_0}}$ 
with $\lim\limits_{k\to\infty}w''_k=\Gamma_ww'$ (with respect to 
$\vert\vert\cdot\vert\vert$).  Then we have 
$$\begin{array}{l}
\displaystyle{\vert\vert\Gamma_ww'\vert\vert^2
=\lim_{k\to\infty}\langle\Gamma_ww',w''_k\rangle_{\pm}
=-\lim_{k\to\infty}\langle w',\Gamma_ww''_k\rangle_{\pm}}\\
\hspace{1.62truecm}\displaystyle{\leq\lim_{k\to\infty}\vert\vert w'\vert\vert
\,\,\vert\vert\Gamma_ww''_k\vert\vert\leq\widehat C_{\alpha}
\vert\vert w\vert\vert\,\,\vert\vert w'\vert\vert\,\,
\vert\vert\Gamma_ww'\vert\vert,}
\end{array}$$
that is, 
$$\vert\vert\Gamma_ww'\vert\vert\leq\widehat C_{\alpha}
\vert\vert w\vert\vert\,\,\vert\vert w'\vert\vert,$$
where $\widehat C_{\alpha}$ is as above.  
Thus $\Gamma_w$ is bounded.  Therefore, $X^w$ is defined on the whole of $V$.  
\begin{flushright}q.e.d.\end{flushright}

\vspace{0.5truecm}

By using Theorem 3.4, its proof (see the proof of Theorem A in [K7])) and Theorem 5.1, 
we shall prove Theorem A.  

\vspace{0.5truecm}

\noindent
{\it Proof of Theorem A.} 
Take any $i\in I$ and any $w_0\in(E_i)_{x_0}$.  
According to Theorem 5.1, $X^{w_0}$ is defined over the whole of $V$, that is, 
$F^{w_0}_1\in I_h^b(V)$.  
On the other hand, $F^{w_0}_1$ preserves $M$ invariantly.  
Hence we have $F^{w_0}_1\in H_b$.  
Since the holomorphic isometries $f_k$'s in the proof of Theorem A in [K7] are given as 
the composition of the holomorphic isometries of $F^{w_0}_1$-type, it is then shown that $f_k$'s 
are elements of $H_b$ and hence so is also the holomorphic isometry $\widehat f$ in Step IV of 
the proof of Theorem A in [K7] (see the construction of $\widehat f$ in Step IV).  
Therefore we obtain $H_b\cdot x=M$ for any $x\in M$.  
\hspace{5truecm}q.e.d.

\vspace{0.5truecm}

\noindent
{\Large\bf Appendix} 

\vspace{0.3truecm}

\noindent
In this Appendix, we give examples of elements of $I_h(V)\setminus I_h^b(V)$.  
Denote by ${\cal K}^h$ the Lie algebra of all holomorphic Killing fields on 
the whole of $V$.  
Also, denote by $\mathfrak o_{AK}(V)$ the Lie algebra of all continuous 
skew-symmetric complex linear maps from $V$ to oneself.  
Any $X\in{\cal K}^h$ is described as $X_u=Au+b\,\,(u\in V)$ for some 
$A\in\mathfrak o_{AK}(V)$ and some $b\in V$.  
Hence ${\cal K}^h$ is identified with $\mathfrak o_{AK}(V)\times V$.  
Give $\mathfrak o_{AK}(V)$ the operator norm (which we denote by 
$\vert\vert\cdot\vert\vert_{\rm op}$) associated with 
$\langle\,\,,\,\,\rangle_{\pm}$ 
and ${\cal K}^h$ the product norm of 
this norm $\vert\vert\cdot\vert\vert_{\rm op}$ of $\mathfrak o_{AK}(V)$ and 
the norm $\vert\vert\cdot\vert\vert$ of $V$.  The space ${\cal K}^h$ is 
a Banach Lie algebra with respect to this norm.  
The group $I_h^b(V)$ is a Banach Lie group consisting of all holomorphic 
isometry $f$'s of $V$ which admit 
a one-parameter transformation group $\{f_t\,\vert\,t\in{\Bbb R}\}$ of 
$V$ such that each $f_t$ is a holomorphic isometry of $V$, that 
$f_1=f$ and that 
$\displaystyle{\left.\frac{d}{dt}\right\vert_{t=0}(f_t)_{\ast}}$ is 
an element of $\mathfrak o_{AK}(V)$.  
Note that, for a general holomorphic isometry $f$ of $V$, 
$\displaystyle{\left.\frac{d}{dt}\right\vert_{t=0}(f_t)_{\ast}}$ is not 
necessarily defined on the whole of $V$ (but it can be defined on a dense 
linear subspace of $V$).  
It is clear that the Lie algebra of this Banach Lie group $I_h^b(V)$ is equal 
to ${\cal K}^h$.  

\vspace{0.25truecm}

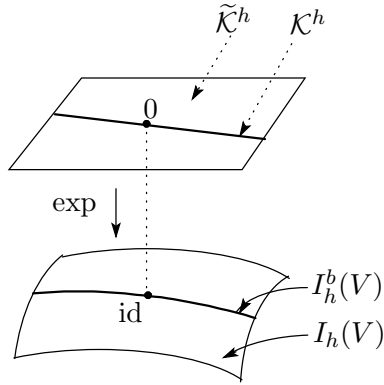
\centerline{
\unitlength 0.1in
\begin{picture}( 23.9200, 25.5000)( 12.7000,-30.6000)
%
\special{pn 20}%
\special{sh 1}%
\special{ar 2790 1160 10 10 0  6.28318530717959E+0000}%
\special{sh 1}%
\special{ar 2790 1160 10 10 0  6.28318530717959E+0000}%
%
\special{pn 8}%
\special{ar 2792 2100 836 286  4.1533611 5.8121274}%
%
\special{pn 8}%
\special{ar 3632 2580 352 648  3.2472582 4.4281071}%
%
\special{pn 8}%
\special{ar 2542 2630 836 288  4.1492541 5.8187133}%
%
\special{pn 8}%
\special{ar 2530 2554 440 774  3.3561088 4.2927807}%
%
\special{pn 20}%
\special{sh 1}%
\special{ar 2798 2058 10 10 0  6.28318530717959E+0000}%
\special{sh 1}%
\special{ar 2798 2058 10 10 0  6.28318530717959E+0000}%
%
\special{pn 13}%
\special{pa 2196 2048}%
\special{pa 2228 2046}%
\special{pa 2260 2044}%
\special{pa 2292 2042}%
\special{pa 2324 2040}%
\special{pa 2356 2038}%
\special{pa 2388 2038}%
\special{pa 2420 2038}%
\special{pa 2452 2038}%
\special{pa 2484 2038}%
\special{pa 2516 2038}%
\special{pa 2548 2040}%
\special{pa 2580 2042}%
\special{pa 2610 2044}%
\special{pa 2642 2046}%
\special{pa 2674 2048}%
\special{pa 2706 2052}%
\special{pa 2738 2056}%
\special{pa 2770 2058}%
\special{pa 2802 2062}%
\special{pa 2834 2066}%
\special{pa 2866 2068}%
\special{pa 2898 2074}%
\special{pa 2928 2078}%
\special{pa 2960 2084}%
\special{pa 2992 2090}%
\special{pa 3024 2096}%
\special{pa 3054 2102}%
\special{pa 3086 2108}%
\special{pa 3118 2114}%
\special{pa 3148 2122}%
\special{pa 3180 2130}%
\special{pa 3210 2136}%
\special{pa 3242 2144}%
\special{pa 3272 2152}%
\special{pa 3304 2160}%
\special{pa 3334 2170}%
\special{pa 3360 2178}%
\special{sp}%
\put(36.6200,-21.8700){\makebox(0,0)[lt]{$I_h(V)$}}%
\put(36.5200,-20.9800){\makebox(0,0)[lb]{$I_h^b(V)$}}%
\put(27.7000,-20.9000){\makebox(0,0)[rt]{${\rm id}$}}%
\put(28.5000,-11.3000){\makebox(0,0)[rb]{$0$}}%
\put(25.3000,-15.5000){\makebox(0,0)[rt]{$\exp$}}%
%
\special{pn 8}%
\special{ar 3626 2206 386 188  3.4640539 4.7272885}%
%
\special{pn 8}%
\special{pa 3296 2118}%
\special{pa 3278 2148}%
\special{fp}%
\special{sh 1}%
\special{pa 3278 2148}%
\special{pa 3330 2102}%
\special{pa 3306 2102}%
\special{pa 3296 2080}%
\special{pa 3278 2148}%
\special{fp}%
%
\special{pn 8}%
\special{ar 3680 2522 578 256  3.7848369 4.6324281}%
%
\special{pn 8}%
\special{pa 3258 2354}%
\special{pa 3194 2384}%
\special{fp}%
\special{sh 1}%
\special{pa 3194 2384}%
\special{pa 3264 2374}%
\special{pa 3242 2362}%
\special{pa 3246 2338}%
\special{pa 3194 2384}%
\special{fp}%
%
\special{pn 8}%
\special{pa 2450 920}%
\special{pa 2070 1400}%
\special{pa 3290 1400}%
\special{pa 3620 920}%
\special{pa 3620 920}%
\special{pa 3620 920}%
\special{pa 3620 920}%
\special{pa 2450 920}%
\special{fp}%
%
\special{pn 8}%
\special{pa 3230 700}%
\special{pa 3030 1030}%
\special{dt 0.045}%
\special{sh 1}%
\special{pa 3030 1030}%
\special{pa 3082 984}%
\special{pa 3058 984}%
\special{pa 3048 964}%
\special{pa 3030 1030}%
\special{fp}%
%
\special{pn 13}%
\special{pa 2310 1110}%
\special{pa 3410 1240}%
\special{fp}%
%
\special{pn 8}%
\special{pa 3550 740}%
\special{pa 3280 1220}%
\special{dt 0.045}%
\special{sh 1}%
\special{pa 3280 1220}%
\special{pa 3330 1172}%
\special{pa 3306 1174}%
\special{pa 3296 1152}%
\special{pa 3280 1220}%
\special{fp}%
\put(35.4000,-6.9000){\makebox(0,0)[lb]{${\cal K}^h$}}%
\put(31.5000,-6.8000){\makebox(0,0)[lb]{$\widetilde{{\cal K}}^h$}}%
%
\special{pn 8}%
\special{pa 2790 1160}%
\special{pa 2790 2060}%
\special{dt 0.045}%
\put(12.7000,-26.0000){\makebox(0,0)[lt]{$\widetilde{{\cal K}}^h:$ the space of all holomorphic Killing vector fields }}%
\put(16.3000,-28.3000){\makebox(0,0)[lt]{defined on dense linear subspaces of $V$}}%
\put(12.7000,-30.6000){\makebox(0,0)[lt]{$\exp:$ the exponential map of $I_h(V)$}}%
%
\special{pn 8}%
\special{pa 2630 1500}%
\special{pa 2630 1750}%
\special{fp}%
\special{sh 1}%
\special{pa 2630 1750}%
\special{pa 2650 1684}%
\special{pa 2630 1698}%
\special{pa 2610 1684}%
\special{pa 2630 1750}%
\special{fp}%
\end{picture}%
\hspace{1.6truecm}
}

\vspace{0.5truecm}

\centerline{{\bf Figure 2.}}

\vspace{0.3truecm}

\noindent
{\it Example.} We shall give an example of an element of 
$I_h(V)\setminus I_h^b(V)$.  
Let $V$ be a complex linear topological space consisting of all 
complex number sequences $\{z_k\}_{k=1}^{\infty}$'s satisfying 
$\sum_{k=1}^{\infty}\vert z_k\vert^2\,<\,\infty$, and 
$\langle\,\,,\,\,\rangle$ a non-degenerate inner product of $V$ 
defined by 
$$\langle\{z_k\}_{k=1}^{\infty},\,\{w_k\}_{k=1}^{\infty}\rangle
:=2{\rm Re}\left(\sum_{k=1}^{\infty}z_kw_k\right)\,\,\,\,\,\,\,\,
(\{z_k\}_{k=1}^{\infty},\,\{w_k\}_{k=1}^{\infty}\in V).$$
The pair $(V,\langle\,\,,\,\,\rangle)$ is an infinite dimensional anti-Kaehler 
space.  Define a complex linear transformation $A_t$ ($t\in{\Bbb R}$) of $V$ 
by assigning $\{w_k\}_{k=1}^{\infty}$ defined by 
$$\left(
\begin{array}{c}
w_{2k-1}\\
w_{2k}
\end{array}\right):=
\left(
\begin{array}{cc}
\cos\,2k\pi t&-\sin\,2k\pi t\\
\sin\,2k\pi t&\cos\,2k\pi t
\end{array}
\right)
\left(
\begin{array}{c}
z_{2k-1}\\
z_{2k}
\end{array}\right)
\quad\,\,(k\in{\Bbb N})$$
to each $\{z_k\}_{k=1}^{\infty}\in V$.  
It is clear that each $A_t$ is a holomorphic linear isometry of $V$.  
Define $f_t\in I_h(V)$ by $f_t(u):=A_tu+b_t\,\,(u\in V)$, where $b_t$ is a curve in $V$ with $b_0=0$.  Set 
$$B:=\left.\frac{d}{dt}\right\vert_{t=0}f_{t\ast}
=\left.\frac{d}{dt}\right\vert_{t=0}A_t.$$
It is easy to show that $B$ is a skew-symmetric complex linear map from a 
dense linear subspace $U$ of $V$ to $V$ assigning $\{w_k\}_{k=1}^{\infty}$ 
defined by 
$$\left(
\begin{array}{c}
w_{2k-1}\\
w_{2k}
\end{array}\right):=
\left(
\begin{array}{cc}
0&-2k\pi\\
2k\pi&0
\end{array}
\right)
\left(
\begin{array}{c}
z_{2k-1}\\
z_{2k}
\end{array}\right)
\quad\,\,(k\in{\Bbb N}),$$
to each $\{z_k\}_{k=1}^{\infty}\in U$, 
where $U$ is the set of all elements $\{z_k\}_{k=1}^{\infty}$'s of $V$ 
satisfying $B(\{z_k\}_{k=1}^{\infty})$\newline
$\in V$.  
Let $\{a_k\}_{k=1}^{\infty}$ be an element of $V$ defined by 
$a_k:=\frac{1}{[\frac{k+1}{2}]}\,\,\,(k\in{\Bbb N})$, where $[\cdot]$ is the 
Gauss's symbol of $\cdot$.  
Then we can show $B(\{a_k\}_{k=1}^{\infty})\notin V$, that is, 
$\{a_k\}_{k=1}^{\infty}\notin U$.  
Thus $B$ is not an element of $\mathfrak o_{AK}(V)$ and hence 
$f_t$ does not belong to $I_h^b(V)$ for positive numbers $t$'s sufficiently close to $0$, where 
we note that $f_1={\rm id}\in I_h^b(V)$.  

\vspace{0.2truecm}

\centerline{
\unitlength 0.1in
\begin{picture}( 61.0300, 18.6300)(-30.9000,-23.7000)
%
\special{pn 8}%
\special{pa 1798 858}%
\special{pa 1352 1286}%
\special{pa 2584 1286}%
\special{pa 3014 858}%
\special{pa 3014 858}%
\special{pa 1798 858}%
\special{fp}%
%
\special{pn 13}%
\special{pa 1516 1126}%
\special{pa 2872 1002}%
\special{fp}%
%
\special{pn 20}%
\special{sh 1}%
\special{ar 2176 1072 10 10 0  6.28318530717959E+0000}%
\special{sh 1}%
\special{ar 2176 1072 10 10 0  6.28318530717959E+0000}%
%
\special{pn 8}%
\special{ar 2168 1984 920 270  4.1537249 5.8094672}%
%
\special{pn 8}%
\special{ar 3094 2436 386 612  3.2474455 4.4282147}%
%
\special{pn 8}%
\special{ar 1894 2484 920 270  4.1514482 5.8150534}%
%
\special{pn 8}%
\special{ar 1880 2412 482 730  3.3531717 4.2901013}%
%
\special{pn 8}%
\special{pa 2176 1072}%
\special{pa 2176 1976}%
\special{dt 0.045}%
%
\special{pn 20}%
\special{sh 1}%
\special{ar 2176 1944 10 10 0  6.28318530717959E+0000}%
\special{sh 1}%
\special{ar 2176 1944 10 10 0  6.28318530717959E+0000}%
%
\special{pn 13}%
\special{ar 2428 2538 1616 604  4.0605721 4.9839437}%
\put(27.8500,-14.5000){\makebox(0,0)[lb]{$I_h(V)$}}%
\put(28.3500,-7.0000){\makebox(0,0)[lb]{${\cal K}^h$}}%
\put(29.8600,-16.2800){\makebox(0,0)[lb]{$I_h^b(V)$}}%
%
\special{pn 8}%
\special{pa 2980 1652}%
\special{pa 2786 1938}%
\special{dt 0.045}%
\special{sh 1}%
\special{pa 2786 1938}%
\special{pa 2840 1894}%
\special{pa 2816 1894}%
\special{pa 2806 1872}%
\special{pa 2786 1938}%
\special{fp}%
\put(24.7000,-6.7700){\makebox(0,0)[lb]{$\widetilde{\cal K}^h$}}%
\put(21.1900,-19.9700){\makebox(0,0)[lt]{${\rm id}$}}%
\put(22.0800,-10.2600){\makebox(0,0)[rb]{$0$}}%
%
\special{pn 8}%
\special{pa 1896 1412}%
\special{pa 1896 1654}%
\special{fp}%
\special{sh 1}%
\special{pa 1896 1654}%
\special{pa 1916 1586}%
\special{pa 1896 1600}%
\special{pa 1876 1586}%
\special{pa 1896 1654}%
\special{fp}%
\put(18.2400,-14.6200){\makebox(0,0)[rt]{$\exp$}}%
%
\special{pn 8}%
\special{pa 2844 726}%
\special{pa 2764 1010}%
\special{dt 0.045}%
\special{sh 1}%
\special{pa 2764 1010}%
\special{pa 2800 952}%
\special{pa 2778 960}%
\special{pa 2762 940}%
\special{pa 2764 1010}%
\special{fp}%
%
\special{pn 13}%
\special{pa 2180 1060}%
\special{pa 2408 952}%
\special{fp}%
\special{sh 1}%
\special{pa 2408 952}%
\special{pa 2340 962}%
\special{pa 2360 974}%
\special{pa 2356 998}%
\special{pa 2408 952}%
\special{fp}%
\put(18.6000,-8.9000){\makebox(0,0)[rb]{$\displaystyle{\left.\frac{d}{dt}\right\vert_{t=0}f_t}$}}%
%
\special{pn 13}%
\special{ar 2488 2110 428 262  3.8494646 4.5475129}%
%
\special{pn 8}%
\special{pa 2428 1858}%
\special{pa 2460 1858}%
\special{pa 2492 1860}%
\special{pa 2524 1862}%
\special{pa 2554 1866}%
\special{pa 2586 1870}%
\special{pa 2618 1876}%
\special{pa 2634 1880}%
\special{sp -0.045}%
\put(27.4500,-15.6900){\makebox(0,0)[rb]{{\small $t\mapsto f_t$}}}%
%
\special{pn 8}%
\special{pa 2836 1468}%
\special{pa 2686 1840}%
\special{dt 0.045}%
\special{sh 1}%
\special{pa 2686 1840}%
\special{pa 2730 1786}%
\special{pa 2706 1790}%
\special{pa 2692 1772}%
\special{pa 2686 1840}%
\special{fp}%
%
\special{pn 8}%
\special{pa 2468 1612}%
\special{pa 2318 1874}%
\special{dt 0.045}%
\special{sh 1}%
\special{pa 2318 1874}%
\special{pa 2368 1826}%
\special{pa 2344 1828}%
\special{pa 2334 1806}%
\special{pa 2318 1874}%
\special{fp}%
%
\special{pn 8}%
\special{pa 2550 696}%
\special{pa 2530 952}%
\special{dt 0.045}%
\special{sh 1}%
\special{pa 2530 952}%
\special{pa 2556 886}%
\special{pa 2534 898}%
\special{pa 2516 884}%
\special{pa 2530 952}%
\special{fp}%
%
\special{pn 8}%
\special{pa 2160 696}%
\special{pa 2290 990}%
\special{dt 0.045}%
\special{sh 1}%
\special{pa 2290 990}%
\special{pa 2282 922}%
\special{pa 2268 942}%
\special{pa 2246 938}%
\special{pa 2290 990}%
\special{fp}%
%
\special{pn 8}%
\special{pa 1910 690}%
\special{pa 2160 690}%
\special{dt 0.045}%
\end{picture}%
\hspace{11truecm}
}

\vspace{0.3truecm}

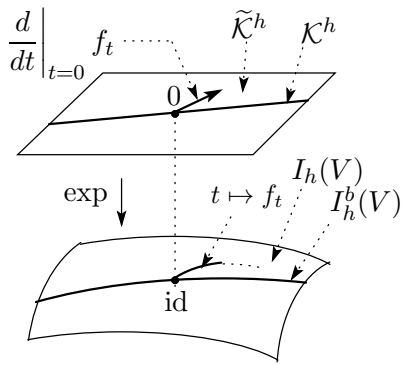
\centerline{{\bf Figure 3.}}

\vspace{1truecm}

\centerline{{\bf References}}

\vspace{0.5truecm}

\small{
\noindent
[Be] M. Berger, Les espaces sym$\acute e$triques non compacts, 
Ann. Sci. $\acute E$c. Norm. Sup$\acute e$r. III. S$\acute e$r. 

{\bf 74} (1959) 85-177.

\noindent
[Br] M. Br$\ddot u$ck, Equifocal famlies in symmetric spaces of compact type, J.reine angew. 
Math. 

{\bf 515} (1999), 73-95.

\noindent
[BCO] J. Berndt, S. Console and C. Olmos, Submanifolds and holonomy, Research 
Notes in 

Mathematics 434, CHAPMAN $\&$ HALL/CRC Press, Boca Raton, London, New York 

Washington, 2003.

%
%

\noindent
[Ch] U. Christ, 
Homogeneity of equifocal submanifolds, J. Differential Geom. {\bf 62} (2002), 
1--15.

\noindent
[Cox] H. S. M. Coxeter, Discrete groups generated by reflections, 
Ann. of Math. (2) {\bf 35} (1934),

588--621.

%

\noindent
[E] H. Ewert, 
Equifocal submanifolds in Riemannian symmetric spaces, Doctoral thesis.

\noindent
[G1] L. Geatti, 
Invariant domains in the complexfication of a noncompact Riemannian 
symmetric 

space, J. Algebra {\bf 251} (2002), 619--685.

\noindent
[G2] L. Geatti, 
Complex extensions of semisimple symmetric spaces, manuscripta math. {\bf 120} 

(2006) 1-25.

\noindent
[GG] L. Geatti and C. Gorodski, 
Polar orthogonal representations of real reductive algebraic grou-

ps, J. Algebra {\bf 320} (2008) 3036-3061.

\noindent
[GH] C. Gorodski and E. Heintze, 
Homogeneous structures and rigidity of isoparametric subman-

ifolds in Hilbert space, J. Fixed Point Theory Appl. {\bf 11} (2012) 93-136.

\noindent
[Ha] J. Hahn, Isotropy representations of semisimple symmetric spaces 
and extrinsically homoge-

neous hypersurfaces, J. Math. Soc. Japan {\bf 40} (1988) 271-288.  

\noindent
[HL1] E. Heintze and X. Liu, 
A splitting theorem for isoparametric submanifolds in Hilbert space, 

J. Differential Geom. {\bf 45} (1997), 319--335.

\noindent
[HL2] E. Heintze and X. Liu, 
Homogeneity of infinite dimensional isoparametric submanifolds, 

Ann. of Math. {\bf 149} (1999), 149-181.

\noindent
[HLO] E. Heintze, X. Liu and C. Olmos, 
Isoparametric submanifolds and a Chevalley type restric-

ction theorem, Integrable systems, geometry, and topology, 151-190, 
AMS/IP Stud. Adv. 

Math. 36, Amer. Math. Soc., Providence, RI, 2006.

\noindent
[HOT] E. Heintze, C. Olmos and G. Thorbergsson, 
Submanifolds with constant prinicipal curva-

tures and normal holonomy groups, Int. J. Math. {\bf 2} (1991), 167-175.

\noindent
[HPTT] E. Heintze, R. S. Palais, C. L. Terng and G. Thorbergsson, 
Hyperpolar actions on sym-

metric spaces, Geometry, topology and physics, 214--245 Conf. Proc. Lecture 
Notes Geom. 

Topology {\bf 4}, Internat. Press, Cambridge, Mass., 1995.

\noindent
[He] S. Helgason, 
Differential geometry, Lie groups and symmetric spaces, Pure Appl. Math. 80, 

Academic Press, New York, 1978.

\noindent
[Hu] M. C. Hughes, 
Complex reflection groups, Comm. Algebra {\bf 18} (1990), 3999--4029.

\noindent
[KN] S. Kobayashi and K. Nomizu, 
Foundations of differential geometry, Interscience Tracts in 

Pure and Applied Mathematics 15, Vol. II, New York, 1969.



\noindent
[K1] N. Koike, 
Submanifold geometries in a symmetric space of non-compact 
type and a pseudo-

Hilbert space, Kyushu J. Math. {\bf 58} (2004), 167--202.

\noindent
[K2] N. Koike, 
Complex equifocal submanifolds and infinite dimensional 
anti-Kaehlerian isopara-

metric submanifolds, Tokyo J. Math. {\bf 28} (2005), 
201--247.

\noindent
[K3] N. Koike, 
Actions of Hermann type and proper complex equifocal submanifolds, Osaka J. 

Math. {\bf 42} (2005) 599-611.

\noindent
[K4] N. Koike, 
A splitting theorem for proper complex equifocal submanifolds, Tohoku Math. J. 

{\bf 58} (2006) 393-417.

%

\noindent
[K5] N. Koike, 
The homogeneous slice theorem for the complete complexification of a 
proper 

complex equifocal submanifold, Tokyo J. Math. {\bf 33} (2010), 1-30.

%

%

%

\noindent
[K6] N. Koike, Collapse of the mean curvature flow for equifocal 
submanifolds, Asian J. Math. 

{\bf 15} (2011), 101-128.

\noindent
[K7] N. Koike, Homogeneity of infinite dimensional anti-Kaehler isoparametric submanifolds, 

Tokyo J. Math. {\bf 37} (2014) 159-178.  

\noindent
[M] I.G. Macdonald, Affine root systems and Dedekind's $\eta$-function, Invent. Math. {\bf 15} (1972), 

91-143.

%

\noindent
[O] C. Olmos, Isoparametric submanifolds and their homogeneous structures, 
J. Differential 

Geom. {\bf 38} (1993), 225-234.

\noindent
[OW] C. Olmos and A. Will, Normal holonomy in Lorentzian space and 
submanifold geometry, 

Indiana Univ. Math. J. {\bf 50} (2001), 1777-1788.

\noindent
[O'N] B. O'Neill, 
Semi-Riemannian Geometry, with applications to relativity, 
Pure Appl. Math. 

103, Academic Press, New York, 1983.

\noindent
[OS] T. Ohshima and J. Sekiguchi, The restricted root system of 
a semisimple symmetric pair, 

Group representations and systems of 
differential equations (Tokyo, 1982), 433--497, Adv. 

Stud. Pure Math. {\bf 4}, North-Holland, Amsterdam, 1984.

\noindent
[Pa] R. S. Palais, 
Morse theory on Hilbert manifolds, Topology {\bf 2} (1963), 299--340.

%

\noindent
[PT2] R. S. Palais and C. L. Terng, Critical point theory and submanifold 
geometry, Lecture Notes 

in Math. {\bf 1353}, Springer-Verlag, Berlin, 1988.

\noindent
[Pe] A. Z. Petrov, Einstein spaces, Pergamon Press, 1969.

\noindent
[R] W. Rossmann, The structures of semisimple symmetric spaces, 
Canad. J. Math. {\bf 31} (1979), 

157-180.

\noindent
[S1] R. Sz$\ddot{{{\rm o}}}$ke, Involutive structures on the 
tangent bundle of symmetric spaces, Math. Ann. {\bf 319} 

(2001), 319--348.

\noindent
[S2] R. Sz$\ddot{{{\rm o}}}$ke, Canonical complex structures associated to 
connections and complexifications of 

Lie groups, Math. Ann. {\bf 329} (2004), 553--591.

\noindent
[Te1] C. L. Terng, 
Isoparametric submanifolds and their Coxeter groups, 
J. Differential Geom. {\bf 21} 

(1985), 79--107.

\noindent
[Te2] C. L. Terng, 
Proper Fredholm submanifolds of Hilbert space, 
J. Differential Geom. {\bf 29} 

(1989), 9--47.

\noindent
[Te3] C. L. Terng, 
Polar actions on Hilbert space, J. Geom. Anal. {\bf 5} (1995), 
129--150.

\noindent
[TT] C. L. Terng and G. Thorbergsson, 
Submanifold geometry in symmetric spaces, J. Differential 

Geom. {\bf 42} (1995), 665--718.

\noindent
[Th] G. Thorbergsson, Isoparametric 
foliations and their buildings, Ann of Math. {\bf 133} (1991), 

429--446.






\end{document}